\documentclass[a4paper,11pt]{article}
\usepackage[width=15cm]{geometry}
\usepackage[british]{babel}
\usepackage[utf8]{inputenc}
\usepackage{amsmath,amssymb,amsfonts,stmaryrd,enumerate,bbm}
\usepackage{mathrsfs}
\usepackage{theorem}

\catcode`\<=\active \def<{
\fontencoding{T1}\selectfont\symbol{60}\fontencoding{\encodingdefault}}
\catcode`\>=\active \def>{
\fontencoding{T1}\selectfont\symbol{62}\fontencoding{\encodingdefault}}
\catcode`\|=\active \def|{
\fontencoding{T1}\selectfont\symbol{124}\fontencoding{\encodingdefault}}
\newcommand{\assign}{:=}
\newcommand{\mathd}{\mathrm{d}}
\newcommand{\nobracket}{}
\newcommand{\nocomma}{}
\newcommand{\noplus}{}
\newcommand{\nosymbol}{}
\newcommand{\tmdate}[1]{\today}
\newcommand{\tmem}[1]{{\em #1\/}}
\newcommand{\tmop}[1]{\ensuremath{\operatorname{#1}}}
\newcommand{\tmscript}[1]{\text{\scriptsize{$#1$}}}
\newcommand{\tmtextit}[1]{{\itshape{#1}}}
\newenvironment{enumeratenumeric}{\begin{enumerate}[1.] }{\end{enumerate}}
\newenvironment{enumerateroman}{\begin{enumerate}[i.] }{\end{enumerate}}

\newtheorem{theorem}{Theorem}[section]

\newtheorem{conjecture}[theorem]{Conjecture}
\newtheorem{corollary}[theorem]{Corollary}

\newtheorem{definition}[theorem]{Definition}

\newtheorem{lemma}[theorem]{Lemma}

\newtheorem{proposition}[theorem]{Proposition}

\newtheorem{remark}[theorem]{Remark}

{\theorembodyfont{\rmfamily}}

\newenvironment{proof}{{\sc Proof.}}{\hfill$\square$}

\newcommand{\CC}{\mathscr C}
\newcommand{\CF}{\mathscr F}

\newcommand{\CQ}{\mathscr Q}

\newcommand{\CN}{\mathscr N}

\newcommand{\RR}{\mathbb R}

\begin{document}

\title{
  Averaging along irregular curves\\
  and regularisation of ODEs.
}

\author{  R. Catellier$^{1}$\footnote{Present address : IRMAR, Centre Henri Lebesgue, Universit\'e de Rennes 1.} and
  M. Gubinelli$^{1,2}$
\\[.1cm] {\small 
(1) CEREMADE  UMR 7534,  Universit\'e Paris--Dauphine \& CNRS, France
}
\\ {\small
(2) Institut Universitaire de France
}
\\ {\small \texttt{\{catellier,gubinelli\}@ceremade.dauphine.fr} }
}

%\date{\today}

\maketitle

\begin{abstract}
  We consider the ordinary differential equation (ODE) $\mathd x_{t} =b
  (t,x_{t} ) \mathd t+ \mathd w_{t}$ where $w$ is a continuous driving
  function and $b$ is a time-dependent vector field which possibly is only a
  distribution in the space variable. We quantify the regularising properties
  of an arbitrary continuous path $w$ on the existence and uniqueness of
  solutions to this equation. In this context we introduce the notion of
  $\rho$-\tmtextit{irregularity} and show that it plays a key role in some
  instances of the regularisation by noise phenomenon. In the particular case
  of a function $w$ sampled according to the law of the fractional Brownian
  motion of Hurst index $H \in (0,1)$, we prove that almost surely the ODE
  admits a solution for all $b$ in the Besov-Hölder space $B^{\alpha
  +1}_{\infty , \infty}$ with $\alpha >-1/2H$. If $\alpha >1-1/2H$ then the
  solution is unique among a natural set of continuous solutions. If $H>1/3$
  and $\alpha >3/2-1/2H$ or if $\alpha >2-1/2H$ then the equation admits a
  unique Lipschitz flow. Note that when $\alpha <0$ the vector field $b$ is
  only a distribution, nonetheless there exists a natural notion of solution
  for which the above results apply.

\end{abstract}

\tableofcontents

\section{Introduction}\label{sec:intro}

In~{\cite{MR2377011}} A. M. Davie showed that the integral equation
\begin{equation}
  \label{eq:ODE-0} x_{t} =x_{0} + \int_{0}^{t} b (s,x_{s} ) \mathd s+w_{t} ,
  \hspace{2em} t \in [0,1],
\end{equation}
with $x,w \in C ([0,1]; \RR^{d} )$ and $b: \RR \times \RR^{d} \to \RR^{d}$
bounded and measurable has a unique continuous solution for almost every path
$w$ sampled from the law of the $d$-dimensional Brownian motion. This result
can be interpreted as a phenomenon of regularisation by noise, in the sense
that it is well known that the same equation without $w$ can show
non-uniqueness.

Regularisation by noise in the case of stochastic differential equations
(SDEs) driven by Brownian motion is nowadays a well understood subject: see
for example Veretennikov, Krylov and Roeckner~{\cite{MR2117951}}, Flandoli,
Gubinelli and Priola~{\cite{FGP}}, Zhang, Flandoli and Da~Prato~{\cite{DF}}.
All these work are essentially based of the use of Itô calculus to highlight
the regularising properties of Brownian paths. Meyer-Brandis and
Proske~{\cite{MR2606880}} use Malliavin calculus to derive similar
conclusions. Davie's contribution~{\cite{MR2377011}} is more subtle in the
sense that it is a result for an ordinary differential equation (ODE) and not
for the related SDE, i.e. the existence and uniqueness of solutions is studied
in the space of continuous paths and not in the more common probabilistic
framework of continuous adapted processes on a given filtered probability
space. This has been clearly pointed out by
Flandoli~{\cite{flandoli_random_2011}} which called these more general
solutions {\tmem{path-by-path}}. In this respect Davie's contribution is
purely analytical and one of the aim of the present work is to
{\tmem{analytically}} caracterize the regularisation effect for general
continuous perturbation $w$ (whether random or not) to the evolution dictated
by an irregular vector field.

Regularisation by ``fast'' or ``dispersive'' motions is an interesting
phenomenon which appears also in some deterministic PDE situations, for
example for Korteweg-de-Vries equation~{\cite{kdv,BIT}} and for fast-rotating
Euler and Navier-Stokes equations~{\cite{BMN}}. In particular the technique of
Young integration we employ in the present work is essentially the same used
in the paper~{\cite{kdv}} to study the periodic Korteweg-de-Vries equation and
take inspiration in the theory of rough paths {\cite{Lyons,controlling,FV}}.

In a recent paper~{\cite{random-modulation-1,random-modulation-2}} \ Chouk
and Gubinelli analyse the regularisation phenomenon in the context of
non-linear dispersive PDEs modulated by an irregular signal. In particular
they considered equations of the form
\begin{equation}
  \label{eq:abs} \frac{\mathd}{\mathd t} \varphi_{t} =A \varphi_{t} 
  \frac{\mathd w_{t}}{\mathd t} +\CN ( \varphi_{t} ) , \hspace{2em} t
  \ge 0
\end{equation}
where $w: \RR_{+} \to \RR$ is an arbitrary continuous function, $A$ is an
unbounded linear operator (like the Schrödinger operator $i \partial^{2}$ or
the Airy operator $\partial^{3}$ acting on periodic or non-periodic functions)
\ and $\CN$ some local polynomial non-linearity with possibly
derivative terms. The unifying theme of this last study and the present one is
the fact that the regularising properties of $w \in C ([0,1]; \RR^{d} )$ are
analysed in terms of the \tmtextit{averaging operator} $T^{w}_{t}$ defined as
\begin{equation}
  T^{w}_{t} f (x) = \int_{0}^{t} f (x+w_{r} ) \mathd r, \hspace{2em} x \in
  \mathbbm{R}^{d} \label{eq:averaging-op}
\end{equation}
for any measurable functions $f: \RR^{d} \to \RR$. Characterising the mapping
properties of $T^{w}$ for various kind of perturbations $w$ seems very
interesting and not straightforward. Mapping properties of $T^{w}$ for
deterministic smooth curves $w$ are, for reasons not related to the
regularisation by noise phenomenon, an interesting subject in analysis: we
have in mind, for example the work of Tao and Wright~{\cite{tao}} on $L^{p}$
improving bounds for averages along curves (we thank F.~Flandoli and
V.~M.~Tortorelli for having pointed us the existence of these results). \

The averaging operator can be seen as the convolution against the
\tmtextit{occupation measure} $L^{w}_{t}$ of the path $w$ defined as
\[ L^{w}_{t} ( \mathd y ) = \int_{0}^{t } \delta_{w_{u}} ( \mathd y ) \mathd u. \]
Indeed, for continuous $b$, the following computation holds
\[ T^{w}_{t} b ( x ) = \int_{0}^{t} b ( x+w_{u} ) \mathd u = \int_{0}^{t}
   \mathd u \int_{\mathbbm{R}^{d}}  b ( x-y ) \delta_{w_{u}} ( \mathd y ) = (
   b \ast L_{t}^{w} ) ( x ) . \]

The basic observation contained in Davie's paper~{\cite{MR2377011}} is that if
$b:\mathbbm{R}^{d} \rightarrow \mathbbm{R}^{d}$ is a given bounded function
then for almost every $d$-dimensional Brownian path $w: [ 0,T ] \rightarrow
\mathbbm{R}^{d}$ and for all $0 \leqslant t \leqslant T$ the function $x
\mapsto T^{w}_{t} b ( x )$ has almost Lipschitz regularity (its modulus of
continuity is of the type $| x | \log^{1/2}_{} ( 1/ | x | )$). Morally this is
a gain of almost $1$ degree of the regularity and one of the key steps to prove
uniqueness of the ODE~{\eqref{eq:ODE-0}} for a bounded measurable drift $b$.

In this paper we analyse the behaviour of the averaging operator $T^{w}$ in
the scale of Hölder-Besov spaces $\CC^{\alpha} = \CC^{\alpha} ( \RR^{d} ,
\RR^{n} ) =B^{\alpha}_{\infty , \infty} ( \RR^{d} , \RR^{n} )$ for arbitrary
regularity $\alpha \in \mathbbm{R}$. We consider a class of perturbations $w$
given by the sample paths of the $d$-dimensional fractional Brownian motion
(fBm) of Hurst index $H \in ( 0,1 )$, that is the unique centered Gaussian
process $( B^{H}_{t} )_{t \geqslant 0}$ with values on $\mathbbm{R}^{d}$ and
covariance function
\[ \mathbbm{E} [ B^{H}_{t} B_{s}^{H} ] =c_{H} ( | t+s |^{2H} - | t |^{2H} - |
   s |^{2H} )Id \]
for all $t,s \geqslant 0$.

As an application of the averaging properties we obtain various existence and
uniqueness results for solutions of the ODE~{\eqref{eq:ODE-0}} and relative
flow properties for distributional vector field $b$.

The choice of fBm has the advantage of being a simple process for which many
other results about existence and uniqueness of associated SDE are
available~{\cite{MR1934157}}. More interestingly, the approach based on Itô
calculus, used in most of the papers on the regularisation effect for Brownian
motion, does not easily extend to the fBm case, nor does the explicit computations
of Davie~{\cite{MR2377011}}. The freedom in the choice of the Hurst parameter
gives us the possibility to explore the effect of different degrees of
irregularity of the perturbation on the regularisation phenomenon and the
quasi-invariance of the law of the fBm will allow us to study the effect of
perturbations on the the averaging properties of the paths.

Returning to the averaging behaviour of fBm paths we obtain the following
result

\begin{theorem}
  \label{th:regularization-fbm}Take $H \in ( 0,1 )$ and $\rho <1/2H$. Then
  there exists $\gamma >1/2$ such that for all $f \in C ( \mathbbm{R}^{d}
  ;\mathbbm{R} )$ there exists a Borel set $\CN_{f, \gamma} \subseteq
  C ( [ 0,1 ] ,\mathbbm{R}^{d} )$ (which depends on $f, \gamma$) of zero
  measure with respect to the law of the $d$-dimensional fractional Brownian
  motion (fBm) of Hurst index $H$ such that for all $w \notin
  \CN_{f, \gamma}$ we have for $\alpha>-\rho$,
  \[ \| T^{w}_{t} f-T^{w}_{s} f \|_{\CC^{\alpha + \rho , \psi}} \lesssim_{w}
     \| f \|_{\CC^{\alpha}} | t-s |^{\gamma} \]
  for all $0 \leqslant s<t \leqslant 1$.
\end{theorem}

In this statement the weighted space $\CC^{\alpha , \psi}$ is a subspace of
the space of local Hölder continuous functions with given grow at infinity
described by the weight $\psi$, and its precise definition is given in Definition \ref{def:weighted_holder} below. The space $\CC^\alpha$ is the usual Besov-Hölder define below in \eqref{eq:besov_holder}.

Letting for a moment aside the time regularity, this result shows that the
averaging against fBm paths gains almost $1/2H$ derivatives in the space
variable. Unfortunately the result stated in
Theorem~\ref{th:regularization-fbm} is not very satisfying since one would
really like to have the almost sure boundedness of $T^{w}_{t} : \CC^{\alpha}
\rightarrow \CC^{\alpha + \rho , \psi}$. The difficulty is, of course, the
fact that the exceptional set $\CN_{f}$ of
Theorem~\ref{th:regularization-fbm} depends itself on the function $f$. Using
the Littewood--Paley decomposition of Besov--Hölder distributions and the
scaling of the fractional Brownian motion, the problem of finding a version of
$T^{w}$ which is almost surely continuous can be related to the following
conjecture:

\begin{conjecture}
  Let $(B^{H}_{t} )_{t \geqslant 0}$ be a $d$-dimensional fBm of Hurst index
  $H \in ( 0,1 )$. Let $K:\mathbbm{R}^{d} \rightarrow \mathbbm{R}$ be a smooth
  function such that
  \[ | K ( x ) | \lesssim ( 1+ | x | )^{-N} , \hspace{2em}
     \int_{\mathbbm{R}^{d}} K ( x ) \mathd x=0, \]
  where $N>d$ can be chosen arbitrarily large. Then
  \[ \mathbbm{E} ( \| T^{B^{H}}_{0,t} K \|_{L^{1} ( \mathbbm{R}^{d} )}^{p} )
     =\mathbbm{E} \left[ \left( \int_{\mathbbm{R}^{d}} \Big\lvert \int_{0}^{t} K (
     x+B^{H}_{s} ) \mathd s \Big\lvert \mathd x \right)^{p} \right] \lesssim
     t^{p/2} \]
  as $t \rightarrow + \infty$.
\end{conjecture}

If the function $K$ has a bounded support the estimation is true as an easy
consequence of our results, however currently we are unable to prove or
disprove this conjecture.

On the positive side if we replace $\CC^{\alpha}$ by the Fourier--Lebesgue
spaces $\mathcal{F}L^{\alpha ,p}$ defined as
\[ \mathcal{F}L^{\alpha ,p} ( \RR^{d} ) = \{f \in \mathcal{S}' ( \RR^{d}
   ):N_{\alpha ,p} (f)^{p} = \int_{\RR^{d}} | \hat{f} ( \xi )|^{p} (1+| \xi
   |)^{\alpha p} \mathd \xi < \infty\}, \]
with $\mathcal{F}L^{\alpha} =\mathcal{F}L^{\alpha ,1}$, then it is easy to
see that for $0 \leqslant \gamma \leqslant 1$ and $\rho \in \mathbbm{R}$:
\[ \| T^{w}_{t} -T^{w}_{s} \|_{\mathcal{L}(\mathcal{F}L^{\alpha} ;
   \mathcal{F}L^{\alpha + \rho})} = \sup_{f \in \mathcal{F}L^{\alpha}}
   \frac{N_{\alpha + \rho} ( T^{w}_{t} f-T^{w}_{s} f )}{N_{\alpha} ( f )}
   \lesssim \| \Phi^{w} \|_{\mathcal{W}^{\rho , \gamma}_{T}} | t-s |^{\gamma},
\]
where $\Phi^{w}_{t} (a) = \int_{0}^{t} e^{i \langle a,w_{r} \rangle} \mathd
r=e^{-i \langle a,x \rangle} T^{w}_{t} (e^{i \langle a, \cdot \rangle} ) ( x
)$ and where we introduced the norm
\[ \| \Phi^{w} \|_{\mathcal{W}^{\rho , \gamma}_{T}} = \sup_{a \in \RR^{d}}  
   \sup_{0 \leqslant s<t \leqslant T} (1+ |a|)^{\rho} \frac{|
   \Phi^{w}_{t} (a)- \Phi^{w}_{s} (a) |}{|s-t|^{\gamma}} . \]
This observation reduces the question of the boundedness of $T^{w}$ to that of
the decay of the Fourier transform $a \mapsto \Phi^{w}_{t} (a)$ of the
{\tmem{occupation measure}} of $w$~(for generalities about occupation measures
and densities for deterministic and random functions see for example the
review of Geman and Horowitz~{\cite{geman_occupation_1980}}). This suggests to
introduce the following novel notion of "irregularity" of the perturbation $w$
:

\begin{definition}
  \label{def:irregularity}Let $\rho >0$ and $\gamma >0$. We say that a
  function $w \in C ([0,T]; \RR^{d} )$ is $( \rho , \gamma )$-irregular if
  \[ \| \Phi^{w} \|_{\mathcal{W}^{\rho , \gamma}_{T}} <+ \infty . \]
  Moreover we say that $w$ is $\rho$-irregular if there exists $\gamma >1/2$
  such that $w$ is $( \rho , \gamma )$-irregular.
\end{definition}

The time regularity of this Fourier transform, measured by the Hölder
exponent $\gamma$, will also be crucial in our analysis. The notion of
$\rho$-irregularity is also relevant to the boundedness of $T^{w}$ in other
functional spaces, for example we easily see that for all $\alpha \in
\mathbbm{R}$:
\[ \|T^{w}_{t} f-T^{w}_{s} f\|_{H^{\alpha + \rho} ( \RR^{d} )} \leqslant \|
   \Phi^{w} \|_{\mathcal{W}^{\rho , \gamma}_{T}} |t-s|^{\gamma}
   \|f\|_{H^{\alpha} ( \RR^{d} )}, \]
where $H^{\alpha} ( \mathbbm{R}^{d} ) = \CF L^{\alpha ,2}$ are the usual
Sobolev spaces on $\mathbbm{R}^{d}$ and in general similar inequalities holds
in Fourier--Lebesgue spaces $\CF L^{\alpha ,p}$ of arbitrary integrability $p
\in [ 1,+ \infty ]$. However the notion of $\rho$-irregularity does not seem enough
 to control the boundedness of the averaging operator in Besov spaces.

The limiting value $1/2$ for $\gamma$ does not seem to have any special meaning, as
far as the occupation measure is concerned, however if $\gamma >1/2$ we are
able to develop a quite simple integration theory for the averaging operator
using Young integral techniques and quite surprisingly it turns out that this
is sufficient for the purpose of this paper. Indeed a proof similar to that of
Theorem~\ref{th:regularization-fbm} gives the existence of (plenty of)
perturbations $w$ which are $\rho$-irregular :

\begin{theorem}
  \label{th:irr-for-fbm}Let $(B^{H}_{t} )_{t \geqslant 0}$ be a fractional
  Brownian motion of Hurst index $H \in ( 0,1 )$ then for any $\rho <1/2H$
  there exist $\gamma >1/2$ so that with probability one the sample paths of
  $B^{H}$ are $( \rho , \gamma )$-irregular.
\end{theorem}

In particular there exist continuous paths which are $\rho$-irregular for
arbitrarily large $\rho$ and thus paths which deliver an arbitrary degree of
regularisation. Using well known properties of support of the law of the
fractional Brownian motion it is also possible to show that there exists $\rho$-irregular
trajectories which are arbitrarily close in the supremum norm to any smooth
path.

As a direct corollary of Theorem~\ref{th:irr-for-fbm} we have the boundedness
of $T^{w}$ in the Fourier--Lebesgue spaces $\CF L^{\alpha}$:

\begin{corollary}
  Let $H \in ( 0,1 )$ and $\rho <1/2H$. Then almost surely with respect to the law of the
  fBm of Hurst index $H$ we have that for all $0 \leqslant s \leqslant t
  \leqslant T$ \ the averaging operator $T^{w}$ is bounded from $\CF
  L^{\alpha}$ to $\CF L^{\alpha + \rho}$ and satisfy
  \[ \| T^{w}_{t} -T^{w}_{s} \|_{\mathcal{L} \left( \CF L^{\alpha} ; \CF
     L^{\alpha + \rho} \right)} \leqslant C_{w, \gamma , \rho} | t-s
     |^{\gamma} \]
  for some constant $C_{w, \gamma , \rho}$ which depends only on $w, \gamma ,
  \rho$. This means that $$T^{w}_{} \in \CC^{\gamma} \left( [ 0,T ] ;\mathcal{L}
  \left( \CF L^{\alpha} , \CF L^{\alpha + \rho} \right) \right).$$
\end{corollary}

One of the contributions of our work is the observation that the regularity of
the occupation measure of $w$ seems to play a major role in the understanding
of the regularising properties of $w$ in a non-linear context and it would be
desirable to understand more deeply the link of the notion of
$\rho$-irregularity with the pathwise properties of $w$, for example linking
them to the notion of true roughness appearing in the literature on densities
for differential equations driven by rough paths~{\cite{true-roughness}}.

It would also be interesting to study more deeply the notion of irregularity
for ``generic'' continuous paths (for example in the class of Hölder
continuous paths). Indeed, set aside the classic contribution of Geman and
Horowitz~{\cite{geman_occupation_1980}} mentioned above, the authors are not
aware of any systematic study of occupation measures of random processes from
the point of view of their action on spaces of functions or distributions,
topic which seems central to our analysis.

An open problem is, for example, understanding what happens if we replace $w$
with a regularised version $w^{\varepsilon}$ or with a perturbed version. In
this respect we conjecture that if $w$ is $( \rho , \gamma )$-irregular then
for any smooth function $\varphi \in C ( [ 0,1 ] ;\mathbbm{R}^{d} )$ the
perturbed path $w^{\varphi} =w+ \varphi$ is still $( \rho , \gamma
)$-irregular. In relation to this last problem we have obtained the following
general result:

\begin{theorem}
  \label{th:perturbed-young}Let $\rho \in \mathbbm{R}$ and $\varphi \in
  \CC^{\beta} ( [ 0,T ] ;\mathbbm{R}^{d} )$ with $1/2 \leqslant \beta <1$.
  Then if $w$ is $\rho$-irregular the path $w^{\varphi} =w+ \varphi$ is $(
  \rho -1/2 \beta )$-irregular. Moreover for $\gamma >1/2$ we have
  \[ \| T^{w+ \varphi}_{} f \|_{\CC^{\gamma} \left( [ 0,T ] ; \CC^{\alpha + \rho
     -1/2 \beta} \right)} \lesssim_{T, \beta , \gamma} \| T^{w} f
     \|_{\CC^{\gamma , \psi} \left( [ 0,T ] ; \CC^{\alpha + \rho} \right)} \|
     \varphi \|_{\CC^{\beta}}. \]
In particular if $T^{w} \in \CC^{\gamma} \left( [ 0,T ] ;\mathcal{L}
  \left( \CC^{\alpha} ; \CC^{\alpha + \rho , \psi} \right) \right)$ then
  $$T^{w+ \varphi} \in \CC^{\gamma} \left( [ 0,T ] ;\mathcal{L} \left(
  \CC^{\alpha} ; \CC^{\alpha + \rho -1/2 \beta} \right) \right).$$
\end{theorem}

In particular the irregularity property is preserved at the price of a loss at
least $1/2$ in regularity (which happens when $\beta$ is close to $1$).

If $w$ is sampled according to the law of a fBm and if \ the perturbation
$\varphi$ is adapted to the natural filtration of $w$ then it is possible to
exploit the quasi-invariance of the fBm measure with respect to adapted shifts to prove
the irregularity of the perturbed path without any loss on the irregularity
exponent:

\begin{theorem}
  \label{th:perturbed-girsanov}Let $B^{H}$ be a fBm of Hurst index $H \in (
  0,1 )$ and let $\Phi : [ 0,T ] \rightarrow \mathbbm{R}^{d}$ be an Hölder
  continuous process which is adapted to the natural filtration of $B^H$. Then,
  for all $\rho <1/2H$ \ almost surely the process $B^H+ \Phi$ is
  $\rho$-irregular and for any $f \in \CC^{\alpha}$
  \[ \| T^{B^H+ \Phi}_{} f \|_{\CC^{\gamma} \left( [ 0,T ] ; \CC^{\alpha + \rho ,
     \psi} \right)} <+ \infty \quad \text{almost surely.}\]

\end{theorem}

The disadvantage of this result is that the exceptional set where the
irregularity property fails depends a priori on $\Phi$ and this poses problems
in applications to pathwise results valid for a large class of perturbations
(for example smooth and adapted $\Phi$).

One of our aims is to apply these results on the averaging properties of
paths $w$ and of its perturbations to the study of existence and uniqueness of
solutions to the ODE~{\eqref{eq:ODE-0}} for distributional $b$. Two main
situations will be considered:
\begin{enumeratenumeric}
  \item $b \in \CC^{\alpha}$ (or $b \in \CF L^{\alpha}$) for some $\alpha >0$.
  In this case $b$ will be a bounded continuous function and the
  ODE~{\eqref{eq:ODE-0}} has a natural meaning and allows for continuous
  solution, we will then consider the related uniqueness problem and the
  existence of a Lipshitz flow.
  
  \item $b \in \CC^{\alpha}$ (or $b \in \CF L^{\alpha}$) for some $\alpha <0$.
  In this case even the appropriate meaning to give to the
  ODE~{\eqref{eq:ODE-0}} is not clear and we will investigate this problem and
  the related well-posedness and continuity issues.
\end{enumeratenumeric}
In the case $\alpha \geqslant 0$ we have the following results:

\begin{theorem}
  \label{th:ode-flow-0}Let $b \in \CC ( \mathbbm{R}^{d} )$ and assume that
  $\|T^{w} b\|_{\CC^{\gamma} ( [ 0,T ] ; \CC^{3/2, \psi} )} <+ \infty$. Then for
  any $x_{0} \in \mathbbm{R}^{d}$ there exists a unique continuous solution $x
  \in C ([0,1];\mathbbm{R}^{d} )$ of the ODE~{\eqref{eq:ODE-0}} and the flow
  map $x_{0} \mapsto x_{t}$ of the equation is locally Lipshitz continuous in
  space uniformly in $t \in [0,1]$.
\end{theorem}

\begin{theorem}
  \label{th:ode-better-0}Let $b \in \CC^{\alpha}$ and assume that $\alpha
  >1-1/2H$. Then for any $x_{0} \in \mathbbm{R}^{d}$ there exists a measurable
  set of perturbations $\CN_{b,x_{0}} \subseteq C ( [ 0,1 ]
  ;\mathbbm{R}^{d} )$ which is of zero measure with respect to the law of the
  fBm with index $H \in ( 0,1 )$ and such that, for all $w \notin
  \CN_{b,x_{0}}$ there exists a unique continuous solution $x \in C (
  [ 0,1 ] ;\mathbbm{R}^{d} )$ of the ODE~{\eqref{eq:ODE-0}}.
\end{theorem}

As we already remarked, in the case where $b \in \CC^{\alpha}$ for $\alpha
<0$, the ODE~{\eqref{eq:ODE-0}} is not well defined since in general the
evaluation of the distribution $b$ along a continuous curve is not possible.
However if we take into account a particular class of continuous paths we can
show that this coupling has a meaning. A suitable class of continuous
functions is given by a space of paths which are perturbations of $w$:

\begin{definition}
  Let $\gamma\in(0,1)$. The space $\CQ^{w}_{\gamma}$ of $(w, \gamma )$-controlled paths is
  the space
  \[ \CQ^{w}_{\gamma} = \{x \in C ( [ 0,1 ] ;\mathbbm{R}^{d} )  
     \text{  : } x-w \in \CC^{\gamma} ( [ 0,1 ] ;\mathbbm{R}^{d} ) \} . \]
\end{definition}

Then for \tmtextit{controlled paths} we can prove the following result.

\begin{theorem}
  \label{th:def-of-int-b}Let $b \in \mathcal{S}' ( \mathbbm{R}^{d}
  ;\mathbbm{R}^{d} )$, $\gamma>\frac12$  and assume that $\| T^{w} b \|_{\CC^{\gamma} \left( [
  0,T ] ; \CC^{0, \psi} \right)} <+ \infty$. Let $\rho \in \mathcal{S} (
  \mathbbm{R}^{d} )$ be a positive function with $\rho ( 0 ) =1$ and let
  $\rho_{\varepsilon} ( x ) = \varepsilon^{-d} \rho ( x/ \varepsilon )$. Then,
  for all $x \in \CQ^{w}_{\gamma}$,
  \begin{equation}
    \lim_{\varepsilon \rightarrow 0} \int_{0}^{t} ( \rho_{\varepsilon} \ast b
    ) ( x_{s} ) \mathd s=: \int_{0}^{t} b ( x_{s} ) \mathd s
    \label{eq:limiting-procedure}
  \end{equation}
  exists uniformly in $t \in [ 0,T ]$, is independent of $\rho$ and extends
  the usual definition of the right hand side for continuous $b$. Moreover the function $t
  \mapsto \int_{0}^{t} b ( x_{s} ) \mathd s$ is Hölder continuous of exponent
  $\gamma$. 
\end{theorem}

Theorem~\ref{th:def-of-int-b} allows to give a natural meaning to
$\int_{0}^{t} b ( x_{s} ) \mathd s$ for all $x \in \CQ^{w}_{\gamma}$
and from this we can say that $x \in \CQ^{w}_{\gamma}$ is a solutions
of the ODE~{\eqref{eq:ODE-0}} if
\[ x_{t} -w_{t} = \int_{0}^{t} b ( x_{s} ) \mathd s \]
for all $t \in [ 0,1 ]$. That is the ODE has a meaning not in the space of all
continuous functions, as it was when $b$ is a function, but in the more
restricted space of functions which can be seen as ``not too irregular''
additive modifications of $w$. In this context we have natural generalisations
of the Theorems~\ref{th:ode-flow-0} and~\ref{th:ode-better-0} provided we
restrict the space of allowed functions to $\CQ^{w}_{\gamma}$:

\begin{theorem}
  \label{th:ode-flow}Assume that $\|T^{w} b\|_{\CC^{\gamma} ( [ 0,T ] ;
  \CC^{3/2, \psi} )} <+ \infty$. Then for any $x_{0} \in \mathbbm{R}^{d}$
  there exists a unique continuous solution $x \in \CQ^{w}_{\gamma}$
  of the ODE~{\eqref{eq:ODE-0}} and the flow map $x_{0} \mapsto x_{t}$ of the
  equation is Lipshitz continuous uniformly in $t \in [0,1]$.
\end{theorem}

\begin{theorem}
  \label{th:ode-better}Let $b \in \CC^{\alpha +1}$ and assume that $\alpha
  >-1/2H$. Then for any $x_{0} \in \mathbbm{R}^{d}$ there exists a measurable
  set of perturbations $\CN_{b,x_{0}} \subseteq C ( [ 0,1 ]
  ;\mathbbm{R}^{d} )$ which is of zero measure with respect to the law of the
  fBm with index $H \in ( 0,1 )$ and such that, for all $w \notin
  \CN_{b,x_{0}}$ there exists a unique continuous solution $x \in
  \CQ^{w}_{\gamma}$ of the ODE~{\eqref{eq:ODE-0}}.
\end{theorem}

Note that Theorem~\ref{th:ode-flow} \ and~\ref{th:ode-better} are applicable
also when $\alpha \geqslant 0$. In this case existence of solutions is simply
a result of a compactness argument in $C ( [ 0,1 ] ;\mathbbm{R}^{d} )$ and
given a continuous solution it belongs necessarily to
$\CQ^{w}_{\gamma}$ so, in this case, Theorem~\ref{th:ode-flow} \
and~\ref{th:ode-better} are natural generalisations of
Theorems~\ref{th:ode-flow-0} and~\ref{th:ode-better-0}.

When $w$ is sampled according to the law of the fBm with Hurst parameter $H$
Theorem~\ref{th:ode-flow} give the following corollary

\begin{theorem}
  \label{th:flow-H}Fix $H \in ( 0,1 )$ and assume that $b \in \CC^{\alpha
  +3/2}$ for some $\alpha >-1/2H$. Then there exists a measurable set of
  perturbations $\CN_{b} \subseteq C ( [ 0,1 ] ;\mathbbm{R}^{d} )$
  which is of zero measure with respect to the law of the fBm with index $H
  \in ( 0,1 )$ and such that, for all $w \notin \CN_{b}$ and for
  all $x_{0} \in \mathbbm{R}^{d}$ there exists a unique continuous solution $x
  \in \CQ^{w}_{\gamma}$ of the ODE~{\eqref{eq:ODE-0}} and the
  corresponding flow map $\Phi_{t} :x_{0} \mapsto x_{t}$ is globally Lipshitz.
  Moreover the exceptional set $\CN_{b}$ can be chosen to be the same
  for all $b \in \CF L^{\alpha}$. 
\end{theorem}

An interesting consequence of Theorem~\ref{th:flow-H} is the fact that if one
consider the ODE~{\eqref{eq:ODE-0}} as a strong SDE (that is an equation for
stochastic processes adapted to the filtration generated by the process $w$)
and if $w$ has the law of the fBm of index $H$ then we can allow general
random $b \in \CF L^{\alpha}$ and still retain uniqueness under the regularity
conditions of the theorem. This was one of our main motivation to introduce
the scale of Fourier-Lebesgue regularities $\left( \CF L^{\alpha}
\right)_{\alpha}$. Similar results for the Besov scale $\left( \CC^{\alpha}
\right)_{\alpha}$ are not known since we are not able to prove the
corresponding mapping properties for the averaging operator $T^{w}$. Note that
even in the case of the Brownian motion this was an open problem~{\cite{FGP}}
since the standard approach using stochastic calculus cannot be applied in
this case. Allowing random $b$ could open the way to the study of a general
class of stochastic transport equations where the drift itself depends on the
solution.

The key to obtain these results (the existence part when $\alpha <0$ and the
uniqueness part for $\alpha \geqslant 0$ or $\alpha <0$) lies in the fact that
in all cases the ODE~{\eqref{eq:ODE-0}} is equivalent to an equation of Young
\ type (YE) of the form
\begin{equation}
  \theta_{t} = \theta_{0} + \int_{0}^{t} X_{\mathd s} ( \theta_{s} ),
  \label{eq:young-ode-theta}
\end{equation}
where here $X: [ 0,1 ] \times \mathbbm{R}^{d} \rightarrow \mathbbm{R}^{d}$
plays the role of a time-varying, integrated, vector field and $\theta_{t}
=x_{t} -w_{t}$ is the perturbation which by the hypothesis $x \in
\CQ^{w}_{\gamma}$ belongs to $\CC^{\gamma} ( [ 0,1 ] ;\mathbbm{R}^{d}
)$. The integral operation featuring in~{\eqref{eq:young-ode-theta}} has to be
understood as a natural non-linear generalisation of the Young
intergal~{\cite{Young}} defined as limit of Riemman sums:
\[ \int_{0}^{t} X_{\mathd s} ( \theta_{s} ) = \lim_{| \Pi | \rightarrow 0}
   \sum X_{t_{i} ,t_{i+1}} ( \theta_{t_{i}} ), \]
where $X_{s,t} ( x ) =X_{t} ( x ) -X_{s} ( x )$. In the case of the
ODE~{\eqref{eq:ODE-0}} the integrated vector field $X$ corresponds to the
average of the original vector field $b$ given by $X_{t} ( x ) =T_{0,t}^{w} b
( x )$ for all $t \in [ 0,1 ]$ and $x \in \mathbbm{R}^{d}$. Young differential
equations of the type~{\eqref{eq:young-ode-theta}} are used also
in~{\cite{random-modulation-1,random-modulation-2}} to study the regularisation phenomenon for some
non-linear dispersive equations. The theory of such equations is very similar
to the theory for standard Young-type equation but for the sake of the reader
we rederive here the main results in our slightly non standard setting.

This paper is then divided naturally into two parts: in the first we study the
non-linear Young integral and the YE~{\eqref{eq:young-ode-theta}} and derive
the results announced above about existence and uniqueness for the
ODE~{\eqref{eq:ODE-0}}. In the second we analyse the averaging properties of
fBm sample paths and apply the results to the study of the regularisation
phenomenon for eq.~{\eqref{eq:ODE-0}} driven by fBm paths.

\subsection{Notations}

Several function spaces are involved in the rest of the article. In this
section we define those spaces, and specify some notations.
Let $\psi , \varphi \in \mathcal{D}$ be a nonnegative radial functions such
that
\begin{enumerate}
  \item The support of $\psi$ is contained in a ball and the support of
  $\varphi$ is contained in an annulus;
  
  \item $\psi ( \xi ) + \sum_{j \ge 0} \varphi (2^{-j} \xi ) =1$ for all $\xi
  \in \RR^{d}$;
  
  \item $\tmop{supp} ( \psi ) \cap \tmop{supp} ( \varphi (2^{-j} . )) =
  \emptyset$ for $i \ge 1$ and if $|i-j| >1$, then $\tmop{supp} ( \varphi
  (2^{-i}. )) \cap \tmop{supp} ( \varphi (2^{-j} .)) = \emptyset$.
\end{enumerate}
For the existence of $\psi$ and $\varphi$ see {\cite{bahouri_fourier_2011}}.
The Littlewood-Paley blocks are now defined as
\[ \begin{array}{crlc}
     \Delta_{-1} u= \CF^{-1} ( \psi \CF u) & \text{and for $j\ge 0$}
 & \Delta_{j} u= \CF ( \varphi (2^{-j} .) \CF u) .
   \end{array} \]
The $\Delta_{j} u$ are smooth function with Fourier transform with compact
support. We define the Hölder-Besov space $\CC^{\alpha}$ by
\begin{equation}\label{eq:besov_holder}
 \CC^{\alpha} ( \RR^{d} ; \RR^{n} ) =B^{\alpha}_{\infty , \infty} ( \RR^{d}
   ; \RR^{n} ) = \left\{ u \in \mathcal{S}' ( \RR^{d} , \RR^{n}
   ):\|u\|_{\alpha} =\| ( 2^{j \alpha} \| \Delta_{j} u\|_{\infty} )_{j}
   \|_{\infty} < \infty \right\} .
   \end{equation}
While the norm $\|.\|_{\alpha}$ depends on the choice of $\psi$ and $\varphi$,
the space $\CC^{\alpha}$ does not and each choice of $\psi , \varphi$
correspond to an equivalent semi-norm on $\CC^{\alpha}$. If $\alpha \in
\RR_{+} -\mathbbm{N}$, then the space $\CC^{\alpha}$ is the space of $[ \alpha
]$ times differentiable functions , whose partial derivatives up to order $[
\alpha ]$ are bounded, and whose partial derivatives of order $[ \alpha ]$
are ($\alpha - [ \alpha ]$)-Hölder continuous. Note that we have the
following continuous embedding, for $\alpha' \leqslant \alpha$ then
$\CC^{\alpha} \subset \CC^{\alpha'}$ and $\|u\|_{\alpha'} \lesssim
\|u\|_{\alpha}$. When $f \in C \left( [ 0,T ] , \CC^{\alpha} \right)^{}$, we
denote abusively $\| f \|_{\alpha} = \sup_{t \in I} \|u(t,.)\|_{\alpha}$. When
$\alpha >0$, the space $\CC^{\alpha} =B^{\alpha}_{\infty , \infty}$ is the
space of bounded Hölder continuous functions, indeed, for $m \in
\mathbbm{N}\backslash \{ 0 \}$ and $m-1 \leqslant \alpha <m$, when we define
$\llbracket f \rrbracket_{\nu} = \sup_{x \neq y} | f ( x ) -f ( y ) | / | x-y
|^{\nu}$ for $\nu \in ( 0,1 ]$ and
\[ \CC^{\alpha} = \{ f:\mathbbm{R}^{d} \rightarrow \mathbbm{R}^{d} : \| f
   \|_{\infty ,\mathbbm{R}^{d}} + \llbracket D^{m-1} f \rrbracket_{m- \nu} <+
   \infty \} . \]
Furthermore $\| . \|_{\alpha}$ and $\| f \|_{\infty} + \llbracket f
\rrbracket_{m- \alpha}$ are equivalent norms. We will equally use either
one or the other. We will also need some localised Hölder spaces described as
follows:

\begin{definition}
  \label{definition:nu_weight}Let $\nu \in [ 0,1 )$. A weight is a continuous
  non-decreasing function $\psi :\mathbbm{R}_{+} \rightarrow \mathbbm{R}_{+}$
  such that for $c>0$, there exists a constant $C_{c, \psi} >0$ such that
  \[ \psi ( c x ) \leqslant C_{c, \psi} \psi ( x ) . \]
  A $\nu$-weight is a weight such that
  \[ x^{- ( 1- \nu )} \psi ( x ) \xrightarrow[x \rightarrow + \infty]{} 0. \]
\end{definition}

Hence, in that setting we define some weighted Hölder spaces as

\begin{definition}\label{def:weighted_holder}
  Let $\psi$ be a weight, $\nu \in ( 0,1 ]$ and $V$ and $W$ be two Banach spaces.
  The $\psi$-weighted Hölder space of index $\nu$ is the space $\CC^{\nu ,
  \psi} ( V;W )$ defined by
  \[ \CC^{\nu , \psi} ( V;W ) = \left\{ f: \CC^{\nu}_{\tmop{loc}} : \llbracket
     f \rrbracket_{\nu , \psi} = \sup_{x \neq y \in V} \frac{| f ( x ) -f ( y
     ) |_{W}}{| x-y |_{V} \psi ( | x |_{V} + | y |_{V} )} <+ \infty \right\}.
  \]
  When $n\in \mathbbm{N}$ we say that a continuously $n$-times  (Fréchet) differentiable function $f\in C^n(V;W)$ is in the $\psi$-weighted Hölder space of order $n+\nu$ if  $D^n f \in \CC^{\nu,\psi}(V; \mathcal{L}^n(W;W))$, where $\mathcal{L}^n(W;W)$ denote the space of $n$-linear continuous applications from $W$ into itself.
\end{definition}

To simplify the notation, we introduce also the following spaces related to
time dependent nonlinear mappings between Banach spaces $V$ and $W$.

\begin{definition}
  Let $0< \gamma , \nu \leqslant 1$ and $\psi$ be a weight. Let $I= [0,T]$ and
  $V$ and $W$ be two Banach spaces. For all $n \in \mathbbm{N}$ and any $G:I \times V
  \rightarrow W$ we define
  \[ \llbracket G \rrbracket_{\gamma , \nu , \psi} = \sup_{s \neq t}  
     \sup_{x \neq y} 
     \frac
     		{| (G_{t} ( x ) - G_s(x)) -(G_{t} ( y ) - G_s(y))|}
	     	{| t-s |^{\gamma}
     | x-y |^{\nu} \psi ( | x | + | y | )}, \]
 
  \[ \| G \|_{\gamma ,n+ \nu , \psi} = \llbracket D^{n} G \rrbracket_{\gamma ,
     \nu , \psi} + \sum_{k=0}^{n} \sup_{s \neq t} \frac{| D^{k} G_{t} ( 0
     ) -D^k G_s(0)|}{| t-s |^{\gamma}} \]
     and
  \[ \CC^{\gamma ,n+ \nu , \psi} (I,V,W) = \left\{ G \in L^{\infty} \left( I;
     \CC^{\nu , \psi} ( V;W ) \right) :\|G\|_{\gamma ,n+ \nu , \psi} < \infty
     \right\} . \]
  When $V=W$ we write $\CC^{\gamma , \nu , \psi} (I,V,V) = \CC^{\gamma , \nu ,
  \psi} (I,V)$. Furthermore, when it is not ambiguous we only use $\CC^{\gamma
  , \nu , \psi}$. When $\psi =1$ and there is no ambiguity, we only write
  $\CC^{\gamma , \nu}$.
\end{definition}

As stated in the introduction, in order to have estimates for the averaging
operator $T^{w}$ which will not depend on the functions $f$, we introduce the
following Fourier--Lebesgue spaces

\begin{definition}
  Let $\alpha \in \RR$,
  \[ N_{\alpha ,p} (f) = \left(\int_{\RR^{d}} | \hat{f} ( \xi )|^{p}  ( 1+|
     \xi | )^{p \alpha} \mathd \xi\right)^{1/p} \]
  and $\CF L^{\alpha ,p}  ( \RR^{d} ) = \{f \in \mathcal{S}' ( \RR^{d}
  ):N_{\alpha} (f)< \infty\}$. Then $N_{\alpha ,p}$ is a norm on $\CF
  L^{\alpha ,p}  ( \RR^{d} )$. When $p=1$ we only write $\CF
  L^{\alpha}$ and $N_{\alpha}$. 
\end{definition}

When $\alpha \geqslant 0$ and $f \in \CF L^{\alpha}$ implies that $\hat{f}$ is
in $L_{1}$ and $f$ is bounded continuous function. Furthermore if $\alpha \ge
1$, $f \in \CF L^{\alpha}$ is globally Lipschitz continuous in the second
variable. Furthermore for $\alpha \in (0,1)$, $f \in \CF L^{\alpha}$ is
globally Hölder continuous in the second variable. Note that when $\alpha <0$
the vector fields are only distributions.

\begin{remark}
  An easy computation gives $\CF L^{\alpha} \subset \CC^{\alpha}$ \ for all
  $\alpha \in \mathbbm{R}$, and for $\alpha >0$ and $\psi$ a weight,
  $\CC^{\alpha} \subset \CC^{\alpha , \psi}$.
\end{remark}

It is natural to make some approximations in $\CC^{\alpha}$ and in $\CF
L^{\alpha}$. Although the quantity $\sum_{i} \Delta_{i} f$ does not converge in
$\CC^{\alpha}$, it converges in all $\CC^{\alpha'}$ with $\alpha' < \alpha$,
which gives the following lemma :

\begin{lemma}
  Let $\alpha\in\RR$ and $u\in\CC^\alpha(\RR^d)$. The sequence $( \pi_{\leqslant N} u )_{N \ge -1} = ( \sum_{j \leqslant N}
  \Delta_{j} u)_{N}$ converges to $u$ in $\CC^{\alpha'}$ for all $\alpha' <
  \alpha$. Furthermore, for all $\alpha$, $\pi_{\leqslant N} f
  \xrightarrow{\CF L^{\alpha}} f$ for $f \in \CF L^{\alpha}$.
\end{lemma}

Finally if $G: [ 0,T ] \times \mathbbm{R}^{d} \rightarrow \mathbbm{R}^{d}$ we
write $G_{s,t} ( x ) =G_{t} ( x ) -G_{s} ( x )$.

\section{The non-linear Young integral and Young-type
equations}\label{sec:young-ode}

As already said, we intend to study the ODE~{\eqref{eq:ODE-0}} where
$w:\mathbbm{R}_{+} \to \mathbbm{R}^{d}$ is a continuous function (with $w_{0}
=0$) and $b:\mathbbm{R} \times \mathbbm{R}^{d} \to \mathbbm{R}^{d}$ is a
(time-dependent, distributional) vector field. We think $w$ as a very rough
function whose oscillations dominate in small time scales the effects of the
integrated vector field $b$. In this situation the function $x$ behaves at
small scales very much like $w$ and the effects of $b$ are seen only via a
average over these fast oscillations. All this will cooks up some
regularisation effect which will allow to prove existence and uniqueness even
when the vector field $b$ does not enjoys sufficient space regularity.

To highlight the effect of the translations induced by $w$ on the flow of $b$
let us introduce the change of variables $\theta_{t} =x_{t} -w_{t}$ so that
the above equation now reads:
\[ \theta_{t} = \theta_{0} + \int_{0}^{t} b (w_{s} + \theta_{s} ) \mathd s. \]
If we believe that $w$ oscillate faster than $\theta$ then it seems reasonable
to approximate the integral in the right hand side by a sum over a partition $t_{0} =0,
\ldots ,t_{n} =t$ of $[0,t]$ where we have fixed the $\theta$ parameter at the
initial time of each segment:
\begin{equation}
  \label{eq:approx} \int_{0}^{t} b (s,w_{s} + \theta_{s} ) \mathd s \simeq
  \sum_{i=0}^{n-1} \int_{t_{i}}^{t_{i+1}} b (s,w_{s} + \theta_{t_{i}} ) \mathd
  s= \sum_{i=0}^{n-1} (T^{w}_{t_{i} ,t_{i+1}} b) ( \theta_{t_{i}} ) .
\end{equation}
where $T^{w}_{s,t} =T^{w}_{t} -T^{w}_{s}$.

Under appropriate conditions the expression on the right hand side of
eq.{\eqref{eq:approx}} will have a well defined limit as the size of the
partition goes to zero and it defines a kind of integral which we naturally
denote by
\[ \int_{0}^{t} (T^{w}_{\mathd s} b) ( \theta_{s} ) = \lim   \sum_{i=0}^{n-1}
   (T^{w}_{t_{i} ,t_{i+1}} b) ( \theta_{t_{i}} ) . \]
and will enable us to set up an alternative formulation of the above ODE as an
integral equation involving the time-dependent integrated vector field $G_{t}
=T^{w}_{t} b$ which is an averaged version of $b$. The integral appearing in
this equation is a kind of non-linear Young integral~{\cite{Young}}. Existence
and uniqueness of solutions for equations involving Young integrals are by now
standard~{\cite{Lyons,controlling,FV}} and easily extended to this context as
shown below. In particular the equation
\[ \theta_{t} = \theta_{0} + \int_{0}^{t} G_{\mathd s} ( \theta_{s} ) \]
will have a solution $\theta \in \CC^{\gamma} ([0,1],\mathbbm{R}^{d} )$ (the
space of $\gamma$-Hölder continuous functions from $[0,T]$ to $\mathbbm{R}^{d}$) provided
$(x,t) \mapsto G_{t} (x)$ is a $\gamma$-Hölder function of time, locally
Lipshitz in space with $\gamma >1/2$, that is
\[ |G_{s,t} (x)-G_{s,t} (y)| \lesssim |x-y|  |t-s|^{\gamma} \psi ( | x | + | y
   | ) \]
for all $x,y \in \mathbbm{R}^{d}$ and $0 \leqslant s \leqslant t \leqslant 1$.
Note that some space regularity is already needed to have existence (to be
compared with the classical setup where bounded vector fields are sufficient for existence).

A strategy to prove uniqueness is to consider the difference between a
solutions $\theta$ and a solution $\theta'$ of a similar equations
\[ \theta'_{t} =x_{0} + \int_{0}^{t} G_{\mathd u}' (\theta'_{u} ) . \]
It is the necessary to estimate the difference
\[ (\theta_{t} -\theta'_{t} ) - (\theta_{s} -\theta'_{s} ) = \int_{s}^{t} G_{\mathd u} (\theta_{u} )
   -G_{\mathd u}' (\theta'_{u} ). \]
To deal with such an estimates, we will need an averaged translation operator
$\tau_{f}  G_{s,t} ( z ) = \int_{s}^{t} G_{\mathd u} ( f_{u} +z )$ in order to
have an equation on $\theta-\theta'$.

In order for these estimates to be useful we need a way to link the
regularity of the original vector field $b$ with its averaged version
$T^{w}_{s,t} b$ along an arbitrary continuous path $w$.

\begin{theorem}
  Assume that for $\alpha \in \mathbbm{R}$, $f \to T^{w} f$ is defined on the
  whole space $\CF L^{\alpha + \nu}$ for all $\nu \geqslant 0$. Assume also
  that there exists $\gamma >1/2$ such that for all $\nu >0$, there exists a
  $\nu$-weight $\psi$ such that \ $T^{w}$ maps $\CC^{\alpha + \nu}$ into
  $\CC^{\gamma , \nu , \psi}$. Then there exists a solution $\theta (x_{0} )
  \in \CC^{\gamma} ([0,1],\mathbbm{R}^{d} )$ to the Young-type equation
  \[ \theta_{t} (x_{0} ) =x_{0} + \int_{0}^{t} (T^{w}_{\mathd s} b) (
     \theta_{s} (x_{0} )) \]
  for any $b \in \CC^{\alpha + \nu}$ for $\nu$ such that $\gamma ( 1+ \nu )
  >1$. If $b \in \CC^{\alpha +2}$ (or $\alpha +3/2>0$ and $b \in \CC^{\alpha
  +3/2}$) this is the unique $\gamma$-Hölder solution to this equation, and
  for all $t \in [0,1]$, the flow map $x_{0} \to \theta_{t} (x_{0} )$ is well
  defined and locally Lipschitz continuous, uniformly in time.
\end{theorem}

\begin{remark}
  To prove such a theorem, we need the two hypothesis about $T^{w}$. The first
  one is that this map is well defined. This will follow either from
  the definition of the map (when $\alpha \ge 0$) or from section
  \ref{sec:averaging_deterministic}. The second one is to prove that $T^{w}$ maps
  $\CC^{\alpha}$ into $\CC^{\gamma , \nu , \psi}$. We also need a theory of
  integration for vector fields in $\CC^{\gamma , \nu , \psi}$. In the next
  section we will build such a theory.
\end{remark}

This theorem is obtained when we apply Theorem \ref{theorem:existence_young_solution}, Remark \ref{remark:localization_principle}, Theorem \ref{theorem:uniqueness_general_case} and Corollary \ref{corollary:uniqueness_a_priori_case} to the operator $T^w$ with the wanted hypothesis.

When $\alpha \geqslant -1$ the vector field $b \in \CC^{\alpha +1}$ is
continuous and the solutions are simply solutions to the classical ODE
\[ \theta_{t} = \theta_{0} + \int_{0}^{t} b (u,w_{u} + \theta_{u} ) \mathd u.
\]
In the case that $\alpha <-1$ the vector field $b$ is a distribution and the
previous ODE does not make sense. In that situation the natural meaning of
these solution is the following. Let $b_{n} = \pi_{\leqslant n} b$ for $b \in
\CC^{\alpha}$ then
\[ \int_{0}^{t} b_{n}  (u,w_{u} + \theta_{u} ) \mathd u= \int_{0}^{t}
   (T^{w}_{\mathd s} b_{n} ) ( \theta_{s} ) \to \int_{0}^{t} (T^{w}_{\mathd s}
   b) ( \theta_{s} ) \]
by continuity of the Young integral and of the averaging with respect to the
norm of $\CF L^{\alpha +1}$. Then $\theta$ solves the equation
\[ \theta_{t} = \theta_{0} + \lim_{n}   \int_{0}^{t} b_{n}  (u,w_{u} +
   \theta_{u} ) \mathd u ,\]
where the right hand side is well defined for any $\theta \in \CC^{\gamma}
([0,1],\mathbbm{R}^{d} )$. At this point we can identify
\[ \int_{0}^{t} b (u,w_{u} + \theta_{u} ) \mathd u= \lim_{n}   \int_{0}^{t}
   b_{n}  (u,w_{u} + \theta_{u} ) \mathd u \]
and give meaning to the ODE with a distributional drift $b$.

\begin{remark}
  When the vector field $b$ is in $\CF L^{\alpha}$, the limiting procedure
  does not depend on the choice of the sequence. That is the principal reason
  of the introduction of that spaces.
\end{remark}

One of the aims of this paper is to show that the above program can be carried
out successfully in the case of $w$ given by a sample path of a fractional
Brownian motion $B^{H}$ of Hurst parameter $H \in (0,1)$.

\subsection{Definition of the Young integral}

We define now the Young integral~{\cite{Young,Lyons,controlling}} for non
linear operators.

\begin{theorem}
  \label{theorem:young integral}Let $\gamma , \rho , \nu >0$ with $\gamma +
  \nu \rho >1$, a $\nu$-weight $\psi$, and $V$ and $W$ two Banach spaces and
  $I$ a finite interval on $\mathbbm{R}$. Let $G \in \CC^{\gamma , \nu , \psi}
  ( I,V,W )$ and $f \in \CC^{\rho} (I;V)$. Let $s,t \in I$ with $s \leqslant
  t$. Then the following limit exists and is independent of the partition
  \[ \int_{s}^{t} G_{\mathd u} (f_{u} ) \assign
     \lim_{\tmscript{\begin{array}{c}
         \text{$\Pi$ partition of $[s,t]$}\\
       | \Pi | \to 0
     \end{array}}}   \sum_{i} G_{t_{i} ,t_{i+1}} (f_{t_{i}} )
. \]
  Furthermore
  \begin{enumerate}
    \item For all $s \leqslant u \leqslant t$ with $s,u,t \in I$ we have
    \[ \int_{s}^{t} G_{\mathd r} (f_{r} ) = \int_{s}^{u} G_{\mathd r} (f_{r} )
       + \int_{u}^{t} G_{\mathd r} (f_{r} ) . \]
    \item
    \[ \left\lvert \int_{s}^{t} G_{\mathd r} (f_{r} )-G_{s,t} (f_{s} ) \right\lvert_{W}
       \leqslant C_{\gamma , \rho , \nu} \llbracket G \rrbracket_{\gamma , \nu
       , \psi} \llbracket f \rrbracket_{\rho ,I}^{\nu}  |t-s|^{\gamma + \nu
       \rho} \psi ( \| f \|_{\infty ,I} ) . \]
    \item For all $s \leqslant t \in I$ and $R>0$, the map $(f,G) \mapsto
    \int_{s}^{t} G_{\mathd r} (f_{r} )$ is continuous as a function of $\left(
    \left\{ g \in \CC^{\rho} (I,V), \| g \|_{\gamma , [ s,t ]} \leqslant R
    \right\} ,\|.\|_{\infty ,[s,t]} \right) \times ( \CC^{\gamma , \nu , \psi} (
    I,V,W ) ,\|.\|_{\gamma , \nu , \psi} )$ to $W$.
  \end{enumerate}
\end{theorem}

\begin{proof}
  Let $s,t \in I$ with $s \leqslant t$ be fixed until the end of the proof.
  Suppose first that $G$ is differentiable (in time) and $G' \in \CC^{\gamma ,
  \nu , \psi} (I,V,W)$ and $G \in \CC^{\gamma , \nu , \psi}$. For simplicity,
  in all the proof we write $\| G \|$ and $\llbracket G \rrbracket$ instead of
  $\| G \|_{\gamma , \nu , \psi}$ and $\llbracket G \rrbracket_{\gamma , \nu ,
  \psi}$. Then we define for $s \leqslant t$
  \[ \int_{s}^{t} G_{\mathd u} (f_{u} ) \assign \int_{s}^{t} G'_{u} (f_{u} )
     \mathd u \assign I_{s,t} (f,G) \]
  and also define $J_{s,t} (f,G) \assign I_{s,t} (f,G) -G_{s,t} (f_{s} )$. For
  $u \in [s,t]$ we have
  \[ J_{s,t} (f,G) =J_{s,u} (f,G) +J_{u,t} (f,G) +G_{u,t} (f_{u} ) - \delta
     G_{u,t} (f_{s} ), \]
  hence, for $n \ge 1$, $i \in \{0, \ldots ,2^{n} \}$ and $t_{i}^{n} =s+ (t-s)
  i2^{-n}$,
  \[ J_{s,t} (f,G) = \sum_{i=0}^{2^{n} -1} J_{t_{i}^{n} ,t_{i+1}^{n}} (f,G) +
     \sum_{k=1}^{n} \sum_{i=1}^{2^{k} -1} (G_{t_{2i-1}^{k} ,t_{2i}^{k}}
     (f_{t_{2i-1}^{k}} )-G_{t_{2i-1}^{k} ,t_{2i}^{k}} (f_{t_{2(i-1)}^{k}} )).
  \]
  But, as $G$ is differentiable, the following computation holds
  \begin{eqnarray*}
    | J_{t_{i}^{n} ,t_{i+1}^{n}} (f,G) |_{W} & \leqslant &
    \int_{t_{i}^{n}}^{t_{i+1}^{n}} | G_{u}' (f_{u} )-G_{u}' (f_{t_{i}^{n}} )
    |_{W} \mathd u\\
    & \leqslant & \int_{t_{i}^{n}}^{t_{i+1}^{n}}  \llbracket G'_{u}
    \rrbracket_{\nu , \psi} | f_{u} -f_{t_{i}^{n}} |^{\nu}_{V} \psi ( | f_{u}
    | + | f_{t_{i}^{n}} | ) \mathd u\\
    & \leqslant & 2 \| \llbracket G' \rrbracket_{\nu , \psi} \|_{\infty}
    \int_{t_{i}^{n}}^{t_{i+1}^{n}} | f |^{\nu}_{\rho}  |u-t_{i}^{n} |^{\nu
    \rho} \psi ( \| f \|_{\infty , [ s,t ]} ) \mathd u\\
    & \lesssim & 2^{- (1+ \nu \rho ) n}.
  \end{eqnarray*}
  Hence
  \[ \left\lvert \sum_{i=0}^{2^{n} -1} J_{t_{i}^{n} ,t_{i+1}^{n}} (f,G) \right\lvert
     \lesssim 2^{-n \nu \rho} \underset{n \to \infty}{\to} 0, \]
  and then
  \[ | J_{s,t} (f,G) | \leqslant \sum_{k=1}^{\infty} \sum_{i=1}^{2^{k} -1}
     |G_{t_{2i-1}^{k} ,t_{2i}^{k}} (f_{t_{2i-1}^{k}} )-G_{t_{2i-1}^{k}
     ,t_{2i}^{k}} (f_{t_{2(i-1)}^{k}} )|_{W} . \]
  For $k \ge 1$ and $i \in \{1, \ldots ,2^{k} -1\}$, we have
  \begin{equation*} \begin{split}
    | G_{t_{2i-1}^{k} ,t_{2i}^{k}} (f_{t_{2i-1}^{k}} )& -G_{t_{2i-1}^{k}
    ,t_{2i}^{k}} (f_{t_{2(i-1)}^{k}} ) |  \\ & \leqslant  \llbracket G
    \rrbracket | t_{2i-1}^{k} -t^{k}_{2i} |^{\gamma} | f_{t_{2i-1}^{k}}
    -f_{t_{2i-2}^{k}} |^{\nu} \psi ( | f_{t_{2i-1}^{k}} | + | f_{t_{2i-2}^{k}}
    | )\\
    & \lesssim  \llbracket G \rrbracket \llbracket f \rrbracket^{\nu}_{\rho
    , [s,t]} \psi ( \|f\|_{\infty ,[s,t]} ) |t-s|^{\gamma + \nu \rho} 2^{- (
    \gamma + \nu \rho ) k}.
    \end{split}
    \end{equation*}
  Hence, the following bound holds
  \begin{equation} \begin{split}
       \sum_{k=1}^{\infty} \sum_{i=1}^{2^{k} -1} & |G_{t_{2i-1}^{k}
       ,t_{2i}^{k}} (f_{t_{2i-1}^{k}} )-G_{t_{2i-1}^{k} ,t_{2i}^{k}}
       (f_{t_{2(i-1)}^{k}} )|_{W}\\
       & \lesssim \llbracket G \rrbracket \llbracket f \rrbracket^{\nu}_{\rho
       , [s,t]} \psi ( \| f \|_{\infty ,[s,t]} )  |t-s|^{\nu \rho + \gamma} 
       \sum_{k=1}^{\infty} \sum_{i=1}^{2^{k} -1} 2^{- ( \gamma + \nu \rho )
       k}\\
       & \lesssim \frac{2^{- ( \nu \rho + \gamma -1)}}{1-2^{- ( \nu \rho +
       \gamma -1)}} \llbracket G \rrbracket \llbracket f
       \rrbracket^{\nu}_{\rho , [s,t]} \psi ( \|f\|_{\infty ,[s,t]} ) 
       |t-s|^{\nu \rho + \gamma} .
     \end{split} \end{equation}
  The result is proved for differentiable $G$. Let us now take $G \in \CC^{\gamma , \nu ,
  \psi} (I,V,W)$ and $f$ as wanted. Let $G^{n}$ be differentiable as above such that
  $G_{s,t}^{n} (f_{s} ) \to G_{s,t} (f_{s} )$ as $n \to \infty$; for all
  $\gamma' < \gamma$ $\lim_{n \to \infty} \|G-G^{n} \|_{\gamma' , \nu , \psi}
  =0$ and for all $n \ge 0$, $\|G^{n} \|_{\gamma , \nu , \psi} \leqslant
  \|G\|_{\gamma , \nu , \psi}$. As $I_{s,t}$ is linear in the second variable,
  we have, for $\gamma' < \gamma$
  \begin{eqnarray*}
    |J_{s,t} (f,G^{n} )-J_{s,t} (f,G^{n+m} )|_{W} & = & |J_{s,t} (f,G^{n}
    -G^{n+m} )|_{W}\\
    & \lesssim & \llbracket G^{n} -G^{n+m} \rrbracket_{\gamma' , \nu}\\
    & \underset{n \to \infty}{\to} & 0.
  \end{eqnarray*}
  The sequence $( J_{s,t} (f,G^{n} ) )_{n}$ is Cauchy in $W$ which is a Banach
  space. Let us say it converges to a number $J_{s,t} (f,G)$. Furthermore, the sequence
  $G_{s,t}^{n} (f_{s} )$ converges obviously to $G_{s,t} (f_{s} )$. Then as
  $I_{s,t} (f,G^{n} ) =J_{s,t} (f,G^{n} ) +G_{s,t}^{n} (f_{s} )$ the sequence
  $( I_{s,t} (f,G^{n} ) )_{n}$ converges to a limit called $I_{s,t} (f,G)$.
  Furthermore,
  \begin{eqnarray*}
    |J_{s,t} (f,G^{n} )|_{W} & \lesssim_{\gamma , \rho} & \llbracket G^{n}
    \rrbracket \llbracket f \rrbracket_{\rho , [s,t]} \psi ( \|f\|_{\infty} )
    |t-s|^{\gamma + \nu \rho}\\
    & \lesssim & \llbracket G \rrbracket \llbracket f \rrbracket_{\rho} \psi
    ( \|f\|_{\infty} ) |t-s|^{\gamma + \nu \rho}
  \end{eqnarray*}
  and so does $|I_{s,t} (f,G)-G_{s,t} (f_{s} )|_{W}$. The Chasles property and
  the triangular inequality are obvious with the definition of $I$. Moreover
  since $I (f,G)$ is linear in $G$ it is easy to see that the definition does
  not depend on the particular sequence $G^{n}$.
  
  Let us show that $I_{s,t} (f,G)$ is the limit of Riemann sum. Let $\Pi =
  \{s=t_{0} <t_{1} < \ldots <t_{n} =t\}$ a partition of $[s,t]$. Let
  \[ S_{\Pi} = \sum_{k=0}^{n-1} G_{t_{i} ,t_{i+1}} (f_{t_{i+1}} ) \]
  be the Riemann sum corresponding to this partition. As $G_{t_{i} ,t_{i+1}}
  (f_{t_{i+1}} ) =I_{t_{i} ,t_{i+1}} (f,G) -J_{t_{i} ,t_{i+1}} (f,G)$ the
  following equality holds
  \[ S_{\Pi} -I_{s,t} (f,G) =- \sum_{i=0}^{n-1} J_{t_{i} ,t_{i+1}} (f,G). \]
  Hence
  \begin{equation*}
    |S_{\Pi} -I_{s,t} (f,G)|_{W}  \leqslant  \sum_{i=0}^{n-1} |J_{t_{i}
    ,t_{i+1}} (f,G)|_{W}
    \lesssim   \sum_{i=0}^{n-1} |t_{i+1} -t_{i} |^{\gamma + \nu \rho}
     \lesssim  | \Pi |^{\gamma + \nu \rho -1} \to_{| \Pi | \to 0} 0.
  \end{equation*}

  It remains to show the continuity of the map $(f,G) \mapsto I (f,G)$. Take
  $f,f' ,G,G'$ and assume for simplicity that $G (0) =G' (0) =0$ then
  \begin{eqnarray*}
    I_{s,t} (f,G) -I_{s,t} (f' ,G' ) & = & [I_{s,t} (f,G)-I_{s,t} (f' ,G)]
    +I_{s,t}  (f' ,G-G' )
  \end{eqnarray*}
  and
  \begin{eqnarray*}
    | I_{s,t} (f' ,G-G' ) |_{W} & \leqslant & | (G-G' )_{s,t} (f_{s} )-(G-G'
    )_{s,t} (0) |_{W} + | (G-G' )_{s,t} (0) | + | J_{s,t} (f,G-G' ) |_{W}\\
    & \lesssim & \| G-G' \| | t-s |^{\gamma} ( | f_{s} |^{\gamma} \psi ( |
    f_{s} | ) +1+|t-s|^{\nu \rho} \llbracket f \rrbracket^{\nu}_{\rho} \psi (
    \|f\|_{\infty} ) )\\
    & \lesssim & \| G-G' \| |t-s|^{\gamma}  ( \|f\|_{\infty}^{\nu} \psi (
    \|f\|_{\infty} ) +1+|t-s|^{\nu \rho} \llbracket f \rrbracket^{\nu}_{\rho}
    \psi ( \|f\|_{\infty} ) ) .\\
    & \lesssim & \| G-G' \| |t-s|^{\gamma}  ( 1+ \| f \|^{\nu}_{\rho} \psi (
    \| f \|_{\rho} ) ) .
  \end{eqnarray*}
  Furthermore
  \[ I_{s,t} (f,G) -I_{s,t} (f' ,G) =G_{s,t} (f_{s} ) -G_{s,t} (f'_{s} )
     +J_{s,t} (f,G) -J_{s,t} (f' ,G) . \]
  We have also
  \begin{eqnarray*}
    |I_{s,t} (f,G)-I_{s,t} (f' ,G)|_{W} & \lesssim & \|G\| \|f-f'
    \|^{\nu}_{\infty} \psi ( \| f \|_{\infty} + \| f' \|_{\infty} )
    |t-s|^{\gamma}\\
    &  & + ( \llbracket f \rrbracket^{\nu}_{\rho} \psi ( \| f \|_{\infty} ) +
    \llbracket f' \rrbracket^{\nu}_{\rho} \psi ( \| f' \|_{\infty} ) ) \|G\|
    |t-s|^{\nu \rho + \gamma} .\\
    & \lesssim & \|G\| \|f-f' \|^{\nu}_{\infty} \psi ( \| f \|_{\rho} + \| f'
    \|_{\rho} ) |t-s|^{\gamma}\\
    &  & + ( \| f \|^{\nu}_{\rho} \psi ( \| f \|_{\rho} ) + \nobracket \|
    \nobracket f' \|^{\nu}_{\rho} \psi ( \| f' \|_{\rho} ) ) \|G\| |t-s|^{\nu
    \rho + \gamma} .
  \end{eqnarray*}
  By partitioning the interval $[s,t]$ in subintervals $[t_{i} ,t_{i+1} ]$ of
  size $2^{-n}$ and summing up the contributions according to these bounds we
  obtain an improved estimate
  \begin{eqnarray*}
    |I_{s,t} (f,G)-I_{s,t} (f' ,G)|_{W} & \lesssim & \|G\|  \|f-f'
    \|^{\nu}_{\infty} \psi ( \| f \|_{\rho} + \| f' \|_{\rho} ) 2^{(1- \gamma
    ) n}  |t-s|^{\gamma}\\
    &  & + \|G\| ( \| f \|^{\nu}_{\rho} \psi ( \| f \|_{\rho} ) + \| f'
    \|^{\nu}_{\rho} \psi ( \| f' \|_{\rho} ) )  ( 2^{-n} |t-s| )^{\nu \rho +
    \gamma} 2^{n}.
  \end{eqnarray*}
  Taking $n$ large enough so that
  \[ 2^{- \nu \rho n} \leqslant \frac{\|f-f' \|^{\nu}_{\infty} \psi ( \| f
     \|_{\rho} + \| f' \|_{\rho} )}{( \| f \|^{\nu}_{\rho} \psi ( \| f
     \|_{\rho} ) + \| f' \|^{\nu}_{\rho} \psi ( \| f' \|_{\rho} ) ) |t-s|^{\nu
     \rho}} \leqslant 2^{- \nu \rho ( n-1 )}, \]
  we have
  \[ |I_{s,t} (f,G)-I_{s,t} (f' ,G)|_{W} \lesssim \|G\|  \|f-f'
     \|^{\nu}_{\infty} \psi ( \| f \|_{\rho} + \| f' \|_{\rho} ) 2^{(1- \gamma
     ) n}  |t-s|^{\gamma}, \]
  which means that it is possible to choose $n$ such that
  \begin{equation*}\begin{split} & |I_{s,t} (f,G)-I_{s,t} (f' ,G)|_{W}   \lesssim \|G\|  \|f-f'
     \|^{\nu}_{\infty} \psi ( \| f \|_{\rho} + \| f' \|_{\rho} ) 
     |t-s|^{\gamma}  \\  
     		& \qquad \qquad \qquad  \times \left( \frac{( \| f \|^{\nu}_{\rho} \psi ( \| f \|_{\rho}
     ) + \| f' \|^{\nu}_{\rho} \psi ( \| f' \|_{\rho} ) )|t-s|^{\nu \rho}
     }{\|f-f' \|^{\nu}_{\infty} \psi ( \| f \|_{\rho} + \| f' \|_{\rho} )}
     \right)^{\frac{1- \gamma}{\nu \rho}}
     \\  
     		\qquad & \lesssim \|G\| \{ \|f-f' \|^{\nu}_{\infty} \psi ( \| f \|_{\rho} + \| f'
     \|_{\rho} ) \}^{( \gamma + \nu \rho -1 ) / \nu \rho}  |t-s|  \\ & \qquad \qquad \times ( \| f
     \|^{\nu}_{\rho} \psi ( \| f \|_{\rho} ) + \| f' \|^{\nu}_{\rho} \psi ( \|
     f' \|_{\rho} ) )^{(1- \gamma ) / ( \nu \rho )} 
     \end{split}\end{equation*}
  and this allows us to infer the continuity of $I ( f,G )$.
\end{proof}

\begin{remark}
  It is easy to construct a suitable sequence $(G^{n} )_{n \ge 1}$. Let
  $h:\mathbbm{R} \to \mathbbm{R}$ be a compactly supported, smooth positive
  function with integral $1$. Define $h_{n} (t) =nh (nt)$ and define for all
  $v \in V$ and all $t \in \mathbbm{R}$
  \[ G_{t}^{n} (v) = \int_{\mathbbm{R}} h_{n}  (t-s) G_{s} (v) \mathd
     s=G_{t}^{n} (v) = \int_{\mathbbm{R}} h_{n} (s) G_{t-s} (v) \mathd s .\]
  Then $G^{n}$ is as wanted. Indeed,
  \begin{equation} \begin{split}
       & | G_{s,t}^{n} (v)-G_{s,t}^{n} (w) |_{W}\\
       & \hspace{2em} \leqslant \int_{\mathbbm{R}} h_{n} (r) | ( G_{t-r}
       -G_{s-r} ) (v)- ( G_{t-r} -G_{s-r} ) (w) |_{W} \mathd r\\
       & \hspace{2em} \leqslant \int_{\mathbbm{R}} h_{n} (r) | t-s |^{\gamma}
       |v-w|^{\nu}_{V} \psi ( | v | + | w | ) \mathd r\\
       & \hspace{2em} \leqslant \|G\|_{\gamma , \nu , \psi} |t-s|^{\gamma}
       |v-w|^{\nu}_{V}   \psi ( | v | + | w | )
     \end{split} \end{equation}
  which proves that $G^{n} \in \CC^{\gamma , \nu , \psi} (\mathbbm{R},V,W) )$
  and that $\|G^{n} \|_{\gamma , \nu , \psi} \leqslant \|G\|_{\gamma , \nu ,
  \psi}$. Furthermore $G^{n}$ is differentiable and $(G^{n} )' \in \CC^{\gamma ,
  \nu , \psi} (\mathbbm{R},V,W)$. As we can chose $h_{n}$ to be a good kernel,
  all the properties required on $G^{n}$ are satisfied.
\end{remark}

\begin{definition}
  The limit functional $I$ defined in the last theorem is obviously an
  integral and then we will refer to it as $\int_{s}^{t} G_{du} (f_{u} )$.
\end{definition}

\begin{remark}
  Let $g \in \CC^{\gamma} (I,V' )$ and $f \in \CC^{\rho} (I,V)$ with $\gamma +
  \rho >1$, where $V$ and $V'$ are (finite--dimensional) Banach spaces. Let
  $W=V \otimes V'$ and for all $x \in V$, $G_{t} (x) =x \otimes g_{t}$. Then
  $G \in \CC^{\gamma ,1} (I,V,W)$ and the above integral is the standard Young
  integral.
\end{remark}

\begin{remark}
  The bound in Theorem \ref{theorem:young integral} is
  \[ \left\lvert \int_{s}^{t} G_{\mathd u} ( f_{u} ) -G_{s,t} ( f_{s} ) \right\lvert
     \lesssim \llbracket G \rrbracket_{\gamma , \nu , \psi} | t-s |^{\gamma +
     \rho \nu} \llbracket f \rrbracket_{\gamma}^{\nu} \psi ( \| f \|_{\infty}
     ) . \]
  But as $\llbracket f \rrbracket_{\gamma} \leqslant \| f \|_{\gamma}$ and $\|
  f \|_{\infty} \leqslant ( 1+ | I | ) \| f \|_{\gamma}$ and $\psi$ is a
  weight, we also have this other useful bound
  \[ \left\lvert \int_{s}^{t} G_{\mathd u} ( f_{u} ) -G_{s,t} ( f_{s} ) \right\lvert
     \lesssim \llbracket G \rrbracket_{\gamma , \nu , \psi} | t-s |^{\gamma +
     \rho \nu} \| f \|_{\gamma}^{\nu} \psi ( \| f \|_{\gamma} ) ,\]
  where the new constant depends on the length of the interval $| I |$ and
  $\psi$. In the following, we will exploit these three bounds equally
  and without further notice.
\end{remark}

We intend to solve differential equations driven by such $G$. Thanks to the
definition of the integral and the bound in Theorem \ref{theorem:young
integral}, we are able to define the equation, prove the existence of
solutions and give an a priori bound on the norm of the solutions. Here we
will use the notion of $\nu$-weight, in order to control the growth of the
norm.

\begin{theorem}
  \label{theorem:existence_young_solution}Let $\gamma > \frac{1}{2}$, $\nu \in
  [ 0,1 )$ such that $\gamma ( 1+ \nu ) >1$ and $\psi$ be a $\nu$-weight as in
  Definition $\ref{definition:nu_weight}$. Let $G \in \CC^{\gamma , \nu ,
  \psi} ( [ 0,T ] ,\mathbbm{R}^{d} )$ and $x \in \mathbbm{R}^{d}$. There
  exists a solution $\theta\in\CC^{\gamma}([0,T];\RR^d)$ to the non-linear Young differential equation
  \[ \theta_{t} = \theta_{0} + \int_{0}^{t} G_{\mathd u} ( \theta_{u} ) . \]
  Furthermore, there exists two universal constants $K_{1}$ and $K_{2}$
  depending on $\gamma , \nu , \psi$ and $T$ such that
  \[ \| \theta \|_{\infty , [ 0,T ]} \leqslant K_{1} ( 1+ \| G \|_{\gamma ,
     \nu , \psi} )^{K_{2} \| G \|_{\gamma , \nu , \psi}^{1/ \nu \gamma}} ( |
     \theta_{0} | +1 ) \]
  and
  \[ \| \theta \|_{\gamma , [ 0,T ]} \leqslant \tilde{K}_{1} \| G \|^{1/ \nu
     \gamma} ( 1+ \| G \|_{\gamma , \nu , \psi} )^{\tilde{K}_{2} \| G
     \|_{\gamma , \nu , \psi}^{1/ \nu \gamma}} ( 1+ | \theta_{0} | ) . \]
\end{theorem}

\begin{proof}
  Let us first deal with the existence of the solutions. Let $t_{0} \in I= [
  0,T ]$, $K>0$ and $0<S \leqslant T$ to be specify later. Let $J= [t_{0}
  ,(t_{0} +S) \wedge T]$ and let us define for all $x \in V$,
  \[ \mathcal{C}_{t_{0} ,x} = \{\theta \in \CC^{\gamma} (J): \theta_{t_{0}}
     =x,\| \theta \|_{\gamma ,J} \leqslant K\} \]
  and
  \[ \Phi_{t_{0} ,x} : \begin{array}{ccc}
       \CC^{\gamma} (J,V) & \rightarrow & \CC^{\gamma} (J \nocomma ,V)\\
       \theta & \rightarrow & x+ \int_{t_{0}}^{.} G_{\mathd u} ( \theta_{u} ).
     \end{array} \]
  By Theorem~\ref{theorem:young integral} the map $\Phi_{t_{0} ,x}$ is well defined.
  Furthermore we always have
  \[ | \theta_{s} | \leqslant | \theta_{t_{0}} | +T^{\nu} \llbracket \theta
     \rrbracket_{\gamma ,J} \lesssim_{T} \| \theta \|_{\gamma ,J} .\]
  Hence for $s<t \in J$ we have
  \begin{eqnarray}
    | \delta ( \Phi_{t_{0} ,x} ( \theta ))_{s,t} | & \leqslant & \left\lvert
    \int_{s}^{t} G_{\mathd u} ( \theta_{u} )-G_{s,t} ( \theta_{s} ) \right\lvert +
    |G_{s,t} ( \theta_{s} )-G_{s,t} (0)| + |G_{s,t} (0)| \nonumber\\
    & \lesssim & \|G\|_{\gamma , \nu , \psi}  |t-s|^{\gamma}  ( S^{\nu
    \gamma} \| \theta \|^{\nu}_{\gamma ,J} \psi ( \| \theta \|_{\gamma ,J} ) +
    | \theta_{s} |^{\nu} \psi ( | \theta_{s} | ) +1 ) \nonumber\\
    & \lesssim & \|G\|_{\gamma , \nu , \psi}  |t-s|^{\gamma}  ( S^{\nu
    \gamma} \| \theta \|^{\nu}_{\gamma ,J} \psi ( \| \theta \|_{\gamma ,J} )
    +| \theta_{t_{0}} |^{\nu} \psi ( | \theta_{0} | ) +1 ) . \label{majoration
    bound}
  \end{eqnarray}
  Now take $\theta \in \mathcal{C}_{t_{0} ,x}$,
  \begin{eqnarray*}
    \llbracket \Phi_{t_{0} ,x} ( \theta ) \rrbracket_{\gamma ,J} & \lesssim &
    \|G\|_{\gamma} (S^{\nu \gamma} \| \theta \|^{\nu}_{\gamma ,J} \psi ( \|
    \theta \|_{\gamma ,J} ) +|x|^{\nu} \psi ( | x | ) +1).
  \end{eqnarray*}
  But since $\nu <1$ and $\psi$ is a $\nu$-weight,, there exists a constant
  $C_{\nu , \psi}$ such that $\| \theta \|^{\nu}_{\gamma ,J} \psi ( \| \theta
  \|_{\gamma ,J} ) \leqslant C_{\nu , \psi} ( 1+ | x | )$. Hence, there is a
  universal constant $C>0$ such that
  \[ \| \Phi_{t_{0} ,x} ( \theta ) \|_{\gamma ,J} = \llbracket \Phi_{t_{0} ,x}
     ( \theta ) \rrbracket_{\gamma ,J} + | x | \leqslant | x | +C
     \|G\|_{\gamma} (S^{\nu \gamma} \| \theta \|_{\gamma ,J} +|x|^{\nu} \psi (
     | x | ) +1) . \]
  For $S$ such that $C \| G \|_{\gamma} S^{\nu \gamma} <1/2$, and for $K
  \geqslant 2 \{ | x | +C ( |x|^{\nu} \psi ( | x | ) +1 ) \}$, we have
  \[ \| \Phi_{t_{0} ,x} ( \theta ) \|_{\gamma ,J} \leqslant K. \]
  Then $\Phi_{t_{0} ,x} ( \theta ) \in \mathcal{C}_{t_{0} ,x}$, moreover by
  the property of the Young integral the map $\Phi_{t_{0} ,x}$ is continuous
  on $\mathcal{C}_{t_{0} ,x}$ for the norm $\|.\|_{\infty , [t_{0} ,(t_{0} +T)
  \wedge 1]}$. By its definition $\mathcal{C}_{t_{0} ,x}$ is immediately a
  closed convex set of $C ( J )$. Let us show that $\Phi_{t_{0} ,x} (
  \mathcal{C}_{t_{0} ,x} )$ is relatively compact in $\mathcal{C}^{0}$. It is
  obviously equicontinuous as $\| \Phi_{t_{0} ,x} ( \theta )\|_{\gamma}
  \leqslant K$ and relatively bounded as $| \Phi_{t_{0} ,x} ( \theta )_{t} |
  \leqslant |x| +K^{\nu} \psi ( K ) (t-t_{0} )^{\gamma}$. Hence by Ascoli
  theorem $\Phi_{t_{0} ,x} ( \mathcal{C}_{t_{0} ,x} )$ is relatively compact.
  Thanks to Leray-Schauder-Tychonoff fixed point theorem, there exists
  $\theta^{t_{0} ,x}$ such that $\theta^{t_{0} ,x} = \Phi_{t_{0} ,x} (
  \theta^{t_{0} ,x} ) =x+ \int_{t_{0}}^{t} G_{\mathd u} ( \theta^{t_{0}
  ,x}_{u} )$. We then construct by induction a solution on the whole interval.
  For $n$ such that $n S \leqslant T$ let $\theta^{0} = \theta^{0,x_{0}}$ and
  $\theta^{n} = \theta^{nS, \theta^{n-1}_{S}}$. Let us define $\theta_{t} =
  \theta^{n}_{t}$ if $t \in [nS,(n+1)S]$. By an immediate induction, $\theta$
  is solution of the equation $\theta_{t} =x_{0} + \int_{0}^{t} G_{\mathd u} (
  \theta_{u} )$ and then is obviously in $\CC^{\gamma}$.
  
  We have all the tools to bound the norm of a solution of the equation. Again
  take $t_{0}$ and $S$ to be specify lated, and $\theta$ a solution of the
  non-linear Young differential equation. And take $J= [ t_{0} , ( t_{0} +S )
  \wedge T ]$. We have
  \begin{eqnarray*}
    | \delta \theta_{s,t} | & \leqslant & \left\lvert \int_{s}^{t} G_{\mathd u} (
    \theta_{u} ) -G_{s,t} ( \theta_{s} ) \right\lvert + | G_{s,t} ( \theta_{s} )
    -G_{s,t} ( \theta_{t_{0}} ) | \\
    & & + | G_{s,t} ( \theta_{t_{0}} ) -G_{s,t} ( 0
    ) | + | G_{s,t} ( 0 ) |\\
    & \lesssim & | t-s |^{\gamma} \| G \|_{\gamma , \nu , \psi} ( S^{\gamma
    \nu} \llbracket \theta \rrbracket_{\gamma ,J}^{\nu} \psi ( \| \theta
    \|_{\gamma ,J} ) + | \theta_{t_{0}} |^{\nu} \psi ( | \theta_{t_{0}} | ) +1,
    )
  \end{eqnarray*}
  hence
  \[ \llbracket \theta \rrbracket_{\gamma , [ t_{0} ,t_{0} +S ]} \lesssim \| G
     \|_{\gamma , \nu , \psi} ( S^{\gamma \nu} \| \theta \|_{\gamma ,J}^{\nu}
     \psi ( \| \theta \|_{\gamma ,J} ) + | \theta_{t_{0}} |^{\nu} \psi ( |
     \theta_{t_{0}} | ) +1 ) . \]
  Let $S$ be such that $C \| G \|_{\gamma , \nu , \psi} S^{\gamma , \nu}
  \leqslant 1$, and we have
  \[ \| \theta \|_{\gamma , [ t_{0} ,t_{0} +S ]} \lesssim \| \theta \|_{\gamma
     ,J}^{\nu} \psi ( \| \theta \|_{\gamma ,J} ) + | \theta_{t_{0}} | +C \| G
     \|_{\gamma , \nu , \psi} | \theta_{t_{0}} |^{\nu} \psi ( | \theta_{t_{0}}
     | ) +C \| G \|_{\gamma , \nu} .\]
  As $x \rightarrow x^{\nu} \psi ( x )$ is sublinear (as before), there
  exists a constant depending on $\nu$ and $\psi$ such that
  \[ \| \theta \|_{\gamma , [ t_{0} ,t_{0} +S ]} \lesssim_{\nu , \psi} 1+ |
     \theta_{t_{0}} | + \| G \|_{\gamma , \nu , \psi} | \theta_{t_{0}} |^{\nu}
     \psi ( | \theta_{t_{0}} | ) + \| G \|_{\gamma , \nu , \psi} . \]
  There also exists a constant $C$ such that
  \[ | \theta_{t_{0}} |^{\nu} \psi ( | \theta_{t_{0}} | ) \leqslant C+ |
     \theta_{t_{0}} | \]
  and
  \begin{equation}
    \| \theta \|_{\gamma , [ t_{0} ,t_{0} +S ]} \leqslant ( 1+ \| G \|_{\gamma
    , \nu , \psi} ) | \theta_{t_{0}} | +C ( \| G \|_{\gamma , \nu , \psi} +1 )
    . \label{eq:estimation-holder-norm-colutions}
  \end{equation}
  From this we deduce
  \[ | \theta_{t} | \leqslant | \theta_{t} - \theta_{t_{0}} | + |
     \theta_{t_{0}} | \lesssim_{T} \| \theta \|_{\gamma ,J} + | \theta_{t_{0}}
     | \lesssim ( 1+ \| G \|_{\gamma , \nu , \psi} ) | \theta_{t_{0}} | + \| G
     \|_{\gamma , \nu , \psi} +1 \]
  and then
  \[ \| \theta \|_{\infty ,J} \lesssim ( 1+ \| G \|_{\gamma , \nu , \psi} ) |
     \theta_{t_{0}} | + \| G \|_{\gamma , \nu , \psi} +1. \]
  Now let $n$ be such that $S=T/n$ and $1/2 \leqslant C  \| G \| T^{\nu \gamma}
  n^{- \nu \gamma} \leqslant 1$ hence $n \geqslant T ( C \| G \| )^{1/ \nu
  \gamma}$ and we have for $J_{i} = [ i T/n, ( i+1 ) T/n ]$
  \[ \| \theta \|_{\infty ,J_{i}} \lesssim ( 1+ \| G \|_{\gamma , \nu , \psi}
     ) \| \theta \|_{\infty ,J_{i-1}} + \| G \|_{\gamma , \nu , \psi} +1 \]
  and
  \[ \| \theta \|_{\infty ,J_{0}} \lesssim ( 1+ \| G \|_{\gamma , \nu , \psi}
     ) | \theta_{0} | + \| G \|_{\gamma , \nu , \psi} +1. \]
  Hence
  \[ \| \theta \|_{\infty ,J_{i}} \lesssim C^{i} ( 1+ \| G \|_{\gamma , \nu ,
     \psi} )^{i} ( ( 1+ \| G \|_{\gamma , \nu , \psi} ) | \theta_{0} | + \| G
     \|_{\gamma , \nu , \psi} +1 ) \]
  and finally
  \[ \| \theta \|_{\infty , [ 0,T ]} \leqslant K_{1} ( 1+ \| G \|_{\gamma ,
     \nu , \psi} )^{K \| G \|_{\gamma , \nu , \psi}^{1/ \nu \gamma}} ( |
     \theta_{t_{0}} | +1 ), \]
  where $K_{1}$ and $K_{2}$ are two universal constants depending on $\nu ,
  \gamma , \psi  $ and $T$. From the
  equation~{\eqref{eq:estimation-holder-norm-colutions}}, we can deduce,
  with the same induction argument, that
  \[ \| \theta \|_{\gamma , [ 0,T ]} \leqslant | \theta_{0} | +C \| G \|^{1/
     \nu \gamma} T ( ( 1+ \| G \|_{\gamma , \nu , \psi} ) \| \theta \|_{\infty
     , [ 0,T ]} + \| G \|_{\gamma , \nu , \psi} +1 ) \]
  and the result follows.
\end{proof}

\begin{remark}
  \label{remark:localization_principle}The bounds on the solutions of the
  differential equation allows us to get rid of the $\nu$-weight $\psi$.
  Indeed, we have, for a solution $\theta$ of the non-linear Young equation we
  have
  \[ \| \theta \|_{\infty , [ 0,T ]} \leqslant K_{1} ( 1+ \| G \|_{\gamma ,
     \nu , \psi} )^{K_{2} \| G \|_{\gamma , \nu , \psi}^{1/ \nu \gamma}} ( |
     \theta_{t_{0}} | +1 ) . \]
  Then for $R>0$, and $\theta_{0} \in B ( 0,R )$ $\| \theta \|_{\infty , [ 0,T
  ]} \leqslant K_{1} ( 1+ \| G \|_{\gamma , \nu , \psi} )^{K \| G \|_{\gamma ,
  \nu , \psi}^{1/ \nu \gamma}} ( R+1 ) .$ Hence, it is enough to consider the
  localised norm of $G$
  \[ \| G \|^{R}_{\gamma , \nu} = \sup_{t \neq s \in I} \sup_{x \neq y \in B (
     0,R^{G,R} )} \frac{| G_{s,t} ( x ) -G_{s,t} ( y ) |}{| x-y |^{\nu} | t-s
     |^{\gamma}} + \sup_{s \neq t} \frac{| G_{s,t} ( 0 ) |}{| t-s |^{\gamma}},
  \]
  where
  \[ R^{G,R} =K_{1} ( 1+ \| G \|_{\gamma , \nu , \psi} )^{K \| G \|_{\gamma ,
     \nu , \psi}^{1/ \nu \gamma}} ( R+1 ) . \]
\end{remark}

From now we will consider only bounded $G$, namely $G \in \CC^{\gamma , \nu}$,
and we will extend the results to $\CC^{\gamma , \nu , \psi}$ thanks to the
previous remark.

\subsection{Uniqueness of solutions}\label{subsec:uniqueness}

\subsubsection{Comparison Principle}\label{subsec:comparison_principle}

From now, thanks to remark \ref{remark:localization_principle} we can restrict
the study of the properties of the solutions, their uniqueness and their
regularity with respect to the parameters for bounded $G$. Hence, we define the
 space
$ \CC^{\gamma ,n+ \nu}_{b} = \{ G \in \CC^{\gamma , \nu} : \| G
   \|^{b}_{\gamma ,n+ \nu} <+ \infty
   \} $
   with $$
   \| G
   \|^{b}_{\gamma ,n+ \nu} \assign \llbracket D^{n} G \rrbracket_{\gamma ,
   \nu} + \sum_{k=0}^{n} \sup_{s \neq t \in [ 0,T ]} \sup_{x \in
   \mathbbm{R}^{d}} | D^{k} G_{s,t} ( x ) | / | t-s |^{\gamma}.
   $$
As there will be no ambiguity in the following, we will usually avoid to
mention explicitly the $b$ in the norm on that space. Those spaces are nicer
than the whole space $\CC^{\gamma , \nu}$ as there are natural embeddings:

\begin{lemma}
  \label{lemma:embedding}Let $0< \gamma \leqslant 1$, $0 \leqslant
  \nu' \leqslant \nu$ and $G \in \CC^{\gamma , \nu}_{b}$ then
  $ \| G \|^{b}_{\gamma , \nu'} \lesssim \| G \|_{\gamma , \nu}^{b} $.
  \end{lemma}

\begin{proof}
  Let $x,y \in \mathbbm{R}^{d}$ \ $s,t \in [ 0,T ]$. For $0 \leqslant \nu
  \leqslant \nu' \leqslant 1$, we have
  \begin{eqnarray*}
    \frac{| G_{s,t} ( x ) -G_{s,t} ( y ) |}{| x-y |^{\nu'}} & \leqslant &
    \left( \frac{| G_{s,t} ( x ) -G_{s,t} ( y ) |}{| x-y |^{\nu}}
    \right)^{\nu' / \nu} | G_{s,t} ( x ) -G_{s,t} ( y ) |^{1- \nu' / \nu}\\
    & \leqslant & 2^{1- \nu'} | t-s |^{\gamma} \llbracket G
    \rrbracket_{\gamma , \nu}^{\nu' / \nu} ( \| G \|^{b}_{\gamma , \nu} )^{1-
    \nu' / \nu}\\
    & \lesssim & 2^{1- \nu' / \nu} | t-s |^{\gamma} \| G \|_{\gamma ,
    \nu}^{b}
  \end{eqnarray*}
  and the following bound holds
  \[ \| G \|_{\gamma , \nu'}^{b} = ( 1+2^{1- \nu' / \nu} ) \| G \|_{\gamma ,
     \nu}^{b} . \]
  Furthermore, we also have
  \begin{eqnarray*}
    | G_{s,t} ( x ) -G_{s,t} ( y ) | & \leqslant & \int_{0}^{1} \mathd r  | D
    G_{s,t} ( r ( x-y ) +y ) |   | x-y |\\
    & \leqslant & \| D G \|_{\gamma ,1+ \nu}^{b} | x-y | | t-s |^{\gamma}
  \end{eqnarray*}
  and
  \[ \| G \|^{b}_{\gamma ,1} \leqslant 2 \| D G \|_{\gamma ,1+ \nu}^{b} . \]
  The general result follows by an easy induction.
\end{proof}

\begin{remark}
  \label{remark:a_priori_bounds}These embeddings allows us to state a result
  for the existence of the solutions when $G \in \CC^{\gamma ,1}_{b}$ with
  $\gamma > \frac{1}{2}$. Indeed, as for all $\nu <1$, $G \in \CC^{\gamma ,
  \nu}_{b}$ and $\| G \|_{\gamma , \nu}^{b} \lesssim \| G \|^{b}_{\gamma ,1}$,
  there exists a solution $\theta$ and the non-linear Young differential
  equation. Furthermore, for all $\nu <1$, there exists a constant $K_{2}$
  such that
  \[ \| \theta \|_{\infty} \lesssim ( 1+ \| G \|_{\gamma , \nu}^{b} )^{K_{2}
     ( \| G \|^{b}_{\gamma , \nu} )^{1/ \nu \gamma}} ( | \theta_{0} | +1 )  
     \lesssim ( 1+ \| G \|_{\gamma ,1}^{b} )^{K_{2} ( \| G \|^{b}_{\gamma ,1}
     )^{1/ \nu \gamma}} ( | \theta_{0} | +1 ) . \]
  In fact, a deeper look at the proof of Theorem
  \ref{theorem:existence_young_solution}, allows us to get rid of the $\nu$,
  and state that there exists a constant $K$ depending on $T$ and $\gamma$
  such that
  \[ \| \theta \|_{\infty} \lesssim ( 1+ \| G \|_{\gamma ,1}^{b} )^{K ( \| G
     \|^{b}_{\gamma ,1} )^{1/ \gamma}} ( | \theta_{0} | +1 )  \]
  and a similar bound holds for $\| \theta \|_{\gamma}$.
\end{remark}

In order to study the properties of the solutions of the non-linear Young
differential equation, we intend to compare two solutions $\theta^{1}$ and
$\theta^{2}$. In the classical case (when $G$ is differentiable in time), we
would have
\begin{eqnarray*}
  \theta^{1}_{t} - \theta^{2}_{t} & = & ( \theta^{1}_{0} - \theta^{2}_{0} ) +
  \int_{0}^{t}\big( G'_{u} ( \theta^{1}_{u} ) -G'_{u} ( \theta^{2}_{u} )\big) \mathd u\\
  & = & ( \theta^{1}_{0} - \theta^{2}_{0} ) + \int_{0}^{t} \big( G'_{u} (
  \theta^{1}_{u} - \theta^{2}_{u} + \theta^{2}_{u} ) -G'_{u} ( \theta^{2}_{u}
  )\big) \mathd u\\
  & = & ( \theta^{1}_{0} - \theta^{2}_{0} ) + \int_{0}^{t} (
  \tau_{\theta^{2}} G )'_{u} ( \theta^{1}_{u} - \theta^{2}_{u} ) \mathd u-
  \int_{0}^{t} G'_{u} ( \theta^{2}_{u} ) \mathd u,
\end{eqnarray*}
where $( \tau_{\theta^{2}} G )_{t} ( x ) = \int_{0}^{t} G'_{u} (
\theta^{2}_{u} +x ) \mathd u$. Hence $\theta^{1} - \theta^{2}$ solve a
differential equation, but with a translated and averaged function $(
\tau_{\theta^{2}} G )$ and a second member. In order to prove some properties
on the solutions, we then have to study this differential equation.

In the case of the Young differential equation, this strategy will be very
profitable, we have to define the averaged translation and to study some
of its properties. Hence, we define the natural action of the additive group
of $\CC^{\rho}$ paths on the integrated vector fields $C \in \CC^{\gamma ,
\nu}_{b}$.

\begin{definition}
  Let $\gamma , \nu , \rho \in [ 0,1 ]$ such that $\gamma + \rho \nu >1$, $G
  \in \CC_{b}^{\gamma , \nu}$ and $f \in \CC^{\rho}$. We define the average
  translation of $G$ by $f$, and we write $\tau_{f} G$ the following quantity
  \[ \tau_{f} G: ( t,x ) \rightarrow \int_{0}^{t} G_{\mathd u} ( f_{u} +x ) .
  \]
\end{definition}

Due to the requirements of Young integration, the estimations for the
translated integrated vector field $\tau_{f} G$ show a loss of regularity
quantified by the next lemma.

\begin{lemma}
  \label{lemma:translation}For $\gamma + \nu \rho >1$ and $\gamma + \eta \rho
  >1$, $f \in \CC^{\rho}$ and $G \in \CC_{b}^{\gamma , \nu + \eta}$ we have
  $\tau_{f} G \in \CC^{\gamma , \nu}$ and
  \[ \| \tau_{f} G \|_{\gamma , \nu} \lesssim \| G \|_{\gamma , \nu + \eta} (
     1+ \| f \|^{\eta}_{\rho} ) .\]
\end{lemma}

\begin{proof}
  Suppose first that $\eta + \nu \leqslant 1$. Let $x,y \in V$ and define
  $\tilde{G} ( z ) =G ( x+z ) -G ( y+z )$. There is two bounds for the
  increments of $\tilde{G}$:
  \[ | \tilde{G}_{s,t} ( z_{1} ) - \tilde{G}_{s,t} ( z_{2} ) | \lesssim
     \llbracket G \rrbracket_{\gamma , \nu + \eta} | t-s |^{\gamma} | x-y
     |^{\nu + \eta} \]
  and
  \[ | \tilde{G}_{s,t} ( z_{1} ) - \tilde{G}_{s,t} ( z_{2} ) | \lesssim
     \llbracket G \rrbracket_{\gamma , \nu + \eta} | t-s |^{\gamma} | z_{1}
     -z_{2} |^{\nu + \eta} . \]
  Hence, by interpolating these two inequalities
  \[ | \tilde{G}_{s,t} ( z_{1} ) - \tilde{G}_{s,t} ( z_{2} ) | \lesssim
     \llbracket G \rrbracket_{\gamma , \nu + \eta} | t-s |^{\gamma} | x-y
     |^{\nu} | z_{1} -z_{2} |^{\eta} . \]
  When $2 \geqslant \eta + \nu >1$, we have
  \begin{eqnarray*}
    | \tilde{G}_{s,t} ( z_{1} ) - \tilde{G}_{s,t} ( z_{2} ) | & = & \left\lvert
    \int_{0}^{1} \mathd r  \{ D G_{s,t} ( x(r) +z_{1}  ) -D G_{s,t} ( x(r) +z_{2} ) \} .(x-y) \right\lvert\\
    & \lesssim & \llbracket G \rrbracket_{\gamma , \nu + \eta} | t-s
    |^{\gamma} | x-y | | z_{1} -z_{2} |^{\nu + \eta -1},
  \end{eqnarray*}
  where $x(r) = y+r(x-y)$,
  
  \[ | \tilde{G}_{s,t} ( z_{1} ) - \tilde{G}_{s,t} ( z_{2} ) | \lesssim
     \llbracket G \rrbracket_{\gamma , \nu + \eta} | t-s |^{\gamma} | x-y
     |^{\nu + \eta -1} | z_{1} -z_{2} | \]
  and again
  \[ | \tilde{G}_{s,t} ( z_{1} ) - \tilde{G}_{s,t} ( z_{2} ) | \lesssim
     \llbracket G \rrbracket_{\gamma , \nu + \eta} | t-s |^{\gamma} | x-y
     |^{\nu} | z_{1} -z_{2} |^{\eta} . \]
  In these two cases, we have
  \[ \llbracket \tilde{G} \rrbracket_{\gamma , \eta} \lesssim \llbracket G
     \rrbracket_{\gamma , \nu + \eta} | t-s |^{\gamma} | x-y |^{\nu} . \]
  Hence
  \begin{eqnarray*}
    \tau_{f} G_{s,t} ( x ) - \tau_{f} G_{s,t} ( y ) & = & \int_{s}^{t}
    G_{\mathd u} ( f_{u} +x ) -G_{\mathd u} ( f_{u} +y )\\
    & \leqslant & \int_{s}^{t} \tilde{G}_{\mathd u} ( f_{u} ) -
    \tilde{G}_{s,t} ( f_{s} ) + \tilde{G}_{s,t} ( f_{s} ).
  \end{eqnarray*}
  Hence,
  \[ | \tau_{f} G_{s,t} ( x ) - \tau_{f} G_{s,t} ( y ) | \lesssim \llbracket
     \tilde{G} \rrbracket_{\gamma , \eta} ( \| f \|^{\eta}_{\gamma} +1 ) + |
     \tilde{G}_{s,t} ( f_{s} ) | \]
  and as $| \tilde{G}_{s,t} ( f_{s} ) | \leqslant 2 | t-s |^{\gamma} \| G
  \|_{\gamma , \nu + \eta}$,
  \[ \llbracket \tau_{f} G \rrbracket_{s,t} \leqslant \| G \|_{\gamma , \nu +
     \eta} ( \| f \|^{\eta}_{\gamma} +1 ) . \]
  Furthermore
  \[ | ( \tau_{f} G )_{s,t} ( 0 ) | \leqslant \left\lvert \int_{s}^{t} G_{\mathd u}
     ( f_{u} ) -G_{s,t} ( f_{s} ) \right\lvert + | G_{s,t} ( f_{s} ) | \leqslant \|
     G \|_{\gamma , \eta} | t-s |^{\gamma} ( \| f \|_{\gamma}^{\eta} +1 ) \]
  and by the embedding of Lemma \ref{lemma:embedding}, the result follows.
\end{proof}

The averaged translation is a suitable tool to control the difference of two
non-linear Young integrals, as soon as we have enough regularity to estimate
the integral. The following lemma states the estimation for generic functions.

\begin{lemma}
  Let $\gamma , \nu , \nu' , \rho \in [ 0,1 ]$ such that $\gamma + \rho \nu
  >1$ and $\gamma + \rho \nu' >1$. Let $f^{1} ,f^{2} \in \CC^{\rho}$, $G \in
  \CC_{b}^{\gamma , \nu}$, and suppose that $\tau_{f^{2}} G \in
  \CC_{b}^{\gamma , \nu'}$. Then
  \[ \left\| \int_{0}^{.} G_{\mathd u} ( f^{1}_{u} ) -G_{\mathd u} ( f^{2}_{u}
     ) \right\|_{\infty , [ 0,T ]} \lesssim \| \tau_{f^{2}} G \|_{\gamma ,
     \nu'} T^{\gamma} ( \llbracket f^{1} -f^{2} \rrbracket^{\nu'}_{\rho} + \|
     f^{1} -f^{2} \|^{\nu'}_{\infty} ) \]
  and
  \[ \left\llbracket \int_{0}^{.} G_{\mathd u} ( f^{1}_{u} ) -G_{\mathd u} (
     f^{2}_{u} ) \right\rrbracket_{\gamma , [ 0,T ]} \lesssim \| \tau_{f^{2}}
     G \|_{\gamma , \nu'} ( \llbracket f^{1} -f^{2} \rrbracket^{\nu'}_{\gamma}
     + \| f^{1} -f^{2} \|^{\nu'}_{\infty} ) . \]
  Furthermore when $1 \geqslant \eta >0$ such that $\rho \eta + \gamma >1$ and
  $G \in \CC^{\gamma , \nu' + \eta}$
  \[ \left\| \int_{0}^{.} G_{\mathd u} ( f^{1}_{u} ) -G_{\mathd u} ( f^{2}_{u}
     ) \right\|_{\infty , [ 0,T ]} \lesssim T^{\gamma} \| G \|_{\gamma , \nu'
     + \eta} ( 1+ \| f^{2} \|^{\eta}_{\rho} ) ( \llbracket f^{1} -f^{2}
     \rrbracket^{\nu'}_{\rho} + \| f^{1} -f^{2} \|^{\nu'}_{\infty} ) . \]
\end{lemma}

\begin{proof}
  It is a direct application of the definition of the averaged translation.
  Let $s,t \in [ 0,T ]$, by definition we have $\int_{s}^{t} G_{\mathd u} (
  f^{2}_{u} ) = \tau_{f^{2}} G_{s,t} ( 0 )$. Hence
  \begin{eqnarray*}
    \left\lvert \int_{s}^{t} G_{\mathd u} ( f^{1}_{u} ) - \int_{s}^{t} G_{\mathd u}
    ( f^{2}_{u} ) \right\lvert & \leqslant & \left\lvert \int_{s}^{t} \tau_{f^{2}}
    G_{\mathd u} ( f^{1}_{u} -f^{2}_{u} ) - \tau_{f^{2}} G_{s,t} ( f^{1}_{s}
    -f^{2}_{s} ) \right\lvert\\
    &  & + | \tau_{f^{2}} G_{s,t} ( f^{1}_{s} -f^{2}_{s} ) - \tau_{f^{2}}
    G_{s,t} ( 0 ) |\\
    & \lesssim & \llbracket \tau_{f^{2}} G \rrbracket_{\gamma , \nu'} | t-s
    |^{\gamma} ( \llbracket f^{1} -f^{2} \rrbracket_{\rho} | t-s |^{\nu'
    \gamma} + | f_{s}^{1} -f_{s}^{2} |^{\nu'} ).
  \end{eqnarray*}
  Hence,
  \[ \left\| \int_{0}^{.} G_{\mathd u} ( f^{1}_{u} ) -G_{\mathd u} ( f^{2}_{u}
     ) \right\|_{\infty , [ 0,T ]} \lesssim \llbracket \tau_{f^{2}} G
     \rrbracket_{\gamma , \nu'} T^{\gamma} [ \| f^{1} -f^{2}
     \rrbracket^{\nu'}_{\rho} + \| f^{1} -f^{2} \|^{\nu'}_{\infty} ] \]
  and
  \[ \left\llbracket \int_{0}^{.} G_{\mathd u} ( f^{1}_{u} ) -G_{\mathd u} (
     f^{2}_{u} ) \right\rrbracket_{\gamma , [ 0,T ]} \lesssim \llbracket
     \tau_{f^{2}} G \rrbracket_{\gamma , \nu'} ( \llbracket f^{1} -f^{2}
     \rrbracket^{\nu'}_{\rho} + \| f^{1} -f^{2} \|^{\nu'}_{\infty} ) . \]
  For the second part of the lemma, we use the bound of
  Lemma~\ref{lemma:translation}.
\end{proof}

We are now ready to prove a comparison principle between two solutions. In order
to keep a high degree of generality, we do not use the estimation of the
Lemma~\ref{lemma:translation}, but prefer to state a general assumption for
the regularity of the averaged translation of the first vector field.

\begin{theorem}
  \label{theorem:comparison}Let $\gamma > \frac{1}{2}$, $\nu \in [ 0,1 ]$ such
  that $\gamma ( 1+ \nu ) >1$. Let $G^{1} ,G^{2} \in \CC_{b}^{\gamma , \nu}$,
  $\theta^{1}$ (respectively $\theta^{2}$) be a solution of the nonlinear Young
  differential equation driven by $G^{1}$ (respectively by $G^{2}$). Suppose
  that $\tau_{\theta^{2}} G^{1} \in \CC^{\gamma ,1}$. Then
  \[ \| \theta^{1} - \theta^{2} \|_{\infty , [ 0,T ]} \leqslant c_{1}
     e^{c_{2} \| \tau_{\theta^{2}} G^{1} \|_{\gamma ,1}^{1/ \gamma}} ( \|
     \theta^{2} \|^{\nu}_{\gamma} +1 ) ( | \theta^{1}_{0} - \theta^{2}_{0} | +
     \| G^{1} -G^{2} \|_{\gamma , \nu} ) . \]
\end{theorem}

\begin{proof}
  Let $t_{0} \in [ 0,T ]$, $S>0$ and define $J= [ t_{0} , ( t_{0} +S ) \wedge
  1 ]$. For $s \leqslant t \in J$ we have
  \begin{eqnarray*}
    \delta ( \theta^{1} - \theta^{2} )_{s,t} & = & \int_{s}^{t}
    \tau_{\theta^{2}} G^{1}_{\mathd u} ( \theta^{1}_{u} - \theta^{2}_{u} ) -
    \tau_{\theta^{2}} G^{1}_{s,t} ( \theta^{1}_{s} - \theta^{2}_{s} )\\
    &  & + \tau_{\theta^{2}} G^{1}_{s,t} ( \theta^{1}_{s} - \theta^{2}_{s} )
    - \tau_{\theta^{2}} G^{1}_{s,t} ( 0 )\\
    &  & + \int_{s}^{t} ( G^{1} -G^{2} )_{\mathd u} ( \theta^{2}_{u} ) - (
    G^{1} -G^{2} )_{s,t} ( \theta^{2}_{s} )\\
    &  & + ( G^{1} -G^{2} )_{s,t} ( \theta^{2}_{s} ) - ( G^{1} -G^{2} )_{s,t}
    ( 0 )\\
    &  & + ( G^{1} -G^{2} )_{s,t} ( 0 ).
  \end{eqnarray*}
  Hence,
  \begin{eqnarray*}
    | \delta ( \theta^{1} - \theta^{2} )_{s,t} | & \lesssim & \|
    \tau_{\theta^{2}} G^{1} \|_{\gamma ,1} | t-s |^{\gamma} ( S^{\gamma}
    \llbracket \theta^{1} - \theta^{2} \rrbracket_{\gamma} + | \theta^{1}_{s}
    - \theta^{2}_{s} | )\\
    &  & + \| G^{1} -G^{2} \|_{\gamma , \nu} | t-s |^{\gamma} ( S^{\gamma
    \nu} \llbracket \theta^{2} \rrbracket_{\gamma}^{\nu} + | \theta^{2}_{s}
    |^{\nu} +1 ) .
  \end{eqnarray*}
  When $C_{1}$ is the universal constant in the previous inequality and for
  $S$ small enough such that $\frac{1}{4} \leqslant C_{1} \| \tau_{\theta^{2}}
  G^{1} \|_{\gamma ,1} S^{\gamma} \leqslant \frac{1}{2}$, there exists another constant $C_{2}$ such that
  \[ | \delta ( \theta^{1} - \theta^{2} )_{s,t} | \leqslant \frac{1}{2} | t-s
     |^{\gamma} ( \llbracket \theta^{1} - \theta^{2} \rrbracket_{\gamma} +S^{-
     \gamma} | \theta^{1}_{s} - \theta^{2}_{s} | ) +C_{2} \| G^{1} -G^{2}
     \|_{\gamma , \nu} | t-s |^{\gamma} ( S^{\gamma \nu} \| \theta^{2}
     \|^{\nu}_{\gamma} +1 ) . \]
  Hence
  \[ \llbracket \theta^{2} - \theta^{2} \rrbracket_{\gamma} \leqslant S^{-
     \gamma} \| \theta^{1} - \theta^{2} \|_{\infty ,J} +C_{3} \| G^{1} -G^{2}
     \|_{\gamma , \nu} ( \| \theta^{2} \|^{\nu}_{\gamma} +1 ) \]
  and
  \begin{eqnarray*}
    \| \theta^{1} - \theta^{2} \|_{\infty ,J} & \leqslant & \frac{1}{2} ( \|
    \theta^{1} - \theta^{2} \|_{\infty ,J} + | \theta^{1}_{s} - \theta^{2}_{s}
    | ) +C_{4} \| G^{1} -G^{2} \|_{\gamma , \nu} S^{\gamma} ( \| \theta^{2}
    \|^{\nu}_{\gamma} +1 ) .
  \end{eqnarray*}
  Finally
  \[ \| \theta^{1} - \theta^{2} \|_{\infty ,J} \leqslant 2 | \theta^{1}_{s} -
     \theta^{2}_{s} | +C_{5} \| G^{1} -G^{2} \|_{\gamma , \nu} S^{\gamma} ( \|
     \theta^{2} \|^{\nu}_{\gamma} +1 ) . \]
  By the same gluing argument as in Therorem
  \ref{theorem:existence_young_solution}, we have
  \[ \| \theta^{1} - \theta^{2} \|_{\infty} \lesssim 2^{1/S} ( |
     \theta^{1}_{0} - \theta^{2}_{0} | +C_{5} \| G^{1} -G^{2} \|_{\gamma ,
     \nu} S^{\gamma} ( \| \theta^{2} \|^{\nu}_{\gamma} +1 ) ) . \]
  Remind that $\frac{1}{4} \leqslant C_{1} \| \tau_{\theta^{2}} G^{1}
  \|_{\gamma ,1} S^{\gamma} \leqslant \frac{1}{2}$, and there exists two
  universal constants (depending on $\gamma , \nu ,T$) $c_{1}$ and $c_{2}$
  such that
  \[ \| \theta^{1} - \theta^{2} \|_{\infty} \leqslant c_{1} e^{c_{2} \|
     \tau_{\theta^{2}} G^{1} \|_{\gamma ,1}^{1/ \gamma}} ( | \theta^{1}_{0} -
     \theta^{2}_{0} | + \| G^{1} -G^{2} \|_{\gamma , \nu} ( \| \theta^{2}
     \|^{\nu}_{\gamma} +1 ) ), \]
  which ends the proof.
\end{proof}

\subsubsection{Uniqueness of solutions}

We prove here the uniqueness of the solutions when $G$ is regular enough using
the comparison principle given in Theorem~\ref{theorem:comparison}. In order
to use the comparison principle, we ask that the vector field is regular
enough (in space) and as stated in Lemma~\ref{lemma:translation}, we are able
to estimate the averaged translation if we accept a additional loss of space
regularity. Furthermore, as we will use uniqueness results in different
contexts, especially in situation when we will have a priori regularity
properties for the solutions, we will give a pretty general theorem of
uniqueness, Theorem \ref{theorem:uniqueness_general_case}, and specialise it in
the two following corollaries.

\begin{theorem}
  \label{theorem:uniqueness_general_case}Let $\gamma >1/2$, $\nu \in
  [ 0,1 ]$ such that $\gamma ( 1+ \nu ) >1$ and $G \in \CC^{\gamma ,
  \nu}_{b}$. Suppose that there exists a sequence $( G^{\varepsilon}
  )_{\varepsilon} \in \CC^{\gamma , \nu}_{b}$ such that
  \begin{enumerateroman}
    \item For all $\gamma' < \gamma$ and all, $\| G-G^{\varepsilon}
    \|_{\gamma' , \nu} \rightarrow 0$.
    
    \item For all $\varepsilon >0$ and $\| G^{\varepsilon} \|_{\gamma , \nu}
    \leqslant \| G \|_{\gamma , \nu}$.
    
    \item For all $\varepsilon >0$ and all $x \in \mathbbm{R}^{d}$ there
    exists a unique solution $\theta^{\varepsilon}$ for the equation
    \[ \theta^{\varepsilon}_{t} ( x ) =x+ \int_{0}^{t}
       G^{\varepsilon}_{\mathd u} ( \theta^{\varepsilon}_{u} ( x ) ) \mathd u.
    \]
    \item For all $\varepsilon >0$, $\tau_{\theta^{\varepsilon}} G \in
    \CC^{\gamma ,1}_{b}$ and
    $ \sup_{\varepsilon >0} \| \tau_{\theta^{\varepsilon}} G \|_{\gamma ,1}
       <+ \infty $. 

    Then the solution of the nonlinear Young equation driven by $G$ is unique.
  \end{enumerateroman}
\end{theorem}

\begin{proof}
  This theorem is a direct consequence of the comparison principle of
  Theorem \ref{theorem:comparison}. Let $x \in \mathbbm{R}^{d}$ be an initial
  condition, and let $\theta$ be a solution of the nonlinear Young differential
  equation with initial condition $x$. Furthermore let $G^{\varepsilon}$ and
  $\theta^{\varepsilon}$ be as in the hypothesis of the theorem. Take
  $\frac{1}{2} < \gamma' < \gamma$ such that $\gamma ( 1+ \nu ) >1$. Remark
  that we can apply the comparison principle to $\theta$ and
  $\theta^{\varepsilon}$ with $\gamma'$ instead of $\gamma$, and as $\|
  G^{\varepsilon} \|_{\gamma' , \nu} \leqslant \| G \|_{\gamma' , \nu}
  \lesssim \| G \|_{\gamma , \nu}$, we have
  \begin{eqnarray*}
    \| \theta^{\varepsilon} \|_{\gamma'} & \lesssim & \| G^{\varepsilon}
    \|_{\gamma' , \nu}^{1/ \nu \gamma'} ( 1+ \| G^{\varepsilon} \|_{\gamma' ,
    \nu} )^{\tilde{K}_{2} \| G^{\varepsilon} \|_{\gamma' , \nu}^{1/ \nu
    \gamma'}} ( 1+ | x | ) \\
    & \lesssim & \| G \|_{\gamma , \nu}^{1/ \nu \gamma'} ( 1+ \| G \|_{\gamma
    , \nu} )^{\tilde{K}_{2} \| G \|_{\gamma , \nu}^{1/ \nu \gamma'}} ( 1+ | x
    | ) ,
  \end{eqnarray*}
  but also
  \begin{eqnarray*}
    \| \theta ( x ) - \theta^{\varepsilon} ( x ) \|_{\infty} & \lesssim &
    e^{c_{2} \| \tau_{\theta^{\varepsilon}} G \|_{\gamma' ,1}^{1/ \gamma}} \|
    G^{1} -G^{\varepsilon} \|_{\gamma' , \nu} ( \| \theta^{\varepsilon}
    \|^{\nu}_{\gamma'} +1 )\\
    & \lesssim & \| G \|_{\gamma , \nu}^{1/ \nu \gamma'} ( 1+ \| G \|_{\gamma
    , \nu} )^{\tilde{K}_{2} \| G \|_{\gamma , \nu}^{1/ \nu \gamma'}} ( 1+ | x
    | ) e^{c_{2} \| \tau_{\theta^{\varepsilon}} G \|_{\gamma ,1}^{1/ \gamma}}
    \| G^{1} -G^{\varepsilon} \|_{\gamma' , \nu} .
  \end{eqnarray*}
  As $\sup_{\varepsilon >0} \| \tau_{\theta^{\varepsilon}} G \|_{\gamma ,1} <+
  \infty$, we have $\| \theta ( x ) - \theta^{\varepsilon} ( x ) \|_{\infty}
  \rightarrow_{\varepsilon \rightarrow 0} 0$. As $\theta^{\varepsilon}$ is
  unique, and since this convergence holds true for every function $\theta ( x
  )$ solution of the equation, the solution is unique.
\end{proof}

We can now use the averaged translation operator to establish uniqueness in
the case where we have a priori informations of the regularities of the
solutions $\theta^{\varepsilon}$.

\begin{corollary}
  \label{corollary:uniqueness_a_priori_case}Let $\gamma >1/2 \nocomma$,
  $\delta >0$ and $G \in \CC^{\gamma ,1+ \eta}_{b}$. Suppose that there exists
  a sequence $( G^{\varepsilon} )_{\varepsilon} \in \CC^{\gamma ,1+ \eta}_{b}$
  such that
  \begin{enumerateroman}
    \item For all $\gamma' < \gamma$ and all, $\| G-G^{\varepsilon}
    \|_{\gamma' ,1} \rightarrow 0$.
    
    \item For all $\varepsilon >0$ and $\| G^{\varepsilon} \|_{\gamma ,1+
    \eta} \leqslant \| G \|_{\gamma ,1+ \eta}$.
    
    \item For all $\varepsilon >0$ and all $x \in \mathbbm{R}^{d}$ there
    exists a unique solution $\theta^{\varepsilon}$ for the equation
    \[ \theta^{\varepsilon}_{t} ( x ) =x+ \int_{0}^{t}
       G^{\varepsilon}_{\mathd u} ( \theta^{\varepsilon}_{u} ( x ) ) \mathd u.
    \]
    \item There exists $\rho >0$ such that $\eta \rho + \gamma >0$ and for
    which for all $\varepsilon >0$, $\theta^{\varepsilon} ( x ) \in
    \CC^{\rho}$ and $\sup_{\varepsilon >0} \| \theta^{\varepsilon} ( x )
    \|_{\rho} <+ \infty$.
  \end{enumerateroman}
  Then solution of the non-linear Young equation driven by $G$ is unique.
  
  Furthermore, when the function $x \rightarrow \sup_{\varepsilon >0} \|
  \theta^{\varepsilon} ( x ) \|_{\rho}$ is locally bounded in time, the flow
  $( t,x ) \rightarrow \theta_{t} ( x )$ of the equation is locally Lipschitz
  continuous in space, uniformly in time and
  \[ \| \theta^{} ( x ) - \theta ( y ) \|_{\infty} \lesssim \exp ( C ( 1+
     \log ( 1+ \| G \|_{\gamma ,1+ \eta} ) + ( \sup_{\varepsilon >0} \|
     \theta^{\varepsilon} ( y ) \|_{\rho} )^{\eta} ) \| G \|_{\gamma ,1+
     \delta}^{1/ \gamma} ) | x-y |^{} ( | y | +1 ) . \]
\end{corollary}

\begin{proof}
  Condition \tmtextit{i.}, \tmtextit{ii.} and \tmtextit{iii.} are the same of
  those of Theorem $\ref{theorem:uniqueness_general_case}$. We only have to
  prove that the point \tmtextit{iv.} of Theorem
  $\ref{theorem:uniqueness_general_case}$ is satisfied. But thanks to Lemma
  \ref{lemma:translation}, we know that $\tau_{\theta^{\varepsilon}} G \in
  \CC^{\gamma ,1}_{b}$. Furthermore
  \[ \| \tau_{\theta^{\varepsilon}} G \|_{\gamma ,1} \lesssim \| G \|_{\gamma
     ,1+ \eta} ( \| \theta^{\varepsilon} \|_{\rho}^{\eta} +1 ) \lesssim \| G
     \|_{\gamma ,1+ \eta} ( ( \sup_{\varepsilon >0} \| \theta^{\varepsilon}
     \|_{\rho} )^{\eta} +1 )  \]
  and the uniqueness follows by Theorem \ref{theorem:uniqueness_general_case}.
  Furthermore, for $y \in \mathbbm{R}^{d}$ since $\| G^{\varepsilon}
  \|_{\gamma ,1+ \eta} \leqslant \| G \|_{\gamma ,1+ \eta}$, we have
  \[ \| \tau_{\theta^{\varepsilon} ( y )} G^{\varepsilon} \|_{\gamma ,1}
     \lesssim \| G^{\varepsilon} \|_{\gamma ,1+ \eta} ( ( \sup_{\varepsilon
     >0} \| \theta^{\varepsilon} \|_{\rho} )^{\eta} +1 ) \lesssim \| G
     \|_{\gamma ,1+ \eta} ( ( \sup_{\varepsilon >0} \| \theta^{\varepsilon}
     \|_{\rho} )^{\eta} +1 ) . \]
  Since $\sup_{\varepsilon >0} \| G^{\varepsilon} \|_{\gamma ,1} \lesssim \| G
  \|_{\gamma ,1+ \eta}$ we have, thanks to the a priori bounds for the
  solutions of Theorem \ref{theorem:existence_young_solution}, and Remark
  \ref{remark:a_priori_bounds}, that
  \[ \| \theta^{\varepsilon} ( y ) \|_{\gamma} \lesssim \| \theta \|_{\infty}
     \lesssim ( 1+ \| G \|_{\gamma ,1+ \eta} )^{K ( \| G \|_{\gamma ,1+ \eta}
     )^{1/ \gamma}} ( | y | +1 ) , \]
  which implies \
  \[ \| \theta^{\varepsilon} ( x ) - \theta^{\varepsilon} ( y ) \|_{\infty}
     \lesssim \exp ( C ( 1+ \log ( 1+ \| G \|_{\gamma ,1+ \eta} ) + (
     \sup_{\varepsilon >0} \| \theta^{\varepsilon} ( y ) \|_{\rho} )^{\eta} )
     \| G \|_{\gamma ,1+ \delta}^{1/ \gamma} ) | x-y |^{} ( | y | +1 ) , \]
  the conclusion then easily follows when we let $\varepsilon$ go to zero.
\end{proof}

\begin{remark}
  \label{remark:bound_approximated_flow}Suppose furthermore that for all $t
  \in [ 0,T ]$, $x \rightarrow \theta^{\varepsilon}_{t} ( x )$ is
  differentiable in space. Then
  \[ \sup_{\varepsilon >0}   | D \theta^{\varepsilon}_{t} ( x ) | \lesssim
     \exp ( C ( 1+ \log ( 1+ \| G \|_{\gamma ,1+ \eta} ) + ( \sup_{\varepsilon
     >0} \| \theta^{\varepsilon} ( x ) \|_{\rho} )^{\eta} ) \| G \|_{\gamma
     ,1+ \delta}^{1/ \gamma} )^{} ( | x | +1 ). \]
\end{remark}

Finally, we state the more general results of uniqueness, where all the needed
informations are the regularity of $G$.

\begin{corollary}
  Let $\gamma > \frac{1}{2}$, $\nu \in [ 0,1 ]$ such that $\gamma ( 1+ \nu )
  >0$ and suppose that $G \in \CC^{\gamma , \nu +1}_{b}$, then there exists a
  unique solution $\theta ( x )$ for the non-linear Young equation with
  initial condition $x$. Furthermore $\theta$ is locally Lipschitz continuous
  in space uniformly in time.
\end{corollary}

\begin{proof}
  We only have to check the conditions of Corollary
  \ref{corollary:uniqueness_a_priori_case} with $\nu = \eta$ and $\gamma =
  \rho$. Let $G^{\varepsilon} \in \CC^{\gamma ,1+ \nu}$ such that the time
  derivative $( G^{\varepsilon} )'$ exists and lies in $\CC^{\gamma ,1+
  \nu}_{b}$, and such that $\| G^{\varepsilon} \|_{\gamma ,1+ \nu} \leqslant
  \| G \|_{\gamma ,1+ \nu}$ and for all $\gamma' < \gamma$, $\|
  G-G^{\varepsilon} \|_{\gamma' ,1+ \nu} \rightarrow 0$. In that case,
  $\theta^{\varepsilon} ( x )$ is the solution of
  \[ \theta^{\varepsilon} ( x ) =x+ \int_{0}^{t} ( G^{\varepsilon}_{r} )' (
     \theta^{\varepsilon}_{r} ( x ) ) \mathd r=x+ \int_{0}^{t}
     G^{\varepsilon}_{\mathd r} ( \theta^{\varepsilon}_{r} ( x ) ) . \]
  As $G^{\varepsilon} \in \CC^{\gamma ,1+ \nu}_{b}$, $\theta^{\varepsilon}$ is
  unique and furthermore $\theta^{\varepsilon}$ is differentiable in space,
  and the differential is the solution of the following equation
  \[ D \theta^{\varepsilon} ( x ) = \tmop{id} + \int_{0}^{t} D (
     G^{\varepsilon}_{r} )' ( \theta^{\varepsilon}_{r} ( x ) ) \mathd r. \]
  Thanks to Remark \ref{remark:a_priori_bounds} we have
  \[ \| \theta^{\varepsilon} ( x ) \|_{\infty} + \| \theta^{\varepsilon} ( x )
     \|_{\gamma} \lesssim ( 1+ \| G \|_{\gamma ,1}^{b} )^{K ( \| G
     \|^{b}_{\gamma ,1} )^{1/ \gamma}} ( | x | +1 ) . \]
  Hence $x \rightarrow \sup_{\varepsilon} \| \theta^{\varepsilon} ( x )
  \|_{\gamma}$ is locally bounded in space. All the conditions of the
  Corollary \ref{corollary:uniqueness_a_priori_case} are fulfilled, and the
  result follows.
\end{proof}

\subsection{Localisation of unbounded vector fields}

In order to give a complete survey of the question, we need to go back to the
weighted spaces $\CC^{\gamma , \nu , \psi}$ and to state the  existence and uniqueness
theorems in that case.

Let $\gamma >1/2$, $\nu <1$ and $\gamma ( 1+ \nu ) >1$ and $\psi >0$ a
$\nu$-weight. Let $r>0$ and $r=K_{1} ( 1+ \| G \|_{\gamma ,1}^{b} )^{K_{2} (
\| G \|^{b}_{\gamma ,1} )^{1/ \gamma}} ( r+1 )$, where $K_{1}$ and $K_{2}$ are
define as in Theorem \ref{theorem:existence_young_solution} and depend on
$\psi$. As we intend to use the averaged translation operator, and since any
solution lies in balls of radius $R$, we need to localize $G$ on balls $B$ of
center $0$ and of radius $2R$. We then let $G_{|R} \in \CC^{\gamma , \nu}_{ b}
( [ 0,T ] ,B )$ the restriction of $G$ on $[ 0,T ] \times B$. We have of
course $\| G_{|R} \|^{b}_{\gamma , \nu} \lesssim \psi ( 2R ) \| \tilde{G}
\|_{\gamma , \nu , \psi}$. Furthermore, as all the arguments hold locally, as
we have done all the estimations for $x,y \in B ( 0,r )$ in the previous section.

When $\nu =1$, it is necessary to have the existence and a bound for the
solution in order to localise. As this holds only for $\nu <1$, the good
hypothesis is that there exists $\nu <1$, $\tilde{\psi}$ a $\nu$-weight such
that $G \in \CC^{\gamma , \nu , \tilde{\psi}} \cap \CC^{\gamma ,1, \psi}$. In
that case, we are again able to localise and to use the result of the
previous section. The following theorem holds:

\begin{theorem}
  Let $\gamma > \frac{1}{2}$, $1 \geqslant \nu >0$ with $\gamma ( 1+ \nu ) >1$
  and $\psi$ a weight. Let $1< \nu' \leqslant \nu$ with $\nu' <1$ such that
  $\gamma ( 1+ \nu' ) >1$, and $\psi'$ a $\nu'$-weight. Let $G \in \CC^{\gamma
  , \nu' , \psi'} \cap \CC^{\gamma , \nu , \psi}$. For all $x \in
  \mathbbm{R}^{d}$ there exists a solution $\theta ( x ) \in \CC^{\gamma} ( [
  0,T ] )$ to the equation
  \[ \theta_{t} ( x ) =x+ \int_{0}^{t} G_{\mathd u} ( \theta_{u} ) . \]
  Furthermore, there exists $K_{1}$ and $K_{2}$ two constants depending on
  $\gamma , \nu'$ and $\psi'$ such that
  \[ _{} \| \theta \|_{\gamma , [ 0,T ]} \leqslant K_{1} \| G \|^{1/ \nu'
     \gamma} ( 1+ \| G \|_{\gamma , \nu' , \psi'} )^{K_{2} \| G \|_{\gamma ,
     \nu' , \psi'}^{1/ \nu' \gamma}} ( 1+ | x | ) . \]
  Let $r>0$, $R=K_{1} \| G \|^{1/ \nu' \gamma} ( 1+ \| G \|_{\gamma , \nu' ,
  \psi'} )^{K_{2} \| G \|_{\gamma , \nu' , \psi'}^{1/ \nu' \gamma}} ( 1+ | r |
  )$ and $B=B ( 0,2R )$. Let us take $\tilde{G} \in \CC^{\gamma , \nu' ,
  \psi'} \cap \CC^{\gamma , \nu , \psi}$ such that $\| \tilde{G} \|_{\gamma ,
  \nu , \psi'} \leqslant \| G \|_{\gamma , \nu , \psi'}$, suppose furthermore
  that for $y \in B ( 0,r )$, $\tau_{\tilde{\theta} ( y )} G \in \CC^{\gamma
  ,1, \psi}$ where $\tilde{\theta}$ is the solution of the following equation
  \[ \tilde{\theta}_{t} ( y ) =y+ \int_{0}^{t} \tilde{G}_{\mathd u} (
     \tilde{\theta}_{u} ) . \]
  Then
  \[ \| \theta ( x ) - \tilde{\theta} ( y ) \|_{\infty} \lesssim \varphi ( R )
     e^{c_{2} \varphi ( R ) \| \tau_{\tilde{\theta} ( y )} G \|_{\gamma ,1,
     \psi}^{1/ \gamma}} ( | x-y | + \| G- \tilde{G} \|_{\gamma , \nu , \psi} ),
  \]
  where $\varphi ( R ) = ( R+1 ) \psi ( R )$.
\end{theorem}

\section{Averaging of paths}\label{sec:averaging_deterministic}

We turn now to the study of the averaging operator $T^{w}$ proper. One of our
main results is a proof that fBm paths are $\rho$-irregular for any $\rho
<1/2H$ and as a consequence that the averaging operator $T^{w}$ is bounded
from the Fourier--Lebesgue space $\CF L^{\alpha}$ to $\CC^{\gamma} \CF L^{\alpha
+ \rho}$ for any $\alpha \in \mathbbm{R}$ and for almost every fBm path $w$.
This result was one of our main reasons to look at the scale of
Fourier--Lebesgue spaces.

For the scale of Besov spaces $\left( \CC^{\alpha} \right)_{\alpha}$ we were
unable to prove similar results and we limited ourselves to study the averaged
vector-fields $T^{w} f$ for fixed $f \in \CC^{\alpha}$.

In this section we will first study the almost-sure irregularity of fBm paths.
This study proceeds in two steps: first we use well known chaining arguments
(essentially going back to Kolmogorov lemma in the form given to it by
Garsia, Rodemich and Rumsey) to go from supremum norm to ``integral'' norms
more suitable to probabilistic estimates and then use Hoeffding inequality to
prove these estimates.

The use of Hoeffding inequality replaces what in Davie's paper~{\cite{MR2377011}} are
explicit and painful computations on Brownian motions (relying on the Markov
property) and what in other works (e.g. in~{\cite{FGP}}) is achieved via
stochastic calculus (and thus martingale properties). In the fBm context
neither technique is applicable and explicit computations using Gaussian
tools, while possible are quite cumbersome and moreover we were unable to use
them to obtain the exponential square integrability we show here to be valid.
So we think that our observation that discrete martingale techniques like
Hoeffding inequality are useful in the fBm context is one of the interesting
points of our research.

\subsection{Chaining lemmas}

To see the average properties of the fractional Brownian path, we will need
some chaining lemmas, to infer global estimates from pointwise ones.

\begin{lemma}
  \label{lemma_chaining_time}Let $X$ from $I^{2}$ to $\mathbbm{R}^{d}$ such
  that for all $s \leqslant u \leqslant t$
  \[ | X_{s,t} | \leqslant | X_{s,u} | + | X_{u,t} |   \tmop{and}  X_{s,s} =0.
  \]
  And let us define for $\mu>0$,
  \[ R_{\mu} ( X ) = \sum_{n \in \mathbbm{N}} \sum_{k=0}^{2^{n} -1} 2^{-2n}
     \exp ( \mu 2^{n} | X_{k2^{-n} , ( k+1 ) 2^{-n}} |^{2} ). \]
  Then there exists a constant $K>0$ such that for all $s \leqslant t$,
  \[ \exp ( \mu | X_{s,t} |^{2} / | t-s | ) \lesssim | t-s |^{-K} R_{\mu K} (
     X ). \]
\end{lemma}

\begin{proof}
  Let $0 \leqslant s<t \leqslant 1$, $n \in \mathbbm{N}$ be the largest $n' \in
  \mathbbm{N}$ such that $2^{- (n' +1)} \leqslant t-s \leqslant 2^{-n'}$. By
  definition of $n$ there exists $l$ such that $l/2^{n} \leqslant s<t
  \leqslant (l+1) /2^{n}$. We can find some sequences $(s_{k} )_{k \ge 1}$ and
  $(t_{k} )_{k \ge 1}$ such that $(s_{k} )$ decreases, $(t_{k} )$ increases,
  $s_{1} =t_{1} = (2l+1) /2^{n+1}$, $\lim_{k \to \infty}  s_{k} =s$, $\lim_{k
  \to \infty}  t_{k} =t$, $s_{k+1} -s_{k} \leqslant 2^{n+k+1}$, $t_{k+1}
  -t_{k} \leqslant 2^{n+k+1}$ and $2^{n+k} t_{k} \in \mathbbm{Z}$ and $2^{n+k}
  s_{k} \in \mathbbm{Z}$. Hence $[s,t) = \cup_{k \ge 1} [s_{k+1} ,s_{k} ) \cup
  \cup_{k \ge 1} [t_{k} ,t_{k+1} )$ and thanks to the definition of the
  sequences, the following inequalities hold for $s_{k}$, but also for
  $t_{k}$.
  
  First, if $s_{k+1} =s_{k}$, $\sqrt{\mu} | X_{s_{k+1} ,s_{k}} | =0$. Now, if
  $s_{k+1} <s_{k}$ then there exists $l_{k} \in \{ 0, \ldots ,n+k \}$ such
  that $s_{k+1} = ( 2l_{k} -1 ) /2^{n+k+1}$ and $s_{k} =l_{k} /2^{n+k}$. Hence
  \begin{eqnarray*}
    \sqrt{\mu} | X_{s_{k+1} ,s_{k}} | & = & 2^{- ( n+k+1 ) /2} \log ( 2^{(
    n+k+1 )} 2^{- ( n+k+1 )} \exp ( \mu 2^{k+n+1} | X_{s_{k+1} ,s_{k}} |^{2} )
    )^{1/2}\\
    & \lesssim & 2^{- ( n+k ) /2} \{ ( n+k ) + \log ( R_{\mu} ( X ) )
    \}^{1/2} .
  \end{eqnarray*}
  But $2^{- ( n+1 )} \leqslant | t-s | \leqslant 2^{-n}$, hence
  \begin{eqnarray*}
    \sqrt{\mu} | X_{s_{k+1} ,s_{k}} | & \lesssim   & | t-s |^{1/2} 2^{-k/2} \{
    k+ \log (1/| t-s |) + \log ( R_{\mu} ( X ) ) \}^{1/2} .
  \end{eqnarray*}
  Thanks to the definition of $( s_{k} )_{k}$ and $( t_{k} )_{k}$, we have
  \begin{eqnarray*}
    \sqrt{\mu} | X_{s,t} | & \leqslant & \sum_{k \geqslant 1} \sqrt{\mu} |
    X_{s_{k+1} ,s_{k}} | + \sqrt{\mu} | X_{t_{k} ,t_{k+1}} |\\
    & \lesssim & | t-s |^{1/2} \{ 1+ \log \big(1/| t-s |\big) + \log ( R_{\mu} ( X )
    ) \}^{1/2} .
  \end{eqnarray*}
  Hence,
  \begin{eqnarray*}
    \exp ( \mu | X_{s,t} |^{2} / | t-s | ) & \lesssim & \exp ( K  \log ( 1/ |
    t-s | ) +K  \log  R_{\mu} ( X ) )\\
    & \lesssim & | t-s |^{-K} R_{\mu K} ( X )
  \end{eqnarray*}
  and by Jensen inequality $R_{\mu} ( X )^{K} \leqslant R_{\mu K} ( X )^{}$.
\end{proof}

In the following, to approach a point of $\mathbbm{R}^{d}$ we will use a
similar argument. Namely we will use the graph $( 2^{-m} \mathbbm{Z} )^{d}$ as
a good approximation of $\mathbbm{R}^{d}$. Hence we need to have an
approximation of the biggest error we can make using such an approximation. It
is well known that for all $d$ and all $m \in \mathbbm{N}$, $\sup_{x \in
\mathbbm{R}^{d}} \inf_{y \in ( 2^{-m} \mathbbm{Z} )^{d}} | x-y | = \sqrt{d}
/2^{m+1}$.

\begin{lemma}
  \label{lemma_chaining_space}Let $X$ be a function from $\mathbbm{R}^{d}$ to
  $\mathbbm{R}^{d}$ and $g$ such that $g \geqslant 1$, $\sup_{| \zeta - \zeta'
  | \leqslant \sqrt{d} /2} g ( \zeta' ) /g ( \zeta ) < \infty$ and with $\|
  g^{-1} \|_{L^{1} ( \mathbbm{R}^{d} )} <+ \infty$. Suppose furthermore that
  the following quantity is finite
  \[ C_{X} \assign \sup_{\tmscript{\begin{array}{c}
       m \in \mathbbm{N}\\
       \zeta : 2 \geqslant g ( \zeta ) /2^{m} \geqslant 1/2\\
       | \zeta - \zeta' | \leqslant \sqrt{d} /2
     \end{array}}} | X ( \zeta ) -X ( \zeta' ) | / ( 2^{m} | \zeta - \zeta' |
     ) <+ \infty . \]
  Let
  \[ S_{\mu} ( X ) = \sum_{n \in \mathbbm{N}} \sum_{\zeta \in ( 2^{-n}
     \mathbbm{Z} )^{d}} 2^{- ( d+1 ) n} g ( \zeta' )^{-1} \exp ( \mu | X (
     \zeta' ) |^{2} ) . \]
  Then, there exists a constant $C \geqslant 1$ be such that
  \[ \exp ( \mu | X ( \zeta ) |^{2} ) \lesssim g ( \zeta )^{-C} \exp ( \mu K
     C_{X}^{2} ) S_{\mu C} ( X ) . \]
\end{lemma}

\begin{proof}
  Let $\zeta \in \mathbbm{R}^{d}$ and $m$ such that $g ( \zeta ) \sim 2^{m}$.
  Let $\zeta' \in ( 2^{-m} \mathbbm{Z} )^{d}$ such that $| \zeta - \zeta' |
  \leqslant 2^{-m} \sqrt{d} /2$ then
  $ | X ( \zeta ) -X ( \zeta' ) | \leqslant C_{X} 2^{m} | \zeta - \zeta' |
     \lesssim C_{X} $. 
  Furthermore, the hypothesis on $g$ gives us that $\log ( g ( \zeta' ) )
  \lesssim 1+ \log ( g ( \zeta ) )$. Hence
  \begin{eqnarray*}
    \sqrt{\mu} | X ( \zeta ) | & \leqslant & \sqrt{\mu} | X ( \zeta ) -X (
    \zeta' ) | + \sqrt{\mu} | X ( \zeta' ) |\\
    & \leqslant & \sqrt{\mu} C_{X} + \{ \log (^{} 2^{m} 2^{-m} g ( \zeta' ) g
    ( \zeta' )^{-1} \exp ( \mu | X ( \zeta' ) |^{2} ) ) \}^{1/2}\\
    & \lesssim & \sqrt{\mu} C_{X} + \{ m+ \log ( g ( \zeta' ) ) + \log (
    S_{\mu} ( X ) ) \}^{1/2}\\
    & \lesssim & \sqrt{\mu} C_{X} + \{ 1+ \log ( g ( \zeta ) ) + \log (
    S_{\mu} ( X ) ) \}^{1/2} .
  \end{eqnarray*}
  Finally we have
  \begin{eqnarray*}
    \exp ( \mu | X ( \zeta ) |^{2} ) & \leqslant & \exp ( K \{ \mu C_{X}^{2}
    +1+ \log ( g ( \zeta ) ) + \log  S_{\mu} ( X ) \} )\\
    & \lesssim & g ( \zeta )^{K} \exp ( \mu C C_{X}^{2} ) S_{\mu C} ( X ) .
  \end{eqnarray*}
\end{proof}

We can think of $g$ as $g ( \zeta ) = ( 1+ | \zeta | )^{d+1}$.

\begin{lemma}
  \label{lemma_holder_Y}For all $\beta \in \mathbbm{R}$ and all $R>0$ there
  exists a constant $C ( \beta ,R )$ such that for all $\zeta' \in B ( \zeta
  ,R )$
  \[ | ( 1+ | \zeta | )^{\beta} - ( 1+ | \zeta' | )^{\beta} | \leqslant C (
     \beta ,R ) ( 1+ | \zeta | )^{\beta -1} | \zeta - \zeta' |. \]
\end{lemma}

\begin{proof}
  Let us suppose first that $| \zeta' | \geqslant | \zeta |$ by the choice of
  $\zeta'$ we have $0 \leqslant \frac{| \zeta' | - | \zeta |}{1+ | \zeta |}
  \leqslant R$. Then
  \begin{eqnarray*}
    | ( 1+ | \zeta | )^{\beta} - ( 1+ | \zeta' | )^{\beta} | & = & ( 1+ |
    \zeta | )^{\beta} \left\lvert \left( 1+ \frac{| \zeta' | - | \zeta |}{1+ |
    \zeta |} \right)^{\beta} -1 \right\lvert\\
    & \leqslant & ( 1+ | \zeta | )^{\beta} \sup_{x \in [ 1,R ]} | f'_{\beta}
    ( x ) | \frac{| | \zeta' | - | \zeta | |}{1+ | \zeta |}\\
    & \leqslant & \sup_{x \in [ 0,R ]} | f'_{\beta} ( x ) | ( 1+ | \zeta |
    )^{\beta -1} | \zeta - \zeta' |,
  \end{eqnarray*}
  where the function $f_{\beta}$ is define from $[ 0,R ]$ to $\mathbbm{R}$ by
  \ $f_{\beta} ( x ) = ( 1+x )^{\beta}$. We have
  \begin{eqnarray*}
    | f_{\beta}' ( x ) | & = & | \beta ( 1+x )^{\beta -1} | \leqslant | \beta
    | ( ( 1+R )^{\beta -1} \vee 1 ) .
  \end{eqnarray*}
  If $| \zeta | > | \zeta' |$, the same computation gives
  \[ | ( 1+ | \zeta | )^{\beta} - ( 1+ | \zeta' | )^{\beta} | \lesssim_{R,
     \beta} ( 1+ | \zeta' | )^{\beta -1} | \zeta - \zeta' | . \]
  When $\beta -1 \geqslant 0$, the result follows. Suppose now that 
  $\beta-1<0$, we have to prove that $( 1+ | \zeta | ) \lesssim ( 1+ | \zeta' | )$.
  When $| \zeta | \leqslant 2R$, then we have
  \begin{eqnarray*}
    ( 1+ | \zeta' | ) / ( 1+ | \zeta | ) & \geqslant & 1/ ( 1+2R ) .
  \end{eqnarray*}
  When $| \zeta | >2R$,
  \[ ( 1+ | \zeta' | ) / ( 1+ | \zeta | ) \geqslant ( 1+ | \zeta | - | \zeta'
     - \zeta |   ) / ( 1+ | \zeta | ) \geqslant 1- | \zeta' - \zeta | / ( 1+ |
     \zeta | ) \geqslant 1/2 \]
  and the result follows.
\end{proof}

\subsection{Application of the chaining lemmas, control of the averaging along
curves.}

The last lemmas allows us to control the average of a function (or a
distribution) along the curve $w$. Indeed, to estimate on the quantity
$\int_{s}^{t} f_{u} ( x+w_{u} ) \mathd u$ it will be enough to have a control
on simpler quantities. We will apply those lemmas in two similar situations,
namely when $f \in \CC^{\alpha}$ and when $f \in \CF L^{\alpha}$. In this latter
case, we will see that it is enough to control $\Phi^{w}$.

\subsubsection{Averaging property of the occupation measure}

Recall that we have already defined $\Phi^{w}_{t} ( \xi ) = \int_{0}^{t} e^{i
\langle \xi ,w_{r} \rangle} \mathd r$ and
\[ \| \Phi^{w} \|_{\mathcal{W}^{\rho , \gamma}_{T}} = \sup_{\xi \in
   \mathbbm{R}^{d}}   \sup_{0 \leqslant s<t \leqslant T} ( 1+ | \xi | )^{\rho}
   \frac{| \Phi^{w}_{t} ( \xi )- \Phi^{w}_{s} ( \xi ) |}{|s-t|^{\gamma}} . \]

\begin{lemma}
  \label{lemma:chaining-Y}For all $- \beta < \alpha$ there exist a constant
  $a>0$ and $\gamma >0$ such that for all $\lambda >0$,
  \[ | \Phi^{w}_{t} ( \xi )- \Phi^{w}_{s} ( \xi ) | \lesssim | t-s |^{\gamma}
     ( 1+ | \xi | )^{- \alpha'} ( 1+ \log^{1/2} ( e^{a \mu \| w \|_{\infty}}
     K^{w}_{\alpha} ( \lambda ) ) ), \]
  where
  \[ K^{w}_{\alpha} ( \lambda ) = \sum_{\tmscript{\begin{array}{c}
       n,m \in \mathbbm{N}\\
       0 \leqslant k \leqslant 2^{n} -1\\
       \xi' \in ( 2^{-m} \mathbbm{Z} )^{d}
     \end{array}}} 2^{-2n+ ( d+1 ) m} ( 1+ | \xi' | )^{- ( d+1 )} \exp (
     \lambda 2^{n} ( 1+ | \xi' | )^{2 \beta} | \Phi^{w}_{k2^{-n} , ( k+1 )
     2^{-n}} ( \xi' ) |^{2} ) . \]
\end{lemma}

\begin{proof}
  We apply the Lemmas~\ref{lemma_chaining_time} and~\ref{lemma_chaining_space}
  to
  \[ X_{s,t} ( \xi ) = ( 1+ | \xi | )^{\beta} | \Phi_{s,t}^{w} ( \xi ) | / |
     t-s |^{1/2}. \]
  Thanks to Lemma~\ref{lemma_holder_Y} and the definition of $\Phi^{w}_{s,t}$,
  for all $\xi \in \mathbbm{R}^{d}$, and all $\xi' \in B \left( \xi , \sqrt{d}
  /2 \right)$, we have
  \begin{eqnarray*}
    | X_{s,t} ( \xi ) -X_{s,t} ( \xi' ) | & \leqslant & | ( 1+ | \xi |
    )^{\beta} - ( 1+ | \xi' | )^{\beta} |   | \Phi^{w}_{s,t} ( \xi' ) | / |
    t-s |^{1/2}\\
    &  & + ( 1+ | \xi | )^{\beta} | \Phi^{w}_{s,t} ( \xi ) - \Phi^{w}_{s,t} (
    \xi' ) | / | t-s |^{1/2}\\
    & \lesssim & ( 1+ | \xi | )^{\beta} | \xi - \xi' | ( 1+ \| w \|_{\infty}
    ) | t-s |^{1/2}
  \end{eqnarray*}
  Here we take $\zeta = \xi$, $g ( \zeta ) = ( 1+ | \zeta | )^{\beta +C d}$
  such that $\beta +C d \geqslant d+1$. With those choices, $X_{s,t}$ and $g$
  verify the hypothesis of lemma \ref{lemma_chaining_space}, furthermore
  $C_{X} \lesssim ( 1+ \| w \|_{\infty} )$,hence
  \begin{eqnarray*}
    \exp ( \mu | X_{s,t} ( \xi ) |^{2} ) & = & \exp ( \mu ( 1+ | \xi | )^{2
    \beta} | \Phi^{w}_{s,t} ( \xi ) |^{2} / | t-s | )\\
    & \lesssim & ( 1+ | \xi | )^{C ( d+1 )} S_{C \mu} ( X_{s,t} ) \exp ( \mu
    C \| w \|_{\infty}^{2} ).
  \end{eqnarray*}

  Now, let us apply Lemma \ref{lemma_chaining_time} to $( 1+ | \xi | )^{\beta}
  \Phi^{w}_{s,t} ( \xi )$, then
  \[ \exp ( \mu | X_{s,t} ( \xi ) |^{2} ) \lesssim | t-s |^{-K} R_{K \mu} ( (
     1+ | \xi | )^{\beta} \Phi^{w}_{.} ( \xi ) ) .\]
  But
  \begin{eqnarray*}
    R_{K \mu} ( ( 1+ | \xi | )^{\beta} \Phi_{.} ( \xi ) ) & = & \sum_{n \in
    \mathbbm{N}} \sum_{k=0}^{2^{n} -1} 2^{-2n} \exp ( \mu K2^{n} |
    \Phi^{w}_{k2^{-n} , ( k+1 ) 2^{-n}} ( \xi ) |^{2} ( 1+ | \xi | )^{2 \beta}
    )\\
    & \lesssim & ( 1+ | \xi | )^{C ( d+1 )} \sum_{n \in \mathbbm{N}}
    \sum_{k=0}^{2^{n} -1} 2^{-2n} S_{\mu C K} ( X_{k2^{-n} , ( k+1 ) 2^{-n}} )
    \exp ( \mu C K \| w \|_{\infty}^{2} )\\
    & \lesssim & ( 1+ | \xi | )^{C ( d+1 )} \exp ( \lambda \| w
    \|_{\infty}^{2} ) K^{w}_{\beta} ( \lambda ) .
  \end{eqnarray*}
  When we take the logarithm, we have
  \[ | \Phi_{s,t}^{w} ( \xi ) | \lesssim \mu^{-1/2} | t-s |^{1/2} ( 1+ | \xi |
     )^{- \beta} ( 1+ \log(1/| t-s |) + \log ( 1+ | \xi | ) + \log ( \exp (
     \lambda \| w \|_{\infty}^{2} ) K^{w}_{\beta} ( a \mu ) ) ) . \]
  Hence, for all $\varepsilon_{1} , \varepsilon_{2} >0$, we have
  \[ | \Phi_{s,t}^{w} ( \xi ) | \lesssim_{\varepsilon_{1} , \varepsilon_{2}}
     \mu^{-1/2} | t-s |^{1/2- \varepsilon_{1}} ( 1+ | \xi | )^{- \beta +
     \varepsilon_{2}} ( 1+ \log ( \exp ( \lambda \| w \|_{\infty}^{2} )
     K^{w}_{\beta} ( a \mu ) ) ) . \]
  Furthermore, by interpolating with the trivial estimate $| \Phi_{s,t} ( \xi
  ) | \leqslant | t-s |$, for all $- \beta < \alpha$, there exists $\gamma
  >1/2$, and a constant $a>0$ such that
  \[ | \Phi_{s,t}^{w} ( \xi ) | \lesssim | t-s |^{\gamma} ( 1+ | \xi |
     )^{\alpha} ( 1+ \log^{1/2} ( \exp ( a \mu \| w \|_{\infty}^{2} )
     K^{w}_{\beta} ( a \mu ) ) ) . \]
  
\end{proof}

\subsubsection{Averaging of Besov functions along paths}

In this section we analyse the averaging effect of paths on functions
belonging to the scale of Besov spaces $( \CC^{\alpha} )_{\alpha}$. Note the following. If we write $\tilde \Delta_i = \sum_{j:|i-j|\le 1} \Delta_j$, we have $\tilde \Delta_i \Delta _i = \Delta _i$ for all  $i\ge -1$  and then
$
T^{w}_{s,t} ( \Delta_{i} f ) ( x ) =T^{w}_{s,t} ( \tilde \Delta_{i} \Delta_{i} f ) ( x ) =(T^{w}_{s,t} ( \tilde K_{i}  ) * \Delta_{i} f) ( x )
$
where $\tilde K_i$ is the integral kernel corresponding to the operator $\tilde \Delta_{i}$. In this case
\[
\|T^{w}_{s,t} ( \Delta_{i} f )\|_{L^\infty} \lesssim \| \Delta_{i} f \|_{L^\infty}  \|T^{w}_{s,t} ( \tilde K_{i} )\|_{L^1} .
\]
So any control of quantities like
$
\sum_{i \ge -1} \Psi(2^{\alpha i} \|T^{w}_{s,t} ( \tilde K_{i} )\|_{L^1})
$
for increasing functions $\Psi$ will imply bounedness properties of $T^w$ in H\"older--Besov spaces. However in the case of fractional Browian sample path (or even just in the case of Brownian motion) we were unable to devise useful estimates for this kind of quantities.
Due to this difficulty which prevents us from having (useful) estimates which are uniform in $\CC^{\alpha}$,  the chaining argument now depends on the chosen function $f$ and the computations follows closely those in the previous section.

\begin{lemma}
  \label{lemma:chaining_besov}For all $- \beta < \alpha$ there exists $\gamma
  >1/2$ such that for all $f \in \mathcal{S}' ( \mathbbm{R}^{d} )$, all
  $\lambda >0$ and all $i$,
  \[ | T^{w}_{s,t} ( \Delta_{i} f ) ( x ) | \lesssim_{\lambda} 2^{\alpha i}
     \| \Delta_{i} f \|_{\infty} | t-s |^{\gamma} ( 1+ \log^{1/2} ( 1+ | x | )
     + \log^{1/2} ( K^{w}_{f, \beta} ( \lambda ) ) ), \]
  where
  \[ K_{f, \beta}^{w} ( \lambda ) = \sum_{\tmscript{\begin{array}{c}
       n,m \in \mathbbm{N}\\
       i \geqslant -1\\
       0 \leqslant k \leqslant 2^{n} -1\\
       x' \in ( 2^{-m} \mathbbm{Z} )^{d}
     \end{array}}} \frac{2^{-c_{\beta} ( m+n+i )}} {( 1+ | x' | )^{( d+1 )}} \exp (
     \lambda 2^{n+2i \beta} | T^{w}_{k/2^{n} , ( k+1 ) /2^{n}} ( \Delta_{i} f
     ) ( x' ) |^{2} / \| \Delta_{i} f \|_{\infty}^{2} ) \]
  and $c_{\beta}$ is a constant depending only of $\beta$ and $d$ such that
  the sum without the exponential is finite.
\end{lemma}

\begin{proof}
  The proof is very similar to the proof of Lemma \ref{lemma:chaining-Y}. We
  will apply the Lemmas~\ref{lemma_chaining_time}
  and~\ref{lemma_chaining_space} to
  \[ X^{i}_{s,t} ( x ) =2^{i \alpha} | T^{w}_{s,t} ( \Delta_{i} f ) ( x ) | /
     ( \| \Delta_{i} f \|_{\infty} | t-s |^{1/2} ), \]
  with the convention that $X^{i} =0$ when $\Delta_{i} f=0$. We have, thank to
  the definition of $T^{w}$,
  \begin{equation}
    | T^{w}_{s,t} ( \Delta_{i} f ) ( x ) | \leqslant \| \Delta_{i} f
    \|_{\infty} | t-s |^{} \label{estimation-interpolation-1}.
  \end{equation}
  Furthermore, as the Fourier transform of $\Delta_{i} f$ is compactly
  supported in an annulus, we have the obvious estimate
  \begin{eqnarray*}
    | X^{i}_{s,t} ( x ) -X^{i}_{s,t} ( x' ) | & \lesssim & 2^{i \beta} \|
    \nabla \Delta_{i} f \|_{\infty} | t-s |^{1/2} | x-x' | / \| \Delta_{i} f
    \|_{\infty}\\
    & \lesssim & 2^{i ( \beta +1 )} | x-x' | .
  \end{eqnarray*}
  Let us take $g_{i} ( x ) =2^{( a+ \beta ) i} ( 1+ | x | )^{d+1}$ with $a
  \geqslant 1$ and $a+ \beta \geqslant c_{\beta}$. Hence, $C_{X^{i}_{s,t}}
  \lesssim 1$. By Lemma \ref{lemma_chaining_space}, there exists $b,c>1$ such
  that
  \begin{eqnarray*}
    \exp ( \mu | X^{i}_{s,t} ( x ) |^{2} ) & \lesssim & 2^{b  ( a+ \beta ) i}
    ( 1+ | x | )^{b ( d+1 )} S_{\mu b} ( X^{i}_{s,t} ) .
  \end{eqnarray*}
  Now, thanks to Lemma~\ref{lemma_chaining_time}, there exists $a' ,b' ,c'
  ,d'$ such that
  \begin{eqnarray*}
    \exp ( \mu | X^{i}_{s,t} ( x ) |^{2} ) & \leqslant & 2^{a' i} ( 1+ | x |
    )^{b'} | t-s |^{-c'} K^{w}_{f, \beta} ( d' \mu ) .
  \end{eqnarray*}
  Hence by taking the logarithm, and by losing a small power of time and on
  $i$, we have
  \begin{equation}
    | T^{w}_{s,t} \Delta_{i} f ( x ) |^{} \lesssim 2^{( - \beta +
    \varepsilon_{1} ) i} \| \Delta_{i} f \|_{\infty} | t-s |^{1/2-
    \varepsilon_{2}} ( 1+ \log^{1/2} ( 1+ | x | ) + \log^{1/2} ( K_{f,
    \beta}^{w} ( d' \mu ) ) ) \label{estimation-interpolation-2}.
  \end{equation}

  Now we interpolate~{\eqref{estimation-interpolation-2}}
  and~{\eqref{estimation-interpolation-1}} and for all $\alpha
  >- \beta$, there exists $\rho >0$ and $\gamma >1/2$ such that
  \[ | T^{w}_{s,t} ( \Delta_{i} f ) ( x ) | \lesssim 2 ^{\alpha i} \|
     \Delta_{i} f \|_{\infty} | t-s |^{\gamma} ( 1+ \log^{1/2} ( 1+ | x | ) +
     \log ( K_{f, \beta}^{w} ( d' \mu ) ) )^{1/2}.  \]
  
\end{proof}

\subsubsection{The operator $T^{w}$}
\label{sec:bounds-on-Tw}

We are now able to define the function $T^{w} f$ for all $f \in \CC^{\alpha}$
(respectively $\CF L^{\alpha}$) for all $\alpha >- \beta$, as soon as there
exists $\lambda >0$ small enough such that $K^{w}_{f, \beta} ( \lambda )$
(respectively $K^{w}_{\beta} ( \lambda )$) is finite. As already mentioned,
it remains an open problem to study the boundedness of $T^{w}$ as an operator
with range in Besov spaces so we restrict ourselves to study the image of
$T^{w} f$ for fixed $f$ and with $w$ in the support of the fBm law without any
attempt to obtain estimates which are uniform in $f$. On the contrary, for
$\CF L^{\alpha}$, the estimate on $\Phi^{w}$ are good enough to define $T^{w}$
as an operator on the whole space.

\begin{definition}
  Let $\beta \in \mathbbm{R}$, $\alpha >- \beta$ and let $f \in \CC^{\alpha}$.
  We define
  \[ T^{w}_{s,t} f ( x ) = \sum_{i \geqslant -1} T^{w}_{s,t} ( \Delta_{i} f )
     ( x )  =  \lim_{N \rightarrow \infty} T^{w}_{s,t} ( \pi_{\leqslant N} f )
     ( x ) \]
  and
  \[ T^{w}_{s,t} g ( x ) = \lim_{\tmscript{\begin{array}{c}
       h \in \CF L^{0 \vee \alpha}\\
       h \xrightarrow{\CF L^{\alpha}} g
     \end{array}}} T^{w}_{s,t} h ( x ). \]
\end{definition}

As these objects are defined by some limiting procedures, it is not
straightforward that they exist. Furthermore, for the consistency of the
definition, we must show that when $f \in \CF L^{\alpha}$ these two limiting
procedures give the same object, and that the limit does not depends of the
choice of sequence $( \Delta_{i} )$. This is the purpose of the following
theorem.

\begin{theorem}
  \label{theorem:definition_T}Let $\beta \in \mathbbm{R}$ and let $\alpha >-
  \beta$.
  \begin{enumerateroman}
    \item Suppose that there exists $\lambda_{0}$ such that $K^{w}_{\beta} (
    \lambda_{0} ) <+ \infty$. Then for all $g \in \CF L^{\alpha}$, $T^{w} g$
    exists and does not depends of the choice of the sequence. Furthermore,
    for all $\lambda \leqslant \lambda_{0}$ we have
    \[ | T^{w}_{s,t} g ( x ) | \lesssim | t-s |^{\gamma} N_{\alpha} ( g ) ( 1+
       \log^{1/2} K^{w}_{\beta} ( \lambda ) ) . \]
    Hence, $( T^{w}_{s,t} )_{0 \leqslant s \leqslant t \leqslant 1}$ is well
    defined as a family of operators on $\CF L^{\alpha}$.
    
    \item For $f \in \CC^{\alpha}$ suppose that there exists $\lambda_{0}$
    such that $K^{w}_{f, \beta} \nocomma ( \lambda_{0} ) <+ \infty$. Then
    $T^{w} f$ exists and the following bound holds
    \[ | T^{w} f_{s,t} ( x ) | \lesssim_{\lambda} | t-s |^{\gamma} \| f
       \|_{\alpha} ( 1+ \log^{1/2} ( 1+ | x | ) + \log^{1/2} K_{f, \beta}^{w}
       ( \lambda ) ) . \]
    Furthermore let us suppose that for $g \in \CC^{\alpha}$, $K^{w}_{g,
    \beta} ( \lambda_{0} )$ is also finite, then for all $\lambda \leqslant
    \lambda_{0} /2$
    \[ | T^{w}_{s,t} ( f-g ) ( x ) | \lesssim | t-s |^{\gamma} \| f-g
       \|_{\alpha} ( 1+ \log^{1/2} ( 1+ | x | ) + ( K^{w}_{f, \beta} ( \lambda
       ) +K^{w}_{g, \beta} ( \lambda ) ) ). \]
    \item These two limiting procedures are compatible when $f \in \CF
    L^{\alpha}$.
  \end{enumerateroman}
\end{theorem}

\begin{proof}
  The proof is quite straightforward when $g \in \CF L^{\alpha}$. Indeed, for
  $h^{1}$ and $h^{2}$ in $\CF L^{0} \cap \CF L^{\alpha}$, we have
  \[ T_{s,t}^{w} ( h_{1} -h_{2} ) ( x ) = \int_{\mathbbm{R}^{d}} \mathd \xi (
     \hat{h}_{1} - \hat{h}_{2} ) ( \xi ) \exp ( i \xi \cdot x ) \Phi_{s,t}^{w} (
     \xi ), \]
  hence
  \[ | T^{w}_{s,t} ( h_{1} -h_{2} ) ( x ) | \lesssim N_{\alpha} ( h_{1}
     -h_{2} ) | t-s |^{\gamma} ( 1+ \log^{1/2} ( K^{w}_{\beta} ( \lambda_{0} )
     ) ) \]
  and the result of $( i )$ follows.

  Let us prove $( i i )$. For $f \in \CC^{\alpha}$ let us show as the quantity
  $T^{w}_{s,t} ( \pi_{\leqslant N} f ) ( x )$ converges when $N \rightarrow +
  \infty$. Indeed, thanks to Lemma \ref{lemma:chaining_besov}, for
  $\varepsilon >0$ such that $- \beta < \alpha - \varepsilon < \alpha$, we
  have
  \begin{eqnarray*}
    | T^{w}_{s,t} ( \pi_{\leqslant N} f ) ( x ) -T^{w}_{s,t} ( \pi_{\leqslant
    N+M} f ) ( x ) | & \lesssim & \sum_{i=N+1}^{N+M} | T^{w}_{s,t} (
    \Delta_{i} f ) ( x ) |\\
    & \lesssim & C_{s,t,x} ( \log^{1/2} K_{f, \beta}^{w} ( \lambda_{0} ) ) \|
    f \|_{\alpha} 2^{- \varepsilon N} .
  \end{eqnarray*}
  Hence, $\big( T^{w}_{s,t} ( \pi_{\leqslant N} f ) ( x ) \big)_{N\geqslant 0}$ is Cauchy,
  and then the limit $T^{w} f$ exists. Furthermore we have the straightforward
  bound for all $\lambda < \lambda_{0}$
  \[ | T^{w}_{s,t} ( \pi_{\leqslant N} f ) ( x ) | \lesssim | t-s |^{\gamma}
     \| f \|_{\alpha} ( 1+ \log^{1/2} ( 1+ | x | )   \noplus \noplus +
     \log^{1/2} ( K^{w}_{f, \beta} ( \lambda ) ) ) \]
  and the same bound holds as $N \rightarrow + \infty$.
  For $f$ and $g$ we have
  \[ | T^{w}_{s,t} ( f-g ) ( x ) | \lesssim | t-s |^{\gamma} \| f-g
     \|_{\alpha} ( 1+ \log^{1/2} ( 1+ | x | )    + \log^{1/2} (
     K^{w}_{f-g, \beta} ( \lambda ) ) ) , \]
  but thanks to the definition of the constants, $K^{w}_{f-g, \beta} ( \lambda
  ) \lesssim K^{w}_{f, \beta} ( 2 \lambda ) +K^{w}_{g, \beta} ( 2 \lambda )$,
  and the result follows.
  For $( i i i )$ let us consider $f \in \CF L^{\alpha}$, we have
  \begin{eqnarray*}
    | T^{w}_{s,t} ( \pi_{\leqslant N} f ) ( x ) -T^{w}_{s,t} ( \pi_{\leqslant
    N+M} f ) ( x ) | & = & | T^{w}_{s,t} ( \pi_{\leqslant N} f- \pi_{\leqslant
    N+M} f ) ( x ) |\\
    & \lesssim & C_{s,t,x} N_{\alpha} ( \pi_{\leqslant N} f- \pi_{\leqslant
    N+M} f )   \log^{1/2} K_{\beta}^{w} ( \lambda_{0} )
  \end{eqnarray*}
  and as the sequence $( \pi_{\leqslant N} f )_{N}$ converges in $\CF
  L^{\alpha}$ the limit also exists, and of course it is the same.
  Furthermore, for two functions in $\CF L^{\alpha \vee 0}$,
  \[ | T^{w}_{s,t} ( \pi_{\leqslant N} f ) ( x ) -T^{w}_{s,t} (
     \pi_{\leqslant N+M} f ) ( x ) | \lesssim_{s,t,x} N_{\alpha} ( f-g ) \]
  and the limiting procedure in the $\CF L^{\alpha}$ case is correct.
\end{proof}

\begin{remark}
  The definition of $T^{w} f$ given above seems to depend on the choice of the
  Littlewood-Paley decomposition $( \Delta_{i} )_{i}$. It is indeed the fact.
  When we will consider $w$ being a stochastic process, this will lead us to a
  choice of a version of this averaging process defined almost surely. In
  fact, if \ $( \widetilde{\Delta_{}}_{i} )_{i}$ is another sequence of
  Littlewood--Paley operators, and $\tilde{K}_{i}$ is the associated integral
  kernels, we have
  \begin{eqnarray*}
    | K_{j} | \ast \log ( 1+ | \cdot | ) ( x ) & \leqslant & \log ( 1+ | x | ) +
    \int_{\mathbbm{R}^{d}} \mathd y  | \tilde{K}_{j} ( y ) | | \log ( 1+ | x-y
    | ) - \log ( 1+ | x | ) |\\
    & \lesssim & \log ( 1+ | x | ) +2^{-j} \int_{\mathbbm{R}^{d}} | 2^{j} y |
    2^{j d} \tilde{K} ( 2^{j} y ) \mathd y\\
    & \lesssim & 1+ \log ( 1+ | x | ).
  \end{eqnarray*}
  Hence, there exist two constants $c<C$ such that for all $\alpha >- \beta$ and for all $\varepsilon >0$ small enough,
  we have
  \begin{eqnarray*}
    | ( T^{w}_{s,t} \tilde{\Delta}_{j} f ) ( x ) | & \leqslant & \sum_{c + j \leqslant i
    \leqslant C + j} | ( \tilde{\Delta}_{j} T^{w} \Delta_{i} f ) ( x ) |\\
    & \lesssim & 2^{-j \varepsilon} \| f \|_{\alpha} | t-s |^{\gamma} \sum_{c + j \leqslant i
    \leqslant C + j} \int_{\mathbbm{R}^{d}} \tilde{K}_{j} \ast ( 1+ \log ( 1+ | . | )
    +K_{f, \beta}^{w} ( \lambda_{0} ) ) ( x )\\
    & \lesssim & 2^{-j \varepsilon} \| f \|_{\alpha} ( 1+ \log ( 1+ | x | )
    +K_{f, \beta}^{w} ( \lambda_{0} ) ),
  \end{eqnarray*}
  which gives the convergence of $T^{w} \tilde{\pi}_{\leqslant N} f$ to a
  limit we called $T^{w \nocomma , \tilde{\Delta}} f$, and with the same
  stochastic constant $K_{f, \beta}^{w} ( \lambda_{0} )$.
\end{remark}

In order
 to apply
  the results
   of the
    section
     related
      to the Young
       integral, it is
necessary to have a better understanding of the space regularity of an average
function. Thanks to the property of the operator $T^{w}$, as soon as we ask
$f$ to be regular, this regularity will holds. Furthermore the definition of
$T^{w}$ allows us to differentiate it whenever $f$ is regular enough, and the
constant is finite. Namely we have the following propositions.

\begin{proposition}
  \label{prop:holder_continuity_of_averaging_path}Let $\nu \in [ 0,1 ]$,
  $\alpha >- \beta$, and $f \in \CF L^{\alpha + \nu}$ (respectively in
  $\CC^{\alpha + \nu}$). Furthermore we suppose that there exists $\lambda_{0}
  >0$ such that $K^{w}_{\beta} ( \lambda_{0} ) < \infty$ (respectively
  $K^{w}_{\nabla f, \beta} ( \lambda_{0} ) <+ \infty$). Then $T^{w} f \in
  \CC^{\gamma , \nu}_{b}$ (respectively $T^{w} f \in \CC^{\gamma , \nu ,
  \psi}$ where $\psi ( r ) =1+ \log^{1/2} ( 1+r )$) and the following bounds
  hold
  \[ | T^{w}_{s,t} f ( x ) -T^{w}_{s,t} f ( y ) | \lesssim N_{\alpha + \nu} (
     b ) | t-s |^{\gamma} | x-y |^{\nu} ( 1+ \log^{1/2} K^{w}_{\beta} (
     \lambda_{0} ) ) \]
  ( respectively  
  \[ | T^{w}_{s,t} f ( x ) -T^{w}_{s,t} f ( y ) | \lesssim | t-s |^{\gamma} |
     x-y |^{\nu} \| f \|_{\alpha + \nu} ( \psi ( | x | + | y | ) + \log^{1/2}
     K^{w}_{\nabla f, \beta} ( \lambda_{0} ) + \log^{1/2} K^{w}_{f, \beta} (
     \lambda_{0} ) \text{)}. \]
\end{proposition}

\begin{proof}
  For $f \in \CC^{\alpha + \nu}$, and for all $\beta < \alpha' < \alpha$
  \begin{equation*}
  \begin{split}
    | T^{w} & \Delta_{i} f ( x ) -T^{w} \Delta_{i} f ( y ) |  =  \left\lvert
    \int_{0}^{1} \nabla T^{w} \Delta_{i} f ( r  ( x-y ) +y ) \cdot ( x-y ) \mathd
    r \right\lvert\\
    & \leqslant  | x-y | \sup_{r \in [ 0,1 ]} | T^{w} \Delta_{i} \nabla f (
    r  ( x-y ) +y ) |\\
    & \lesssim  2^{i \alpha} \| \Delta_{i} \nabla f \|_{\infty} | x-y | |
    t-s |^{\gamma} ( 1+ \sup_{r \in [ 0,1 ]} \log^{1/2} ( 1+ | r  ( x-y ) +y |
    ) +K^{w}_{\nabla f, \beta} ( \lambda_{0} ) )\\
    & \lesssim   2^{i ( \alpha +1 )} \| \Delta_{i} f \|_{\infty} | x-y | |
    t-s |^{\gamma} ( 1+ \log^{1/2} ( 1+ | x | + | y | ) +K^{w}_{\nabla f,
    \beta} ( \lambda_{0} ) ).
    \end{split}
  \end{equation*}
  Furthermore, we also have
  \[ | T^{w} \Delta_{i} f ( x ) -T^{w} \Delta_{i} f ( y ) | \lesssim 2^{i
     \alpha} \| \Delta_{i} f \|_{\infty} | t-s |^{\gamma} ( 1+ \log^{1/2} ( 1+
     | x | ) + \log^{1/2} ( 1+ | y | ) +K^{w}_{f, \beta} ( \lambda_{0} ) ) ,
  \]
  and by interpolation we have the bound for $T^{w} \Delta_{i} f$. The
  argument of Theorem \ref{theorem:definition_T} gives us the result. A
  similar argument holds when $f \in \CF L^{\alpha}$.
\end{proof}

The next proposition shows that the definition of the averaging operator $T$
is compatible with the space differential in the Hölder spaces.

\begin{proposition}
  \label{prop:diffenrentiability_averaging_path}Let $r \in \mathbbm{N}^{d}$,
  $| r | =r_{1} + \ldots .r_{d}$, and $\alpha >- \beta$. Suppose furthermore
  that for $f \in \CC^{\alpha + | r |}$ there exists $\lambda_{0}$ such that
  $K^{w}_{D^{| r |} f, \lambda_{0}} ( \lambda_{0} ) <+ \infty$ (respectively
  there exists $\lambda_{0} >0$ such that $K^{w}_{\beta} ( \lambda ) <+
  \infty$. Then the derivative $\partial^{r} T^{w} f$ is well defined and we
  have $\partial^{r} T^{w} f=T^{w} \partial^{r} f$.
\end{proposition}

\begin{proof}
  First, let us take $f \in \CF L^{\alpha +r}$. For $N \geqslant 0$, the
  projection $\pi_{\leqslant N}$ is a convolution operator, hence
  \[ \partial^{r} T^{w} ( \pi_{\leqslant N} f ) ( x ) =T^{w} \pi_{\leqslant N}
     ( \partial^{r} f ) ( x ) . \]
  But for $f \in \CF L^{\alpha + | r |}$ , $\pi_{\leqslant N} ( \partial^{r} f
  ) \rightarrow^{\CF L^{\alpha}} \partial^{r} f$, which gives the result for
  $f \in \CF L^{\alpha}$. Now take $- \beta < \alpha' < \alpha$. We know that
  for all $f \in \CC^{\alpha +r}$, $\pi_{\leqslant N} \partial^{r} f
  \rightarrow^{\CC^{\alpha'}} \partial^{r} f$, and the result follows for $f
  \in \CC^{\alpha}$.
\end{proof}

\section{Averages along Fractional Brownian paths}

\subsection{Fractional Brownian motion case}

The results of Lemmas~\ref{lemma:chaining-Y} and~\ref{lemma:chaining_besov}
show that in order to control the irregularity constant of fBm paths it is
enough to prove that the there exists $\lambda >0$ and $\alpha \in
\mathbbm{R}$ such that the random variable $e^{\lambda \| B^{H}
\|_{\infty}^{2}} K^{B^{H}}_{\alpha} ( \lambda )$ is almost surely finite when
$B^H$ is a continuous random path with the law of the fBm. Then we only have to
consider the following two quantities:
\[ \exp ( \lambda \| B^H \|_{\infty}^{2} ), \qquad
 \exp ( \lambda ( 1+ | \xi | )^{\alpha} | \Phi_{s,t}^{B^H} ( \xi ) |^{2} / |
   t-s | ) . \]
If the expectation of those quantities are bounded independently of $s,t,x,
\omega , \xi$ then the expectations of $K^{B^{H}}_{f}$ and $e^{\lambda \|
B^{H} \|_{\infty}} K^{B^{H}}_{\alpha} ( \lambda )$ are finite and the variable
are finite almost surely. For $\exp ( \lambda \| B^{H} \|_{\infty}^{2} )$ it
is an application of a well known theorem due to Fernique

\begin{theorem}
  \label{th:fernique}Let $X$ be a Gaussian random variable which takes values in a Banach space $(
  \mathcal{B}, \| . \| )$. Then there exists a constant $\mu >0$ such that
  \[ \mathbbm{E} [ \exp ( \mu \| X \|^{2} ) ] <+ \infty. \]
\end{theorem}

\begin{remark}
  This holds for the fractional Brownian motion of Hurst parameter $H \in (
  0,1 )$ and for the Banach spaces $\left( \CC^{H- \varepsilon} ( [ 0,1 ]
  ,\mathbbm{R}^{d} ) , \| . \|_{0,H- \varepsilon}\right)$ and for $( C ( [ 0,1
  ] \nocomma ,\mathbbm{R}^{d} ) , \| . \|_{\infty , [ 0,1 ]} )$.
\end{remark}

To control the square exponential integrability of $\Phi_{s,t}^{B^{H}} ( \xi
)$ we devised a novel technique based on an elementary application of
Hoeffding inequality for discrete martingale increments. This bypasses the
explicit Gaussian computations or the computations based on Malliavin calculus
usual in the studies involving the fBm. The following theorem then gives
general estimates for additive functionals of the fBm of the form
\[ \int_{s}^{t} f_{u} ( B^{H}_{u} ) \mathd u, \]
where $f: [ 0,1 ] \times \mathbbm{R}^{d} \rightarrow \mathbbm{R}$ is a
measurable and bounded function. Note that the following theorem suggest in
general that such functionals have the same Gaussian deviation behavior of
Brownian martingales.

\begin{theorem}
  \label{th:hoeffding}Let $B^{H} = (B^{H, ( 1 )} , \ldots ,B^{H, ( d )} )$ be a
  $d$-dimensional fractional Brownian motion of Hurst parameter $H \in ( 0,1
  )$ and let $f$ be a function bounded by $1$ and such that
  \[ C_{f} \assign \sup_{u} \int_{0}^{\infty} | P_{t^{2H}} f_{u} |_{\infty}
     \mathd t<+ \infty, \]
  where $P$ is the heat kernel on $\mathbbm{R}^{d}$. Then for $\mu >0$ small
  enough independent of $f$ we have
  \[ \sup_{t \neq s}   \mathbb{E} \left[ \exp   \left( \mu \left\lvert \int_{s}^{t}
     f_{u} ( B^{H}_{u} ) \mathd u \right\lvert^{2} / ( | t-s | C_{f} ) \right)
     \right] <+ \infty . \]
\end{theorem}

\begin{proof}
  The fBm $B^{H}$ can be represented as a stochastic integral over a
  $d$-dimensional standard Brownian motion $W_{\nosymbol} = (W^{(1)} , \cdots
  , \nocomma W^{(d)} )$ defined on the whole $\mathbbm{R}$ (with $W_{0} =0$):
  \[ B^{H, ( i )}_{u} = \int_{- \infty}^{u} ( K ( u,r ) -K ( 0,r ) ) \mathd
     W^{( i )}_{r} , \]
  where $K ( u,r ) = ( u-r )_{+}^{\nosymbol H-1/2} / \Gamma ( H+1/2 )$. Let
  $\left( \CF_{t} \right)_{t \in \mathbbm{R}}$ be the natural filtration of $(
  W_{t} )_{t \in \mathbbm{R}}$. For $v \leqslant u$ we have the decomposition
  \begin{eqnarray*}
    B^{H, ( i )}_{u} & = & \int_{- \infty}^{u} ( K ( u,r ) -K ( 0,r ) ) \mathd
    W^{( i )}_{r}\\
    & = & \int_{v}^{u} K ( u,r ) \mathd W^{( i )}_{r} + \int_{- \infty}^{v} (
    K ( u,r ) -K ( 0,r ) ) \mathd W^{( i )}_{r}\\
    & = & W_{u,v}^{1, ( i )} +W_{u,v}^{2, ( i )},
  \end{eqnarray*}
  where the random variable $W_{u,v}^{1, ( i )}$ is independent of
  $\mathcal{F}_{v}$ and $W_{u,v}^{2, ( i )}$ is $\mathcal{F}_{v}$ measurable.
  We define $W^{j} = (W^{j, ( 1 )} , \ldots ,W^{j, ( d )} )$ for $j \in \{ 1,2
  \}$ and in the following we will note abusively
  $ \tmop{Var} ( W^{j} ) = \tmop{Var} (W^{j, ( 1 )} ) $.
  Now there is two cases we have to consider. Suppose first that $t-s/C_{f}
  \leqslant 1$. Then
  \[ \left\lvert \int_{s}^{t} f_{u} ( B^{H}_{u} ) \mathd u \right\lvert \leqslant | t-s
     |^{2} \leqslant | t-s | C_{f} \]
  and the result follows in that case. Suppose now that $| t-s | C_{f}^{-1}
  \geqslant 1$. Let $N \in \mathbbm{N}$ to be specify later. For $n \in \{ 0,
  \ldots ,N \}$, let us define \ $t_{n} =s+ ( t-s ) n/N$ and
  \begin{eqnarray*}
    Z_{n} & = & \mathbb{E} \left[ \int_{s}^{t} f_{u} ( B^{H}_{u} )  \text{d} u
    | \mathcal{F}_{t_{n+1}} \nobracket \right] - \mathbb{E} \left[
    \int_{s}^{t} f_{u} ( B^{H}_{u} )  \text{d} u | \mathcal{F}_{t_{n}}
    \nobracket \right] .
  \end{eqnarray*}
  Thanks to the previous decomposition of the fractional Brownian motion, we
  are able to bound $Z_{n}$ and to apply Hoeffding Lemma to the sum of the
  martingale increments $(Z_{n} )_{1 \leqslant n \leqslant N}$. Let $S_{N} =
  \sum_{n=0}^{N-1} Z_{n}$, then
  \begin{equation}
    \int_{s}^{t} f_{u} ( B^{H}_{u} ) \mathd u=S_{N} + \mathbb{E} \left[
    \int_{s}^{t} f_{u} ( B^{H}_{u} ) \mathd u \middle\lvert \mathcal{F}_{s} \right]
    . \label{eq:splitting-of-int}
  \end{equation}
  Let us first estimate the conditional expectation in
  Equation~{\eqref{eq:splitting-of-int}} : for all $u \geqslant 0$ we have
  \begin{eqnarray*}
    \left\lvert \mathbb{E} \left[ \int_{s}^{t} f_{u} ( B^{H}_{u} ) \mathd u
    \middle\lvert \mathcal{F}_{s} \right] \right\lvert & = & \left\lvert \int_{s}^{t} P_{Var
    ( W^{1}_{u,s} )} f_{u} ( W^{2}_{u,s} ) \mathd u \right\lvert\\
    & \leqslant & \int_{s}^{t} | P_{\tmop{Var} ( W^{1}_{u,s} )} f_{u}
    |_{\infty} \mathd u \leqslant C_{f} <+ \infty
  \end{eqnarray*}
  since \ $\tmop{Var} ( W^{1}_{u,s} ) =C ( u-s )^{2H}$. But we also have the
  trivial bound
  \[ \left\lvert \mathbb{E} \left[ \int_{s}^{t} f_{u} ( B^{H}_{u} ) \mathd u
     \middle\lvert \mathcal{F}_{s} \right] \right\lvert \leqslant | t-s |. \]
  Hence
  \begin{equation}
    \left\lvert \mathbb{E} \left[ \int_{s}^{t} f_{u} ( B^{H}_{u} ) \mathd u
    \middle\lvert \mathcal{F}_{s} \right] \right\lvert \leqslant | t-s |^{1/2}
    C_{f}^{1/2} . \label{eq:bound-conditional-expect}
  \end{equation}
  Next we bound $Z_{n}$ by decomposing it into three pieces which are easier
  to estimate. We have
  \begin{eqnarray*}
    U_{n} & = & \int^{t}_{t_{n}} \mathbb{E} [ f_{u} ( B^{H}_{u} ) |
    \mathcal{F}_{t_{n}} ] \mathd u= \int^{t}_{t_{n}} \mathbb{E} [ f_{u}  (
    W^{1}_{u,t_{n}} +W_{u,t_{n}}^{2} ) | \mathcal{F}_{t_{n}} ] \mathd u\\
    & = & \int^{t}_{t_{n}} P_{\tmop{Var} (W^{1_{\nosymbol}}_{u,t_{n}} )}
    f_{u} ( W_{u,t_{n}}^{2} ) \mathd u\\
    & = & \int^{t_{n+1}}_{t_{n}} P_{\tmop{Var} (W^{1_{\nosymbol}}_{u,t_{n}}
    )} f_{u} ( W_{u,t_{n}}^{2} ) \mathd u \noplus + \int^{t}_{t_{n+1}}
    P_{\tmop{Var} (W^{1_{\nosymbol}}_{u,t_{n}} )} f_{u} ( W_{u,t_{n}}^{2} )
    \mathd u.
  \end{eqnarray*}
  Hence
  \[ Z_{n} = \int_{t_{n}^{\nosymbol}}^{t_{n+1}} f_{u} ( B^{H}_{u} ) \mathd
     u+U_{n+1} -U_{n} , \]
  moreover
  \[ | U_{n} | \leqslant \int^{t}_{t_{n}} |P_{\tmop{Var}
     (W^{1_{\nosymbol}}_{u,t_{n}} )} f_{u} ( W_{u,t_{n}}^{2} ) | \mathd u
     \leqslant \int_{t_{n}}^{t} |P_{\tmop{Var} ( W^{1_{\nosymbol}}_{u,t_{n}}
     )} f_{u} |_{\infty} \mathd u \leqslant C_{f} <+ \infty \]
  and of course $| \int_{t_{n}^{\nosymbol}}^{t_{n+1}} f_{u} ( B^{H}_{u} )
  \mathd u| \leqslant ( t-s ) /N$, which implies that $| Z_{n} | \lesssim ( t-s
  ) /N+C_{f}$. By the standard Hoeffding inequality we obtain
  \[  \mathbbm{P} ( | S_{N} | > \lambda ) \lesssim \exp ( -2 \lambda^{2} / (
     ( t-s ) N^{-1/2} +N^{1/2} C_{f} )^{2} ). \]
  Hence for $0 \leqslant \nu <1$, we have
  \[ \mathbbm{E} [ \exp ( 2 \nu | S_{N} |^{2} / ( ( t-s ) N^{-1/2} +N^{1/2}
     C_{f} )^{2} ) ] \lesssim \nu / ( 1- \nu ) +1 \]
  Now, we can chose $N= [ 1+ | t-s | /C_{f} ]$, hence
  \[ ( ( t-s ) N^{-1/2} +N^{1/2} C_{f} ) \lesssim | t-s | C_{f} , \]
  and thanks to {\eqref{eq:bound-conditional-expect}} we have
  \[ \mathbbm{E} \left[ \exp \left( \mu \left\lvert \int_{s}^{t} f_{u} ( B^{H}_{u}
     ) \mathd u \right\lvert^{2} / ( | t-s | C_{f} ) \right) \right] \lesssim
     \mathbbm{E} [ \exp ( C \mu | S_{N} |^{2} / ( | t-s | C_{f} ) ) ]
     \lesssim_{\mu} 1.  \]
  
\end{proof}

As an immediate corollary we have the wanted result for the
$\rho$-irregularity constant for the fractional Brownian motion.

\begin{corollary}
  \label{cor:irr-fbm}For $\lambda$ small enough,
  \[ \mathbbm{E}  \exp ( \lambda ( 1+ | \xi | )^{1/H} | \Phi_{s,t}^{B^{H}} (
     \xi ) |^{2} / | t-s | ) \leqslant C<+ \infty \]
  uniformly in $\xi ,t,s$. 
\end{corollary}

\begin{proof}
  When $\xi \leqslant 1$, we have
  \[ \exp ( \lambda ( 1+ | \xi | )^{1/H} | \Phi^{B^{H}}_{s,t} ( \xi ) |^{2} /
     | t-s | ) \leqslant \exp ( \lambda 2^{1/H} ) . \]
  For $| \xi | \geqslant 1$, we have $ ( 1+ | \xi | )^{1/H} \lesssim
  | \xi |^{1/H}$. But we also have
  \begin{eqnarray*}
    | P_{t^{2H}} f ( x ) | & = & | \mathbbm{E} [ \exp ( i \xi ( x+B^{H}_{t} )
    ) ] |
    =  \exp ( - | \xi |^{2} t^{2H} /2 ) ,
  \end{eqnarray*}
  therefore
  $ |\xi|^{-1/H}\lesssim C_{f_{\xi}} \lesssim | \xi |^{-1/H} $. 
  Finally there exists a constant $C>0$ such that
  \begin{eqnarray*}
    \mathbbm{E}  \exp ( \lambda ( 1+ | \xi | )^{1/H} | \Phi_{s,t}^{B^{H}} (
    \xi ) |^{2} / | t-s | ) & \leqslant & \mathbbm{E}  \exp ( C \lambda |
    \Phi_{s,t}^{B^{H}} ( \xi ) |^{2} / ( | t-s | C_{f_{\xi}} ) )
  \end{eqnarray*}
  and for $\lambda$ small enough, thanks to Theorem \ref{th:hoeffding} the right hand side
  is bounded by a constant independent of $\xi ,s,t$.
\end{proof}

We are now in condition to prove Theorem~\ref{th:irr-for-fbm}
($\rho$-irregularity of the fBm paths for all $\rho <1/2H$).

\begin{proof}
  (of \ Theorem~\ref{th:irr-for-fbm}) By Lemma~\ref{lemma:chaining-Y} we have
  \[ \| \Phi^{B^{H}} \|_{\mathcal{W}^{\rho , \gamma}_{1}} = \sup_{\xi \in
     \mathbbm{R}^{d}}   \sup_{0 \leqslant s<t \leqslant 1} ( 1+ | \xi |
     )^{\rho} \frac{| \Phi^{B^{H}}_{t} ( \xi )- \Phi^{B^{H}}_{s} ( \xi
     )|}{|s-t|^{\gamma}} \lesssim 1+ \log^{1/2} ( e^{\lambda \| B^{H}
     \|_{\infty}} K^{B^{H}}_{1/2H} ( \lambda ) ) . \]
  Moreover by Theorem~\ref{th:fernique} the quantity $e^{\lambda \| B^{H}
  \|_{\infty}^{2}}$ is almost surely finite and by Corollary~\ref{cor:irr-fbm}
  we readily have that
  \[ \mathbbm{E} [K^{B^{H}}_{1/2H} ( \lambda ) ] <+ \infty \]
  as soon as $\lambda$ is small enough.
\end{proof}

\begin{remark}
  A byproduct of the proof of Theorem~\ref{th:irr-for-fbm} is that the
  irregularity constant $\| \Phi^{B^{H}} \|_{\mathcal{W}^{\rho , \gamma}_{1}}$
  is exponentially square integrable, as easily shown: for small $\lambda >0$
  and all $\gamma >1/2$ and $\rho <1/2H$ we have
  \[ \mathbbm{E} \left[ e^{\lambda \| \Phi^{W} \|^{2}_{\mathcal{W}^{\rho ,
     \gamma}_{1}}} \right] <+ \infty . \]
\end{remark}

When we consider the spaces $\CC^{\alpha}$ instead of $\CF L^{\alpha}$,the
$\rho$-irregularity of the fractional Brownian path is not enough.
Nevertheless Theorem \ref{th:hoeffding} allows us to give the correct bound
for $T^{B^{H}} \Delta_{i} f$ in order to define $T^{B^{H}} f$ for any $f \in
\CC^{\alpha}$.

\begin{corollary}
  Let $B^{H}$ be a $d$-dimensional fractional Brownian motion of Hurst parameter
  $H \in ( 0,1 )$. There exists $\lambda >0$ such that for $f \in C ( [ 0,T ]
  ;\mathcal{S}' ( \mathbbm{R}^{d} ) )$, all $i \geqslant -1$
  \[ \sup_{x \in \mathbbm{R}^{d} ,i \geqslant -1,s \neq t} \mathbbm{E} [ \exp
     ( \lambda 2^{i/H} | T^{w} ( \Delta_{i} f ) ( x ) |^{2} / ( | t-s | \|
     \Delta_{i} f \|_{\infty}^{2} ) ) ] <+ \infty. \]
\end{corollary}

\begin{proof}
  The function $g_{u} = \Delta_{i} f_{u} / \| \Delta_{i} f \|_{\infty}$ is
  bounded by $1$. Furthermore, for $i \geqslant 0$, $\tmop{supp}   \hat{g}_{u}
  \subset 2^{i} \mathcal{A}$ where $\mathcal{A}$ is an annulus. Hence by the
  Lemma 2.4 page 54 of {\cite{bahouri_fourier_2011}}, there exists a constant
  $c>0$ independent of $g$ such that
  \[ \| P_{t^{2H}} g_{u} \|_{\infty} \lesssim \exp ( -c t^{2H} 2^{2i} ), \]
  hence $C_{g} \lesssim 2^{-i/H}$ and the result follows immediately by
  applying Theorem \ref{th:hoeffding} to $g$.
  
  When $i=-1$, we have
  \[ \left\lvert \int_{s}^{t} \Delta_{-1} f ( B^{H}_{u} ) \mathd u \right\lvert^{2}
     \leqslant 2^{1/H} | t-s | \| \Delta_{-1} f \|_{\infty}  ^{2} \]
  and the result follows.
\end{proof}

The result is exactly the needed hypothesis of Theorem
\ref{theorem:definition_T}, but also Propositions
\ref{prop:holder_continuity_of_averaging_path} and
\ref{prop:diffenrentiability_averaging_path}, depending on the regularity of
$f$. Hence, the averaging operator $T^{B^{H}}$, or its finite dimensional
marginals depending of the space, is well defined and has the right range for
applying results of the Section~\ref{sec:young-ode}. This operator is
defined almost surely, hence, almost surely (depending of $b$ when $b \in
\CC^{\alpha}$) the following equation has a solution
\[ \theta_{u} = \theta_{0} + \int_{0}^{t} T^{B^H}_{\mathd u} b ( \theta_{u} ), \]
when $b \in \CC^{\alpha}$. The existence of such a solution is guaranteed by Theorem~\ref{theorem:existence_young_solution}. Furthermore, when $b \in \CC^{\alpha +2}$,
$\alpha >-1/2H$ (or $\CC^{\alpha + \gamma +1}$) there is uniqueness and the
flow is Lipschitz-continuous, thanks to the
Corollary~\ref{corollary:uniqueness_a_priori_case}. As the set where $T^{w} b$
is not defined does not depend on $b$ when $b \in \CF L^{\alpha}$, those
results does not depend on $b$. The uniqueness for $b \in \CC^{\alpha +1}$ or
$b \in \CF L^{\alpha +1}$ is a more probabilistic argument and is the subject
of the following section.

\begin{remark}
  If the regularity of $b$ is not enough to guarantee uniqueness by the above
  arguments the solution constructed via
  Theorem~\ref{theorem:existence_young_solution} lacks, a priori,
  measurability with respect to $B^H$. If a measurable solution is needed the
  fix--point argument of \ Theorem~\ref{theorem:existence_young_solution} has
  to be repeated in a space of random processes, for example in $L^{p} (
  \Omega ,\CC^{\gamma} ( [ 0,1 ] ;\mathbbm{R}^{d} ) )$.
\end{remark}

\subsection{Averaging for absolutely continuous perturbations of the fBm}\label{sec:girsanov}

In this section, we analyse the properties of the averaging operator along a
path of the form $B^{H} + \theta^{n}$ where $\theta^{n}$ is a solution to the
approximate equation $\mathd \theta^{n} =T^{B^{H}}_{\mathd u} b^{n} (
\theta^{n} )$. In order to do so we will use a version of the Girsanov theorem
for fractional Brownian motion. The results holds of course for  both
types of functions spaces $\CC^{\alpha +1}$ and $\CF L^{\alpha +1}$. We will
only give the proof for $b \in \CC^{\alpha +1}$, and give some comments for
the case $\CF L^{\alpha +1}$. Let $(b^{n} )_{n \ge 1}$ be a sequence of smooth
vector fields such that $\| b^{n} \|_{\alpha} \leqslant C$ uniformly in $n \ge
1$. By a standard fixed point argument, it is well known that the following
equation
\begin{equation}
  X^{n}_{t} =x_{0} + \int_{0}^{t} b^{n}_{s} (X^{n}_{s} ) ds+B^{H}_{t}
\end{equation}
has an adapted solution $X^{n}$ (to the standard filtration of the fractional
Brownian motion).

Here we analyse the averaging constant of $X_{n}$ and we prove that it
satisfy the requirements of Theorem \ref{theorem:uniqueness_general_case}
implying uniqueness of the limit ODE for $b \in \CC^{\alpha +1}$ and
convergence of $X^{n}$ to this unique solution. Furthermore, if we consider,
as in Section~\ref{subsec:comparison_principle}, the averaged translation by
$\theta^{n} =X^{n} -B^{H}$, we only have to check the hypothesis of
Theorem~\ref{theorem:uniqueness_general_case}. Indeed, if $\theta$ is the solution
of
\[ \theta_{t} = \theta_{0} + \int_{0}^{t} T^{B^{H}}_{\mathd u} b ( \theta_{u}
   ) \]
and $\theta^{n}$ is the solution of
\[ \theta_{t}^{n} = \theta_{0} + \int_{0}^{t} T^{B^{H}}_{\mathd u} b^{n} (
   \theta^{n}_{u} ), \]
then $\theta^{n} +B^{H}_{u} =X^{n}$ and by the averaged translation by
$\theta^{n}$, as $\tau_{\theta^{n}} T^{B^{H}} b=T^{X_{n}} b$. Then $\theta -
\theta^{n}$ is the solution of the following Young equation
\[ ( \theta - \theta^{n} )_{t} = \int_{0}^{t} T^{X^{n}}_{\mathd u} b (
   \theta_{u} - \theta^{n}_{u} ) + \int_{0}^{t} T^{B^{H}}_{\mathd u} b^{n} (
   \theta^{n}_{u} ). \]
These considerations are the motivation to introduce the comparison principle
based on averaged translations in the proof of uniqueness in Section
\ref{subsec:uniqueness}.

Below we will take advantage of the absolute continuity of the law of $X^{n}$
w.r.t. the law of the fractional Brownian motion $B^{H}$ to transfer the averaging properties of the
fractional Brownian motion to the stochastic process $X^{n}$. This approach is an extension of an
observation of Davie~{\cite{MR2377011}} to the fractional Brownian motion's context.

A drawback of this approach is that the exceptional set will necessarily
depend on the initial point $x_{0}$ and on the vector field $b$. This prevents us
from easily applying the uniqueness result to the case of random $b$ and to the
analysis of the flow of the ODE.

The computation of the Radon-Nikodym derivative between the law of $X^{n}$ and
the law of $B^{H}$ will result in a Girsanov transform. For technical reasons
we will do this transformation only on a subinterval $[0,T_{Gir} ] \subset
[0,1]$. For $b^{n}$ regular enough, and as $X^{n}$ is regular enough,
according to Nualart and Ouknine~{\cite{MR1934157}}, there exist a Brownian
motion $W$ adapted to the filtration associated with $B^{H}$ and a probability
$\mathbbm{P}_{n}$ such that the process $(X^{n}_{t} )_{t \in [0,T_{Gir} ]}$ is
a fractional Brownian motion of Hurst parameter $H$, where
\[ \frac{\mathd \mathbbm{P}_{n}}{\mathd \mathbbm{P}} = \exp   \left( -
   \int_{0}^{T_{Gir}} H_{t}^{n} \cdot \mathd W_{t} - \frac{1}{2} 
   \int_{0}^{T_{Gir}} | H^{n}_{t} |^{2} \mathd t \right), \]
 where for $H \ge \frac{1}{2}$
\[ H_{t}^{n} = \frac{t^{H- \frac{1}{2}}}{\Gamma ( \frac{3}{2} -H)}  \left(
   t^{1-2H} b^{n} (t,X^{n}_{t} )+ \left( H- \frac{1}{2} \right)  \int_{0}^{t}
   \frac{t^{\frac{1}{2} -H} b^{n}_{t} (X^{n}_{t} )-s^{\frac{1}{2} -H}
   b^{n}_{s} (X^{n}_{s} )}{(t-s)^{H+ \frac{1}{2}}} \mathd s \right) \]
and for $H< \frac{1}{2}$
\[ H_{t}^{n} = \frac{t^{H- \frac{1}{2}}}{\Gamma ( \frac{1}{2} -H)} 
   \int_{0}^{t} ( s(t-s) )^{\frac{1}{2} -H} b^{n}_{s} (X^{n}_{s} ) \mathd s. \]
Thanks to that Girsanov transform, the almost sure bound for $T^{B^{H}}_{b}$
can be used to estimate $T^{X_{n}}_{b}$ since $\mathbbm{P}_{n}$ and
$\mathbbm{P}$ are equivalent.

\begin{lemma}
  \label{majoration Qn}Let $0 \geqslant \alpha >-1/2H$. There exists a
  constant $\lambda >0$ and a constant $C_{\lambda}$ independent of $n$ such
  that for all $b \in \CC^{H- \varepsilon} \left( [ 0,1 ] ; \CC^{\alpha +1}
  \right) \cap \CC^{\alpha +1} \left( \mathbbm{R}^{d} , \CC^{H-1/2+
  \varepsilon} ( [ 0,T ] ) \right)$,
  \[ \mathbbm{E} [ K^{X_{n}}_{b,1/2H} ( \lambda ) ] \leqslant C_{\lambda} .\]
\end{lemma}

Until the end of the section, we will only consider $K^{X_{n}}_{b,1/2H}$. For
simplicity we only write it as $K^{X_{n}}_{b}$.

\begin{proof}
  Let $K>0$. By using the notation above, we have
  \begin{eqnarray*}
    \mathbbm{E} [ K^{X_{n}}_{b} ( \lambda ) ]^{2} & = &
    \mathbbm{E}_{\mathbbm{P}_{n}} \left[ K^{X_{n}}_{b} ( \lambda )
    \frac{\mathd \mathbbm{P}}{\mathd \mathbbm{P}_{n}} \right]^{2}\\
    & \leqslant & \mathbbm{E}_{\mathbbm{P}_{n}} [ K^{X_{n}}_{b} ( \lambda
    )^{2} ] \mathbbm{E}_{\mathbbm{P}_{n}} \left[ \left( \frac{\mathd
    \mathbbm{P}}{\mathd \mathbbm{P}_{n}} \right)^{2} \right]\\
    & \lesssim & \mathbbm{E} [ K^{B^{H}}_{b} ( 2 \lambda ) ]
    \mathbbm{E}_{\mathbbm{P}_{n}} \left[ \exp   \left( 2 \int_{0}^{T}
    H^{n}_{t} \mathd W_{t} + \int_{0}^{T} H^{n}_{t} \mathd t \right) \right].
  \end{eqnarray*}

  Where we have used that under $\mathbbm{P}_{n}$, $X^{n}$ is a fractional
  Brownian motion of same Hurst parameter $H$. If $\rho$ is small enough the
  first term is finite by the above results. To prove the lemma, it is
  sufficient to prove that
  \[ \mathbbm{E}_{\mathbbm{P}_{n}} \left[ \exp   \left( 2 \int_{0}^{T}
     H^{n}_{t} \mathd W_{t} + \int_{0}^{T} | H^{n}_{t} |^{2} \mathd t \right)
     \right] \]
  is bounded by a constant independent of $n$. As $W$ is a Brownian motion, it
  is enough to bound
  \[ \mathbbm{E}_{\mathbbm{P}_{n}} \left[ \exp   \left( \int_{0}^{T} |
     H_{t}^{n} |^{2} \mathd t \right) \right]. \]
  The arguments are quite different depending whether $H>1/2$ or $H<1/2$.
  First suppose that $H< \frac{1}{2}$.
  \begin{eqnarray*}
    |H_{t}^{n} |^{2} & = & \left\lvert K_{H}^{-1} \left( \int_{0}^{.} b^{n}_{s}
    (X^{n}_{s} ) \mathd s \right) (t) \right\lvert^{2}\\
    & = & \left\lvert \frac{t^{H- \frac{1}{2}}}{\Gamma ( \frac{1}{2} -H)} 
    \int_{0}^{t} (s(t-s))^{\frac{1}{2} -H} b^{n}_{s} (X^{n}_{s} ) \mathd s
    \right\lvert^{2}\\
    & = & \left\lvert \frac{- \left( \frac{1}{2} -H \right) t^{H-
    \frac{1}{2}}}{\Gamma ( \frac{1}{2} -H)}  \int_{0}^{t} (s(t-s))^{-(
    \frac{1}{2} +H)} (t-2s) \underbrace{\int_{0}^{s} b^{n}_{u} (X^{n}_{u} )
    \mathd u}_{( T_{s}^{X^{n}} b^{n} ) (0)} \mathd s \right\lvert^{2}\\
    & \lesssim & t^{2H}  \int_{0}^{t} (s(t-s))^{- (1+2H)}  | t-2s | | (
    T_{s}^{X^{n}} b^{n} ) (0) |^{2} \mathd s\\
    & \lesssim & \|b^{n} \|_{\CC^{\alpha}}^{2} ( 1+ \log^{1/2} (
    K^{X_{n}}_{b} ( \lambda ) ) )^{2} t^{2H}  \int_{0}^{t} (s(t-s))^{- (1+2H)}
    | t-2s | s^{2 \gamma} \mathd s\\
    & \lesssim & \|b^{n} \|_{\CC^{\alpha}}^{2} ( 1+ \log^{} ( K^{X_{n}}_{b} (
    \lambda ) ) )^{} t^{2 ( \gamma -H)}  \int_{0}^{1} (u(1-u))^{- (1+2H)}  |
    1-2u | u^{2 \gamma} \mathd s\\
    & \leqslant & C (b,H, \gamma , \lambda ) ( 1+ \log^{} ( K^{X_{n}}_{b} (
    \lambda ) ) ).
  \end{eqnarray*}
  Hence,
  \begin{eqnarray*}
    \mathbbm{E}_{\mathbbm{P}_{n}} \left[ \exp   \left( \int_{0}^{T} |
    H_{t}^{n} |^{2} \mathd t \right) \right] & \lesssim &
    \mathbbm{E}_{\mathbbm{P}_{n}} [ K^{X_{n}}_{b} ( C \lambda ) ]\\
    & \lesssim & \mathbbm{E} [  K^{B^{H}}_{b} ( C \lambda ) ].
  \end{eqnarray*}
  For $\lambda$ small enough, this quantity is bounded, and the Lemma is
  proved in this case.
  
  For $H \geqslant \frac{1}{2}$, $b^{n}$ is $( 1+ \alpha)$-H\"older
  continuous and $\|b^{n} \|_{\CC^{\alpha +1}} \lesssim \|b\|_{\CC^{\alpha
  +1}}$. Furthermore, $|(b^{n} (0,0))_{n} |$ is uniformly bounded, then
  \begin{eqnarray*}
    | H_{t}^{n} | & \lesssim & \left\lvert t^{H-1/2}  \left( t^{1-2H} b^{n}
    (t,X^{n}_{t} )+ \int_{0}^{t} \frac{t^{1/2-H} b^{n}_{t} (X^{n}_{t}
    )-s^{1/2-H} b^{n}_{s} (X^{n}_{s} )}{(t-s)^{H+1/2}} \mathd s \right)
    \right\lvert\\
    & \lesssim & t^{1/2-H} ( \|b^{n} \|_{\alpha +1} + \| b^{n}_{.} ( 0 )
    \|_{\infty} ) (|X^{n}_{t} -X^{n}_{0} |^{1+ \alpha} +1)\\
    &  & +t^{H-1/2}  \int_{0}^{t} (t-s)^{- (H+1/2)}  | ( t^{1/2-H} -s^{1/2-H}
    )  ( b^{n}_{t} (X_{t}^{n} )+b^{n}_{s} (X_{s}^{n} ) ) | \mathd s\\
    &  & +t^{H-1/2}  \int_{0}^{t} | (t-s)^{-(H+1/2)}  ( t^{1/2-H} +s^{1/2-H}
    )  ( b^{n}_{t} (X_{t}^{n} )-b^{n}_{s} (X_{s}^{n} ) ) | \mathd s .
  \end{eqnarray*}
  The first term is bounded by $( \| b^{n} \|_{\alpha} + \| b^{n}_{.} ( 0 )
  \|_{\infty} ) ( \| X^{n} \|^{1+ \alpha}_{H- \varepsilon} +1) t^{1/2-H}$ and
  is integrable. The second is bounded by
  \begin{equation*}
  \begin{split}
   & t^{H-1/2}  \int_{0}^{t} \mathd s (t-s)^{- (H+1/2)} | t^{1/2-H}
     -s^{1/2-H} | \\ & \qquad \times ( | b^{n}_{t} (X_{t}^{n} )-b^{n} (t,0) | + | b^{n}_{s}
     (X_{s}^{n} )-b^{n}_{t} (0) | +2 \| b_{.} ( 0 ) \|_{\infty} )  
  \\ & \lesssim t^{H-1/2+1-H-1/2+1/2-H} \int_{0}^{1} ( 1-u )^{-H-1/2} (
     1-u^{1/2-H} ) ( \| b^{n} \|_{\alpha} + \| b^{n}_{.} ( 0 ) \|_{\infty} )
     \| X^{n} \|_{\infty}^{\alpha +1} 
  \\ & \lesssim_{x_{0}} t^{1/2-H} ( ( \| b^{n} \|_{\alpha} + \| b^{n}_{.} ( 0 )
     \|_{\infty} ) ( \| X^{n} \|_{H- \varepsilon}^{\alpha +1} +1 ) + \| b_{.}
     ( X^{n}_{s} ) \|_{H-1/2+ \varepsilon} | t-s |^{H-1/2+ \varepsilon} ). \end{split}
     \end{equation*}

  The third term is bounded by
  \begin{eqnarray*}
    &  & t^{H- \frac{1}{2}}  \int_{0}^{t} \left\lvert (t-s)^{-(H+ \frac{1}{2} )}
    s^{\frac{1}{2} -H}  ( | b^{n}_{t} (X_{t}^{n} )-b^{n}_{t} ( X_{s}^{n} ) | +
    | b^{n}_{t} ( X_{s}^{n} ) -b^{n}_{s} (X_{s}^{n} ) | ) \right\lvert \mathd s\\
    & \lesssim & ( \| b^{n} \|_{\alpha} + \| b ( .,0 ) \|_{H-1/2+
    \varepsilon} ) ( \|X^{n} \|^{1+ \alpha}_{H- \varepsilon} +1 ) t^{H-
    \frac{1}{2}}   \\ & & \times \int_{0}^{t} (t-s)^{- (H+ \frac{1}{2} )} s^{\frac{1}{2} -H} 
    \left( |t-s|^{H- \frac{1}{2} + \varepsilon} +|t-s|^{( \alpha +1)(H-
    \varepsilon )} \right) \mathd s\\
    & \lesssim & ( \| b^{n} \|_{\alpha} + \| b^{n}_{.} ( 0 ) \|_{H-1/2+
    \varepsilon} + \| b^{n}_{.} ( 0 ) \|_{\infty} ) ( \|X^{n} \|^{1+
    \alpha}_{H- \varepsilon} +1 ) t^{\varepsilon} \int_{0}^{1} ( 1-u )^{-1+
    \varepsilon} u^{-H+1/2} \mathd u .
  \end{eqnarray*}
  Finally, we choose $b^{n}$ such that $( \| b^{n} \|_{\alpha} + \| b^{n}_{.}
  ( 0 ) \|_{H-1/2+ \varepsilon} + \| b^{n}_{.} ( 0 ) \|_{\infty} )
  \lesssim_{b} 1$, hence
  \[ | H_{t}^{n} | \lesssim C_{b,x_{0}} ( \|X^{n} \|^{1+ \alpha}_{H-
     \varepsilon} +1 ) t^{\varepsilon}. \]
  Under $\mathbbm{P}_{n}$, $X^{n}$ is a fractional Brownian motion of
  Hurst parameter $H$. Thanks to Fernique theorem, $\|X^{n} \|^{1+ \alpha}_{H-
  \varepsilon}$ is exponentially integrable in $\mathbbm{P}_{n}$ and the
  result follows.
\end{proof}

This bound in $L^{1}$ is not enough to use the Theorem
\ref{theorem:uniqueness_general_case}, as we need an almost surely, uniformly
in $n$, bound for $\| T^{X^{n}} b \|_{\CC^{\gamma ,1, \psi}}$. Nevertheless, by
using the results of Section~\ref{sec:bounds-on-Tw}, we
already know that
\[ \exp \left( C \| T^{X^{n}} b \|_{\CC^{\gamma ,1, \psi}} \right) \lesssim
   1+K^{X^{n}}_{b} ( \lambda ). \]
We have all the tools to prove the following theorem

\begin{theorem}
  Assume that $b \in \CC^{\alpha +1}$. There \ then there exists $\lambda >0$
  and sequence of smooth vector fields $(b^{n} )_{n}$ such that $b^{n} \to b$
  in $\CC^{\alpha'}$ for all $\alpha' < \alpha$ and almost surely
  \[ K_{b}^{X_{n}} ( \lambda )  \|b-b^{n} \|_{^{\alpha' +1}} \to 0, \]
  which implies uniqueness of the Young equation for $b$ by Theorem
  \ref{theorem:uniqueness_general_case}.
\end{theorem}

\begin{proof}
  By the previous result we have that the $L^{1}$ norm of $K_{b}^{X_{n}} (
  \lambda )$ is uniformly bounded in $n$. Moreover consider $b^{n}$ such that
  $b^{n} = \rho_{n} \ast b$ where $\rho_{n} (x) = \frac{1}{2} n^{d} \exp  (-n
  |x| )$. Then $b^{n}$ is smooth, $b^{n} \to b$ in $\CC^{1+ \alpha'}$ for all
  $\alpha' < \alpha$ by the dominated convergence theorem and there exists a
  subsequence which will still denote with $b^{n}$ such that $\| b-b^{n}
  \|_{\alpha' +1} \lesssim n^{-2}$. On this subsequence (which depends on $b$)
  consider the random variable
  \[ D= \sum_{n \ge 1} K_{b}^{X_{n}} ( \lambda ) \|b-b^{n} \|_{\alpha' +1}. \]
  Then
  \[ \mathbbm{E}D= \sum_{n \ge 1} \mathbbm{E} [ K_{b}^{X_{n}} ( \lambda ) ] 
     \|b-b^{n} \|_{_{\alpha' +1}} \lesssim \sum_{n \ge 1} n^{-2} \lesssim 1, \]
  so that almost surely $D< \infty$ which implies that $K_{b}^{X_{n}} (
  \lambda )  \|b-b^{n} \|_{\alpha' +1} \to 0$.
\end{proof}

Note that this argument give an exceptional set of zero measure which a priori
depends on $b$ (and on the sequence $(b^{n} )_{n}$) and of $x_{0}$. As
remarked previously, this fact prevents straightforward extension of the
uniqueness results in $\CC^{\alpha}$ to random $b$. Furthermore, it also
prevent to consider the regularity of the flow of the equation by pathwise
methods.

\end{document}